\documentclass[12pt]{article}
\usepackage{pictex}
\usepackage{amsmath}
\usepackage{amscd}
\usepackage{amstext}
\usepackage{amssymb}
\usepackage{amsfonts}
\usepackage{latexsym}

\usepackage[dvips]{color}

\title{On the Interplay of Regularity and Decay in Case of Radial Functions I. Inhomogeneous spaces}
\author{Winfried Sickel\thanks{Corresponding author}, Leszek Skrzypczak and Jan Vybiral}

\date{\today}
\begin{document}

\maketitle

\begin{abstract}
We deal with decay and boundedness  properties of radial functions belonging to
Besov and Lizorkin-Triebel spaces. In detail we investigate the
surprising interplay of regularity and decay. Our tools are atomic decompositions
in combination with trace theorems.
\end{abstract}

\newcommand{\open}[1]{\smallskip\noindent\fbox{\parbox{\textwidth}{\color{blue}\bfseries\begin{center}
      #1 \end{center}}}\\ \smallskip}


\newtheorem{T}{Theorem}
\newtheorem{Lem}{Lemma}
\newtheorem{Prop}{Proposition}
\newtheorem{Def}{Definition}
\newtheorem{Cor}{Corollary}
\newtheorem{Rem}{Remark}

\newcommand{\bt}{\begin{T}}
\newcommand{\bl}{\begin{Lem}}
\newcommand{\bp}{\begin{Prop}}
\newcommand{\bc}{\begin{Cor}}
\newcommand{\bd}{\begin{Def}}
\newcommand{\br}[2]{\begin{Rem}\label{#1}{\rm #2}}
\newcommand{\er}{ \end{Rem}}
\newcommand{\et}{\end{T}}
\newcommand{\el}{\end{Lem}}
\newcommand{\ep}{\end{Prop}}
\newcommand{\ec}{\end{Cor}}
\newcommand{\ed}{\end{Def}}


\newcommand{\be}{\begin{equation}}
\newcommand{\ee}{\end{equation}}
\newcommand{\beq}{\begin{eqnarray}}
\newcommand{\eeq}{\end{eqnarray}}
\newcommand{\beqq}{\begin{eqnarray*}}
\newcommand{\eeqq}{\end{eqnarray*}}


\def\supp{\mathop{\rm supp \,}\nolimits}
\def\vol{\mathop{\rm vol}\nolimits}
\def\div{\mathop{\rm div}\nolimits}
\def\dist{\mathop{\rm dist \,}\nolimits}
\def\diam{\mathop{\rm diam \,}\nolimits}
\def\Tr{\mathop{\rm Tr \,}\nolimits}
\def\id{\mathop{\rm id \,}\nolimits}
\def\tr{\mathop{\rm tr \,}\nolimits}
\def\ext{\mathop{\rm ext \,}\nolimits}


\newcommand{\QED}{$\Box$}
\newcommand{\reff}[1]{{\rm (\ref{#1})}}
\newcommand{\bpr}{{\bf Proof.}\ }
\newcommand{\epr}{\hspace*{\fill}\rule{3mm}{3mm}\\}

\newcommand{\cn}{{\mathcal N}}
\newcommand{\bspq}{B^s_{p,q}}
\newcommand{\tbspq}{TB^s_{p,q}}
\newcommand{\tfspq}{TF^s_{p,q}}
\newcommand{\fspqh}{\dot{F}^s_{p,q}}
\newcommand{\bspqh}{\dot{B}^s_{p,q}}
\newcommand{\fspq}{F^s_{p,q}}
\newcommand{\rbspq}{RB^s_{p,q}}
\newcommand{\rfspq}{RF^s_{p,q}}
\newcommand{\Rd}{{\mathbb R}^{d}}
\newcommand{\R}{{\mathbb R}}
\newcommand{\N}{{\mathbb N}}
\newcommand{\Z}{\mathbb Z}
\newcommand{\Zd}{\Z^{d}}
\newcommand{\C}{{\mathbb C}}
\newcommand{\cs}{{\mathcal{S}}}
\newcommand{\cz}{{\mathcal{Z}}}
\newcommand{\ca}{{\mathcal{A}}}
\newcommand{\cx}{{\mathcal{X}}}
\newcommand{\cl}{{\mathcal{L}}}
\newcommand{\wt}{\widetilde}
\newcommand{\aspq}{A^s_{p,q}}
\newcommand{\taspq}{TA^s_{p,q}}
\newcommand{\aspqh}{\dot{A}^s_{p,q}}


\newcommand{\cf}{{\mathcal{F}}\,}
\newcommand{\cfi}{{\mathcal{F}}^{-1}\,}

\def\lsim{\raisebox{-1ex}{$~\stackrel{\textstyle <}{\sim}~$}}

\def\gsim{\raisebox{-1ex}{$~\stackrel{\textstyle >}{\sim}~$}}


\tableofcontents


\section{Introduction}


At the end of the seventies Strauss \cite{St}
was the first who  observed that
there is an interplay between the regularity and  decay properties of radial functions.
We recall his
\\
 {\em Radial Lemma:   Let $d\ge 2$.
Every radial function $f \in H^1 (\Rd)$ is almost everywhere equal to a function
$\wt{f}$, continuous for $x\neq 0$, such that
\be\label{strauss}
|\wt{f} (x)| \le c \, |x|^{\frac{1-d}{2}}\, \| \,f \, |H^1(\Rd)\|\, ,
\ee
where $c$ depends only on $d$.}
\\
Strauss stated (\ref{strauss}) with the extra condition $|x|\ge 1$, but this restriction is not needed.
The {\em Radial Lemma} contains three different assertions:
\begin{itemize}
 \item[(a)]
 the existence of a representative of $f$, which is
continuous outside the origin;
\item[(b)]
the decay of $f$ near infinity;
\item[(c)]
the limited unboundedness near the origin.
\end{itemize}
These three  properties do not extend to all functions in $H^1(\Rd)$, of course.
In particular, $H^1 (\Rd) \not \subset L_\infty (\Rd)$, $d\ge 2$, and consequently,
functions in $H^1 (\Rd)$ can be unbounded in the neigborhood of any fixed point $x\in \Rd$.
It will be our aim in this paper to investigate the specific regularity
and decay properties  of radial functions in a more  general framework
than Sobolev spaces.
In our opinion a discussion of these properties in connection with
fractional order of smoothness results in a better understanding of
the announced interplay of regularity on the one side and
local smoothness, decay at infinity and limited unboundedness near the origin
on the other side.
In the literature there are several approaches to fractional order of smoothness.
Probably most popular are Bessel potential spaces $H^s_p (\Rd)$, $s\in \R$, or
Slobodeckij spaces $W^s_p (\Rd)$ ($s>0$, $s\not\in\N$).
These scales would be enough to explain the main interrelations.
However, for some limiting cases these scales are not sufficient.
For that reason we shall discuss generalizations of the
{\em Radial Lemma} in the framework of
Besov spaces $\bspq(\Rd)$ and Lizorkin-Triebel spaces $\fspq (\Rd)$.
These scales essentially cover the Bessel potential and  the Slobodeckij spaces since
\begin{itemize}
 \item $W^m_p (\Rd) = F^m_{p,2} (\Rd)$, $m \in \N_0$,  $1<p<\infty$;
\item $H^s_p (\Rd) = F^s_{p,2} (\Rd)$, $s \in \R$, $1<p<\infty$;
\item $W^s_p (\Rd) = F^s_{p,p} (\Rd) = B^s_{p,p} (\Rd)$, $s>0$, $s \not\in \N$,
$1\le p\le \infty$,
\end{itemize}
where all identities have to be understood in the sense of
equivalent norms, see, e.g., \cite[2.2.2]{Tr1} and the references given there.
\\
All three phenomena (a)-(c) extend to a certain range of parameters
which we shall characterize exactly. For instance, decay near infinity will take place in spaces with $s\ge 1/p$ (see Theorem \ref{decay4})
and limited unboundedness near the origin
in the sense of (\ref{strauss}) will happen in spaces
such that $1/p \le  s \le d/p$ (see Theorem \ref{decay2}).
For $s=1/p$ (or $s=d/p$) always the microscopic
parameter $q$ comes into play.
We will study the above properties also for spaces
with $p<1$. To a certain extent this is motivated  by the fact, that the decay properties of radial functions
near infinity are determined by the parameter $p$ and the decay rate increases when $p$ decreases, see Theorem \ref{decay4}.
\\
Our main tools here are the following. Based on the atomic decomposition theorem for  inhomogeneous  Besov and Lizorkin-Triebel spaces,
which we  proved in \cite{SS}, we shall deduce  a trace theorem for radial subspaces which is of interest on its own.
Then this trace theorem will be applied to derive the extra regularity properties of radial functions.
To derive the decay estimates and the assertions on controlled
unboundedness near zero  we shall also employ
the atomic decomposition technique. \\
With respect to the decay it makes a difference, whether one deals with inhomogeneous
or homogeneous spaces of Besov and Lizorkin-Triebel type.
Homogeneous spaces (with a proper interpretation)
are larger than their inhomogeneous counterparts
(at least if $s> d \max (0,\frac 1p - 1)$).
Hence, the decay rate of the elements of
inhomogeneous spaces can be better than that one for
homogeneous spaces. This turns out to be true.
However, here in this article we concentrate on inhomogeneous spaces.
Radial subspaces of homogeneous spaces will be subject to the continuation of this paper, see \cite{SSV2}.
In a further paper \cite{SSV3} we shall investigate a few more properties of radial subspaces like complex interpolation
and characterization by differences.
\\
The paper is organized as follows.
In Section 2 we describe our main results.
Here Subsection 2.1 is devoted to the study of traces 
of radial subspaces. In particular, in 2.1.1, we state also trace assertions for radial
subspaces of H\"older-Zygmund classes. 
In Subsection 2.2 the behaviour radial functions near infinity
and at the origin is investigated.
Within the borderline cases in
Subsection 2.2.2 the spaces $BV(\Rd)$ show up. In this context we will also
deal with the  trace problem for the associated radial subspaces.
All proofs will be given in Sections \ref{trace} and \ref{inhom}.
There also additional material is collected, e.g., in Subsection \ref{int}
we deal with interpolation of radial subspaces, in Subsection \ref{rsa}
we recall the characterization of radial subspaces by atoms as given in \cite{SS},
and finally,
in Subsection \ref{testfunctions} we discuss the regularity properties of some
families of test functions.
\\
Besov and Lizorkin-Triebel spaces  are discussed at various places,
we refer, e.g., to the monographs \cite{Pe,RS,Tr1,Tr2,Tr08} and to the fundamental paper 
\cite{FJ2}.
We will not give definitions here and refer for this to the quoted literature.
\\
The present paper is a continuation of \cite{SS}, \cite{LSR} and \cite{KLSS}.
\\
{\bf Acknowledgement.} The authors would like to thank the referee
for a constructive hint to simplify the proof of Theorem \ref{anfang}.


\subsection*{Notation}


As usual, $\N$ denotes the natural numbers, $\N_0:= \N \cup \{0\}$,  $\Z$ denotes the integers and
$\R$ the real numbers.
If $X$ and $Y$ are two quasi-Banach spaces, then the symbol  $X \hookrightarrow Y$
indicates that the embedding is continuous.
The set of all linear  and bounded operators
$T : X \to Y$, denoted by $\cl (X,Y)$, is  equipped with the standard quasi-norm. 
As usual, the symbol  $c $ denotes positive constants
which depend only on the fixed parameters $s,p,q$ and probably on auxiliary functions,
 unless otherwise stated; its value  may vary from line to line.
Sometimes we  will use the symbols ``$ \lsim $''
and ``$ \gsim $'' instead of ``$ \le $'' and ``$ \ge $'', respectively. The meaning of $A \lsim B$ is given by: there exists a constant $c>0$ such that
 $A \le c \,B$. Similarly $\gsim$ is defined. The symbol
$A \asymp B$ will be used as an abbreviation of
$A \lsim B \lsim A$.\\
We shall use the following conventions throughout the paper:
\begin{itemize}
\item
If $E$ denotes a space of functions on $\R^d$ then by $RE$ we mean the subset of radial functions in $E$
and we endow this subset  with the same quasi-norm as the original space.
\item
Inhomogeneous Besov and Lizorkin-Triebel spaces are denoted by
$\bspq$ and $\fspq$, respectively.
If there is no reason to distinguish between these two scales we will use the
notation   $\aspq$. Similarly for the radial subspaces.
\item
If an equivalence class $[f]$ (equivalence with respect to coincidence almost everywhere) contains a continuous representative
then we call the class continuous and speak of values of $f$ at any point
(by taking the values of the continuous representative).
\item
Throughout the paper $\psi \in C_0^\infty (\Rd)$ denotes a  specific radial cut-off function,
i.e.,
$\psi (x)=1$ if $|x|  \le 1$ and $\psi (x)=0$ if $|x|\ge 3/2$.
\end{itemize}


\section{Main results} \label{main}


This section consists of  two parts.
In Subsection \ref{main1a} we concentrate on trace theorems which are the basis for the understanding of the higher regularity of radial functions outside the origin.
Subsection \ref{main5a} is devoted to the
study of decay and boundedness properties of radial functions in dependence
on their regularity.
To begin with we study the decay of radial functions near infinity.
Special emphasize is given to the limiting situation which arises
for $s=1/p$.
Then we continue with an investigation of the  behaviour of radial functions
near the origin.
Also here we investigate the limiting situations $s=d/p$ and $s=1/p$ in some detail.


\subsection{The characterization of the traces of radial subspaces}
\label{main1a}


Let $d\ge 2$.
Let $f:\Rd \to \C$ be a locally integrable radial function.
By using a Lebesgue point argument its restriction
\[
f_0 (t):= f(t,0, \ldots \, , 0)\, , \qquad t  \in \R\, .
\]
is well defined a.e. on $\R$.
However, this restriction need not be locally integrable.
A simple example is given by the function
\[
f(x):= \psi (x)\, |x|^{-1}\, , \qquad x \in \Rd\, .
\]
Furthermore, if we start with a measurable and even  function $g:\, \R \to \C$,
s.t. $g$ is locally integrable on all intervals $(a,b)$, $0 < a < b < \infty$,
then (again using a Lebesgue point argument) the function
\[
f(x):= g(|x|) \, , \qquad x \in \Rd\,
\]
is well-defined a.e. on $\Rd$ and is radial, of course.
In what follows we shall study properties of the associated  operators
\[
\tr \, : f \mapsto f_0 \qquad \mbox{and}\qquad \ext \, : \, g \mapsto f\, .
\]
Both operators are defined pointwise only.
Later on we shall have a short look onto the existence of the trace
in the distributional sense, see Subsection \ref{main30}.
Probably it would be more natural to deal with functions defined on $[0,\infty)$
in this context. However, that would result in more complicated descriptions of the trace spaces.
So, our target spaces will be spaces of even functions defined on $\R$.


\subsubsection{Traces of radial subspaces with $p=\infty$}
\label{main11}


The first result is maybe well-known but we did not find a reference for it.
Let $m \in \N_0$. Then
$C^m (\Rd)$ denotes the collection of all functions  $f: \Rd \to \C$
such that all derivatives $D^\alpha f$ of order $|\alpha| \le m$ exist,
are uniformly continuous and bounded.
We put
\[
\| \, f \, |C^m(\Rd)\| := \sum_{|\alpha|\le m} \| \, D^\alpha f\,
|L_\infty (\Rd)\|\, .
\]

\begin{T}
\label{anfang}
Let $d\ge 2$.
For $m \in \N_0$ the  mapping
$\tr$ is a linear  isomorphism of $RC^m (\Rd)$ onto $RC^m (\R)$
with  inverse  $\ext$.
\end{T}

\begin{Rem}\rm
If we replace uniformly continuous by continuous in the definition of the spaces
$C^m(\Rd)$ Theorem \ref{anfang}  remains true with the same proof.
\end{Rem}

\noindent
Using real interpolation it is not difficult to derive the following result
for the spaces of H\"older-Zygmund type.

\begin{T}\label{infty}
Let $s>0$ and let $0 < q \le \infty$. Then the  mapping $\tr $
is a linear  isomorphism of $RB^s_{\infty,q} (\Rd)$ onto $RB^s_{\infty,q} (\R)$
with inverse  $\ext$.
\end{T}


\subsubsection{Traces of radial subspaces with $p<\infty$}
\label{main12}


Now we turn to the description of the trace classes of radial Besov and
Lizorkin-Triebel spaces with $p<\infty$.
Again we start with an almost trivial result.
We need a further notation.
By $L_p (\R,w)$ we denote the weighted Lebesgue space equipped with the norm
\[
\|\, f \, |L_p (\R,w)\| := \Big(\int_{-\infty}^\infty |f(t)|^p\, w(t)\, dt\Big)^{1/p}
\]
with usual modification if $p=\infty$.

\begin{Lem}\label{simple}
Let $d\ge 2$.\\
{\rm (i)}
Let $0< p< \infty$. Then
$\tr \, :\,  RL_p (\Rd) \to RL_p (\R, |t|^{d-1})$ is a linear isomorphism with
inverse  $\ext$.
\\
{\rm (ii)}
Let $p = \infty$. Then $\tr \, :\,  RL_\infty (\Rd) \to RL_\infty (\R)$ is a linear isomorphism with  inverse $\ext $.
\end{Lem}

\noindent
In particular this means, that whenever  the Besov-Lizorkin-Triebel space
$\aspq (\Rd)$ is contained in $ L_1  (\Rd) +L_\infty (\Rd)$,
then $\tr $ is well-defined on its radial subspace.
This is in sharp contrast to the general theory of traces on
these spaces.
To guarantee that $\tr$ is meaningful on $\aspq (\Rd)$  one has to require
\[
s > \frac{d-1}{p} + \max \Big(0, \frac 1p -1\Big)\, ,
\]
cf. e.g. \cite{Ja}, \cite{FJ1},   \cite[Rem.~2.7.2/4]{Tr1} or
\cite{FJS}.
On the other hand we have
\[
\bspq (\Rd)\, , \, \fspq (\Rd) \hookrightarrow  L_1  (\Rd) +L_\infty (\Rd)
\]
if $s>d \, \max(0,\frac 1p -1)$, see, e.g., \cite{SiTr}.
Since
\[
d \, \max(0,\frac 1p -1) < \frac{d-1}{p} + \max \Big(0, \frac 1p -1\Big)
\]
we have the existence of $\tr$ with respect to $R\aspq (\Rd)$
for a wider range of parameters than for $\aspq(\Rd)$.
\\
Below we shall develop a description  of the traces of
the radial subspaces of $\bspq (\Rd)$ and $\fspq (\Rd)$
in terms of atoms. To explain this we need to introduce first an appropriate notion of an atom and second, adapted sequence spaces.

\bd\label{latom}
Let $L\ge 0$ be an integer. Let $I$ be a set either of the form
$I = [-a,a]$ or of the form $I=[-b,-a] \cup [a,b]$ for some $0 < a < b<\infty$.
An even function
$g \in C^L (\R)$ is called an even  $L$-atom centered at $I$
if
\[
\max_{t \in \R} | b^{(n)} (t)|  \le  | I |^{-n}\, , \qquad 0 \le n\le L\, .
\]
and if either
\[
\supp g \subset
[- \frac{3a}{2}, \, \frac{3a}{2}]  \qquad  \mbox{in case} \quad I=[-a,a]\, ,
\]
or
\[
\supp g \subset
[- \frac{3b-a}{2}, -\frac{3a-b}{2}]
\cup [ \frac{3a- b}{2}, \frac{3b-a}{2}]
\qquad \mbox{in case} \quad I =[-b, -a] \cup [a,b]\, .
\]
\ed

\begin{Def}
Let $0 <p< \infty$, $0 <q \le \infty$ and $s\in \R$.
Let
\[
{\chi}_{j,k}^\# (t):= \left\{\begin{array}{lll}
1 &\qquad & \mbox{if}\quad 2^{-j} k \le |t|\le 2^{-j} (k+1)\, ,\\
0 && \mbox{otherwise}\, .
                          \end{array}\right.
\qquad t \in \R\, .
\]
Then we define
\[
b^s_{p,q,d} := \bigg\{
s= (s_{j,k})_{j,k}\, : \quad
\| \,s \, | b^s_{p,q,d}\| =
\bigg(\sum_{j=0}^\infty 2^{j(s - \frac dp)q}\, \bigg(
\sum_{k=0}^\infty (1+ k )^{d-1}\,
|s_{j,k}|^p\bigg)^{q/p}
\bigg)^{1/q}<\infty \bigg\}\, .
\]
and
\beqq
f^s_{p,q,d} & := & \bigg\{s= (s_{j,k})_{j,k}\, : \\
\nonumber
\\
\| \,s \, | f^s_{p,q,d}\| & = &
\bigg\| \bigg(
\sum_{j=0}^\infty \, 2^{jsq}\,
\sum_{k=0}^\infty  \, |s_{j,k}|^q \, {\chi}^\#_{j,k}(\cdot)
\bigg)^{1/q}|L_p(\R, |t|^{d-1})\bigg\|\, <\infty \bigg\}\, ,
\eeqq
respectively.
\end{Def}

\begin{Rem}
Observe $b^s_{p,p,d}=f^s_{p,p,d}$ in the sense of equivalent quasi-norms.
\end{Rem}

\noindent
Adapted to these sequence spaces we define now function spaces
on $\R$.

\begin{Def}\label{tracedef}
Let $0 < p < \infty$, $0 <q \le \infty$, $s >0$ and $L \in \N_0$.
\\
{\rm (i)}
Then ${T}B^s_{p,q} (\R,L,d)$
is the collection of all functions $g:\R \to \C$ such that
there exists a decomposition
\be\label{new}
g(t)= \sum_{j=0}^\infty \sum_{k=0}^\infty s_{j,k}\, g_{j,k} (t)
\ee
(convergence in $L_{\max(1,p)} (\R, |t|^{d-1})$), where
the sequence $(s_{j,k})_{j,k}$ belongs to $b^s_{p,q,d}$
and the functions $g_{j,k}$ are  even  $L$-atoms centered at
either $ [- 2^{-j}, \,  2^{-j}]$ if $k=0$
or at
\[
 [- 2^{-j}(k+1), -2^{-j}k ] \: \cup \:
 [2^{-j}k, 2^{-j}(k+1)]
\]
if $k>0$. We put
\[
\| \, g \, |{T}B^s_{p,q} (\R,L,d)\|:= \inf \Big\{
\| \, (s_{j,k})\, |b^s_{p,q,d}\|\, : \quad (\ref{new})\quad \mbox{holds}
\Big\}\, .
\]
{\rm (ii)} Then ${T}F^s_{p,q} (\R,L,d)$
is the collection of all functions $g:\R \to \C$ such that
there exists a decomposition (\ref{new}),  where
the sequence $(s_{j,k})_{j,k}$ belongs to $f^s_{p,q,d}$ and
the functions $g_{j,k}$ are as in (i).
We put
\[
\| \, g \, |{T}F^s_{p,q} (\R,L,d)\|:= \inf \Big\{
\| \, (s_{j,k})\, |f^s_{p,q,d}\|\, : \quad (\ref{new})\quad \mbox{holds}
\Big\}\, .
\]
\end{Def}

\noindent
We need a few further notation. In connection with Besov and Lizorkin-Triebel spaces
quite often the following numbers occur:
\be\label{sigmap}
\sigma_p(d) := d\, \max \Big(0, \frac{1}{p}-1\Big)
\qquad \mbox{and}\qquad \sigma_{p,q}(d) := d\, \max \Big(0, \frac{1}{p}-1,
\frac 1q-1\Big)\, .
\ee
For a real number $s$ we denote by $[s]$ the integer part, i.e. the largest
integer $m$ such that $m \le s$.

\begin{T}\label{main5}
Let $d \ge 2$, $0 < p< \infty$ and $0 < q \le \infty$.
\\
{\rm (i)} Suppose $s>\sigma_p(d)$ and $L \ge [s]+1$.
Then the  mapping $\tr $
is a linear isomorphism of $R\bspq (\Rd)$ onto  $TB^s_{p,q}(\R,L,d)$
with inverse $\ext$.
\\
{\rm (ii)} Suppose $s>\sigma_{p,q}(d)$ and $L \ge [s]+1$.
Then the mapping $\tr $
is a linear  isomorphism of  $R\fspq (\Rd)$ onto  ${T}\fspq (\R,L,d)$
with inverse  $\ext$.
\end{T}

\begin{Rem}{\rm
Let $0 < p \le 1 < q \le \infty$. Then the spaces $RB^{\sigma_p}_{p,q}(\Rd)$ contain
singular distributions, see \cite{SiTr}.
In particular, the Dirac delta distribution belongs to $RB^{\frac dp -d}_{p,\infty}(\Rd)$,
see, e.g., \cite[Rem.~2.2.4/3]{RS}.
Hence, our pointwise defined mapping $\tr$ is not meaningful on those spaces,
or,
with other words, Theorem \ref{main5} does not extend to
values $s < \sigma_p (d)$.}
\end{Rem}

Outside the origin radial distributions are more regular.
We shall discuss several examples for this claim.

\begin{T}\label{spur1}
Let $d \ge 2$, $0 < p< \infty$ and $0 < q \le \infty$.
Suppose $s> \max (0, \frac 1p -1)$.
Let $f \in R\aspq (\Rd)$ s.t. $0 \not\in \supp f$. Then $f$ is a regular distribution
in $\cs'(\Rd)$.
\end{T}

\begin{Rem}
\label{smean}
\rm There is a nice and simple example which explains the sharpness of
the restrictions in Thm.~\ref{spur1}.
We consider the singular distribution $f$ defined by
\[
\varphi \: \mapsto \: \int_{|x|=1} \varphi (x)\, dx\, , \qquad \varphi \in \cs (\Rd)\, .
\]
By using  the wavelet characterization of Besov spaces,
it is not difficult to prove
that the spherical mean distribution $f$ belongs to the  spaces $B^{\frac 1p - 1}_{p,\infty} (\Rd)$ for all $p$.
\end{Rem}

\begin{T}\label{spur0}
Let $d \ge 2$, $0 < p< \infty$ and $0 < q \le \infty$.
Suppose $s> \max (0, \frac 1p -1)$.
Let $f \in R\aspq (\Rd)$ s.t. $0 \not\in \supp f$. Then
$f_0 = \tr f$ belongs to $A^s_{p,q} (\R)$.
\end{T}

\begin{Rem}{\rm
As mentioned above
\[
\aspq (\R) \hookrightarrow  L_1  (\R) + L_\infty (\R) \qquad \mbox{if}
\qquad s> \sigma_p (1) = \max \Big(0, \frac 1p -1\Big) \, ,
\]
which shows again that we deal with regular distributions.
However, in Thm. \ref{spur0} some additional regularity is proved.
}
\end{Rem}


\subsubsection{Traces of radial subspaces of Sobolev spaces}
\label{main16}


Clearly, one can expect that the description of the traces of
radial Sobolev spaces can be given in more elementary terms.
We discuss a few examples without having the complete theory.

\begin{T}\label{main10}
Let $d \ge 2$ and  $1 \le p< \infty$.
\\
{\rm (i)} The  mapping $\tr $ is a linear isomorphism (with inverse $\ext$)
of $RW^1_p(\Rd)$ onto
the closure of  $RC^{\infty}_0 (\R)$ with respect to the norm
\[
\|\, g\, |L_p(\R,|t|^{d-1})\| + \|\, g'\, |L_p(\R,|t|^{d-1})\|\, .
\]
{\rm (ii)} The  mapping $\tr $ is a linear isomorphism (with inverse $\ext$)
of $RW^2_p(\Rd)$ onto
the closure of  $RC^{\infty}_0 (\R)$ with respect to the norm
\[
\|\, g\, |L_p(\R,|t|^{d-1})\| + \|\, g'\, |L_p(\R,|t|^{d-1})\|
+ \|\, g'/r \, |L_p(\R,|t|^{d-1})\|
 + \|\, g'' \, |L_p(\R,|t|^{d-1})\|  \, .
\]
\end{T}

\begin{Rem}{\rm
Both statements have elementary proofs, see (\ref{w-100}) for (i). However,
the complete extension to higher order Sobolev spaces is open.}
\end{Rem}

There are several ways to define Sobolev spaces on $\Rd$.
For instance, if $1 < p < \infty$ we have
\be\label{eq-31}
f \in W^{2m}_p (\Rd) \quad \Longleftrightarrow
\quad f\in L_p (\Rd)\quad \mbox{and}\quad \Delta^m f \in L_p (\Rd)\, .
\ee
Such an equivalence does not extend to $p=1$ or $p=\infty$
if $d\ge 2$, see \cite[pp.~135/160]{Stein1}.
Recall that the Laplace operator $\Delta$
applied to a radial function yields a radial function.
In particular we have
\be\label{eq-20}
\Delta f (x) = D_r f_0 (r):= f_0''(r) +
\frac{d-1}{r}\, f_0'(r) \, , \qquad r=|x|\, ,
\ee
in case that $f$ is radial and $\tr f = f_0$.
Obviously, if  $f \in RC_0^\infty (\Rd)$, then
\beq\label{laplace}
\| \, f \, |L_p (\Rd)\| & + & \| \, \Delta^m f \, |L_p (\Rd)\|
\\
& = &
\Big(\frac{\pi^{d/2}}{\Gamma(d/2)}\Big)^{1/p} \,\Big(
\|\, f_0 \, |L_p (\R, |t|^{d-1})\|
+ \|\, D_r^m f_0 \, |L_p (\R, |t|^{d-1})\|\Big)\, .
\nonumber
\eeq
This proves the next characterization.

\begin{T}\label{main1}
Let $1 <p<\infty$ and $m \in \N$. Then the mapping
$\tr$ yields a linear  isomorphism (with inverse $\ext$)
of $RW^{2m}_p (\Rd)$ onto the closure of
$RC^{\infty}_0 (\R)$ with respect to the norm
\[
\|\, f_0 \, |L_p (\R, |t|^{d-1})\| +  \|\, D_r^m f_0 \, |L_p (\R, |t|^{d-1})\| \, .
\]
\end{T}

\begin{Rem}{\rm
By means of Hardy-type inequalities one can simplify the terms
 $\|\, D_r^m f_0 \, |L_p (\R, |t|^{d-1})\|$ to some extent,
 see Theorem \ref{main10}(ii) for a comparison. We do not go into detail.}
\end{Rem}


\subsubsection{The trace in $\cs'(\R)$}
\label{main30}


Many times applications of traces are connected with
boundary value problems. In such a context the continuity of $\tr$
considered as a mapping into $S'$ is essential.
Again we consider the simple situation of the $L_p$-spaces first.

\begin{Lem}\label{dp}
Let $d  \ge 2$ and let $0 < p  < \infty$. Then
$RL_p (\R, |t|^{d-1}) \subset S' (\R)$ if and only if $d<p$.
\end{Lem}

From the known embedding relations of $R\aspq (\Rd)$
into $L_u$-spaces one obtains one  half of the proof of the following
 general result.

\begin{T}
\label{main8}
Let $d\ge 2$, $0 <p<\infty$, and  $0<q\le \infty$.
\\
{\rm (a)} Let   $s > \sigma_p(d)$ and $L\ge [s]+1$.
Then the  following assertions are equivalent:\\
\mbox{\hspace{7mm}{\rm (i)}}
The mapping
$\tr$ maps $R\bspq (\Rd)$ into $S'(\R)$.\\
\mbox{\hspace{7mm}{\rm (ii)}}
The mapping
$\tr: \, R\bspq (\Rd) \to S'(\R)$ is continuous.\\
\mbox{\hspace{7mm}{\rm (iii)}}
We have ${T}\bspq (\R,L,d) \hookrightarrow S'(\R)$.\\
\mbox{\hspace{7mm}{\rm (iv)}}
We have either  $s > d(\frac 1p - \frac 1d)$
or $s= d(\frac 1p - \frac 1d) $ and $q\le 1$.
\\
{\rm (b)}  Let   $s > \sigma_{p,q}(d)$ and $L\ge [s]+1$.
Then following assertions are equivalent:
\\
\mbox{\hspace{7mm}{\rm (i)}}
The mapping $\tr$ maps $R\fspq (\Rd)$ into $S'(\R)$.\\
\mbox{\hspace{7mm}{\rm (ii)}}
The mapping $\tr: \, R\fspq (\Rd) \to S'(\R)$ is continuous.\\
\mbox{\hspace{7mm}{\rm (iii)}} We have ${T}\fspq (\R,L,d) \hookrightarrow S'(\R)$.\\
\mbox{\hspace{7mm}{\rm (iv)}} We have either  $s > d(\frac 1p - \frac 1d)$
or $s = d(\frac 1p - \frac 1d)$ and $0< p \le 1$.
\end{T}


\subsubsection{The trace in $\cs'(\R)$ and weighted function spaces of Besov
and Lizorkin-Triebel type}
\label{main31}


Weighted function spaces of Besov and Lizorkin-Triebel type,
denoted by $\bspq (\R, w)$ and $\fspq (\R, w)$, respectively,
are a well-developed subject in the literature, we refer to \cite{bui-1,bui-2,Ry}.
Fourier analytic definitions as well as characterizations by atoms are given under various restrictions on the weights, see
e.g. \cite{Bow1,bui-1,bui-2,HP,IS,ST}.
In this subsection we are interested in these spaces
with respect to the weights $w_{d-1}(t):= |t|^{d-1}$, $t\in \R$, $d\ge 2$.
Of course, these weights belong to the Muckenhoupt class $\ca_\infty$, more exactly
$w_{d-1} \in {\mathcal A}_r$ for any $r>d$, see \cite{stein}.

\begin{T}\label{main32}
Let $d \ge 2$, $0 <p< \infty$, and  $0 < q \le \infty$.
\\
{\rm (i)} Suppose $s>\sigma_{p}(d)$ and let $L\ge [s]+1$.
If ${T}\bspq (\R,L,d) \hookrightarrow S'(\R)$ (see Theorem \ref{main8}),
then
${T}\bspq (\R,L,d) = R\bspq (\R, w_{d-1})$
in the sense of equivalent quasi-norms.\\
{\rm (ii)} Suppose $s>\sigma_{p,q}(d)$ and let $L\ge [s]+1$.
If
${T}\fspq (\R,L,d) \hookrightarrow S'(\R)$ (see Theorem \ref{main8}),
then ${T}\fspq (\R,L,d) = R\fspq (\R, w_{d-1})$
in the sense of equivalent quasi-norms.
\end{T}

\beginpicture

\setcoordinatesystem units <0.8cm,0.8cm>
\unitlength0.8cm

\put {\vector (0,1) {8} } [Bl] at 1 0
\put {\vector (1,0) {13} } [Bl] at 0 1

\thicklines
\plot 5.5 1  12.5 7.125 /
\plot 10.5 1 14 4.0625 /

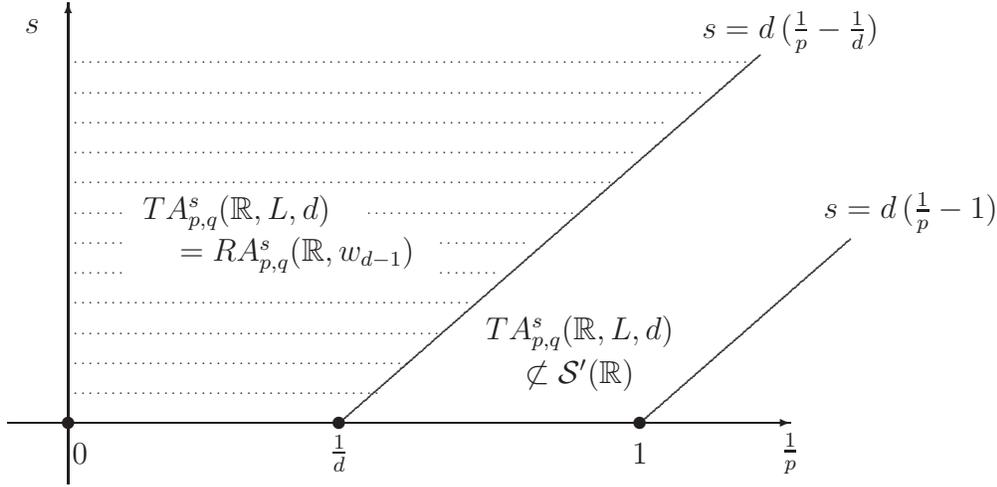
\put {Fig.~1: Existence of the trace} at 6 -1

\put {$0$} at 1.2 0.5
\put {$\frac{1}{p}$} at 13 0.5
\put {$s$} at 0.4 7.6
\put {$TA^{s}_{p,q}(\R,L,d)$} at 3.8 4.5
\put {$= RA^{s}_{p,q}(\R,w_{d-1})$} at 4.8 3.8
\put {$TA^{s}_{p,q}(\R,L,d)$} at 9.5 2.5
\put {$\not\subset \cs'(\R)$} at 9.5 1.8

\put {$\bullet$} at 1 1
\put {$1$} at 10.5 0.5
\put {$\bullet$} at 10.5 1
\put {$\bullet$} at 5.5 1

\put {$\frac{1}{d}$} at 5.5 0.5
\put {$s = d\, (\frac{1}{p} - \frac{1}{d})$} at 13 7.5
\put {$s = d\, (\frac{1}{p} - 1)$} at 15 4.5

{\setdots <0.1cm>

\plot 1 1.5 6.2 1.5 /
\plot 1 2 6.7 2 /
\plot 1 2.5 7.2 2.5 /
\plot 1 3 7.7 3 /
\plot 1 3.5 2 3.5 /
\plot 7.2 3.5 8.2 3.5 /
\plot 1 4 2 4 /
\plot 7.2 4 8.7 4 /
\plot 1 4.5 2 4.5 /
\plot 6 4.5 9.5 4.5 /
\plot 1 5 10.1 5 /
\plot 1 5.5 10.4 5.5 /
\plot 1 6 11 6 /
\plot 1 6.5 11.6 6.5 /
\plot 1 7 12.3 7 /

}

\endpicture
\hfill\\

\begin{Rem}\label{delta}
{\rm
We add some statements concerning the regularity of the
most prominent singular distribution, namely
$\delta ~: \: \varphi \to \varphi (0)$, $\varphi \in S(\Rd)$.
This tempered distribution has the following
regularity properties:
\begin{itemize}
\item First we deal with the situation on $\Rd$. We have
$\delta \in RB^{\frac dp - d}_{p,\infty} (\Rd)$ (but
$\delta \not\in RB^{\frac dp - d}_{p,q} (\Rd)$ for  $q<\infty$ and
$\delta \not\in RF^{\frac dp - d}_{p,\infty} (\Rd)$), see, e.g., \cite[Rem.~2.2.4/3]{RS}.
\item
Now we turn to the situation on $\R$. By using more or less the same arguments as on $\Rd$ one can show
$\delta \in B^{\frac dp - 1}_{p,\infty} (\R,w_{d-1})$
(but  $\delta \not\in B^{\frac dp - 1}_{p,q} (\R,w_{d-1})$ for any $q<\infty$
and $\delta \not\in F^{\frac dp - 1}_{p,\infty} (\R,w_{d-1})$).
\end{itemize}}
\end{Rem}


\subsubsection{The regularity of radial functions  outside
the origin}
\label{main14}


Let $f$ be a radial function such that
$\supp f \, \subset\,  \{x\in \Rd~: \quad |x|\ge \tau \}$
for some $\tau >0$.
Then the following inequality is obvious:
\[
\| \, f_0 \, | L_p (\R)\|  \le  \,  \tau^{-(d-1)/p}\,
\Big(\frac{\Gamma (d/2)}{\pi^{d/2}}\Big)^{1/p}\,
\| \, f \, | L_p (\Rd) \|\, .
\]
An extension to first or second order Sobolev spaces
can be done by using Theorem \ref{main10}.
However, an extension to all spaces under consideration here
is less obvious.
Partly it could be done by interpolation, see Proposition \ref{interpol2},
but we prefer a different way (not to exclude $p<1$).
We shall compare the atomic decompositions in Theorem \ref{main5}
with the known atomic  and wavelet characterizations of
$\bspq (\R)$ and $\fspq (\R)$.

\begin{Cor}\label{nice}
Let $\tau >0$.
Let $d\ge 2$, $0 <p< \infty$,  and $0 < q \le \infty$.\\
{\rm (i)} We suppose $s> \sigma_p(d)$.
If  $f \in R\bspq (\Rd)$  such that
\begin{equation}\label{supp}
\supp f \, \subset\,  \{x\in \Rd~: \quad |x|\ge \tau \}
\end{equation}
then its trace $f_0$ belongs to $\bspq (\R)$.
Furthermore, there exists a constant $c$ (not depending on $f$ and $\tau $)
such that
\begin{equation}\label{LSnewb}
\| \, f_0 \, | \bspq(\R)\|  \le c \,  \tau^{-(d-1)/p}\,
\| \, f \, | \bspq (\Rd) \|
\end{equation}
holds for all such functions $f$ and all $\tau >0$.\\
{\rm (ii)} We suppose $s>  \sigma_{p,q}(d)$.
If  $f \in R\fspq (\Rd)$  such that (\ref{supp}) holds,
then its trace $f_0$ belongs to $\fspq (\R)$.
Furthermore, there exists a constant $c$ (not depending on $f$ and $\tau$) such that
\begin{equation}\label{LSnewf}
\| \, f_0 \, | \fspq(\R)\|  \le c\,  \tau^{-(d-1)/p}\,
\| \, f \, | \fspq (\Rd) \|
\end{equation}
holds for all such functions $f$ all $\tau >0$.
\end{Cor}

We wish to mention that Corollary \ref{nice} has a partial inverse.

\begin{Cor}\label{ok}
Let $d\ge 2$, $0 <p< \infty$,  $0 < q \le \infty$ and $0 < a < b < \infty$.\\
{\rm (i)} We suppose $s>  \sigma_p(d)$.
If  $g \in R\bspq (\R)$  such that
\be\label{suppp}
\supp g \, \subset\,  \{x\in \R~: \quad a \le  |x|\le b\}
\ee
then the radial function
$f:= \ext g$ belongs to $R\bspq (\Rd)$
and there exist positive constants $A,B$ such that
\[
A\, \| \, g \, | \bspq(\R)\|\le
\| \, f \, | \bspq (\Rd) \| \le B\, \| \, g \, | \bspq(\R)\| \, .
\]
{\rm (ii)} We suppose $s>\sigma_{p,q}(d)$.
If  $g \in R\fspq (\R)$  such that (\ref{suppp}) holds,
then the radial function
$f:= \ext g$ belongs to $R\fspq (\Rd)$
and there exist positive constants $A,B$ such that
\[
A\, \| \, g \, | \fspq(\R)\|\le
\| \, f \, | \fspq (\Rd) \| \le B\, \| \, g \, | \fspq(\R)\| \, .
\]
\end{Cor}

For our next result we need H\"older-Zygmund spaces.
Recall, that $C^s (\Rd) = B^s_{\infty,\infty} (\Rd)$ in the sense of equivalent norms
if $s \not\in \N_0$.
Of course, also the spaces  $B^s_{\infty,\infty} (\Rd)$ with $s \in \N$ allow a characterization by
differences. We refer to \cite[2.2.2, 2.5.7]{Tr1} and \cite[3.5.3]{Tr2}.
We shall use the abbreviation
\[
\cz^s (\Rd)=  B^s_{\infty,\infty} (\Rd)\, , \qquad s>0\, .
\]
Taking into account the well-known embedding relations for Besov
as well as for Lizorkin-Triebel spaces, defined on $\R$, Thm. \ref{spur0} implies in particular:

\begin{Cor}\label{smooth}
Let $d\ge 2$, $0 <p< \infty$, $0 < q \le \infty$, and $s> 1/p$.
Let $\varphi$ be a smooth  radial function, uniformly bounded together with all its derivatives, and such that $0 \not\in \supp \varphi$.
If $f \in R\aspq (\Rd)$, then $\varphi \,  f \in \cz^{s-1/p}(\Rd)$.
\end{Cor}

\begin{Rem}
{\rm P.L.~Lions \cite{Li} has proved the counterpart of
Corollary \ref{smooth} for first order Sobolev spaces.
We also dealt in \cite{SS} with these problems.}
\end{Rem}

Finally, for later use, we would like to know when the radial functions are continuous out of the origin.

\begin{Cor}\label{uc}
Let $\tau >0$.
Let $d\ge 2$, $0 <p< \infty$, and $0 < q \le \infty$. \\
{\rm (i)} If either  $s>1/p$ or $s=1/p$ and $q\le 1$ then $f \in R\bspq (\Rd)$
is  uniformly continuous on the set $|x|\ge \tau $.
\\
{\rm (ii)}
If either  $s>1/p$ or $s=1/p$ and $p\le 1$ then $f \in R\fspq (\Rd)$
is  uniformly continuous on the set $|x|\ge \tau $.
\end{Cor}

By looking at the restrictions in Cor. \ref{uc} we introduce the following set
of parameters.

\bd
{\rm (i)}
We say $(s,p,q)$ belongs to the set $U(B)$ if
$(s,p,q)$ satisfies the restrictions in part (i) of Cor. \ref{uc}.
\\
{\rm (ii)}
The triple $(s,p,q)$ belongs to the set $U(F)$ if
$(s,p,q)$ satisfies the restrictions in part (ii) of Cor. \ref{uc}.
\ed

\begin{Rem}\label{border}
{\rm
(a) The abbreviation $(s,p,q) \in U(A)$
will be used  with the obvious meaning.
\\
(b) Let $1\le p=p_0<\infty $ be fixed. Then there is always a  largest space
in the set
\[
\{B^s_{p_0,q}  (\Rd)\, : \: (s,p_0,q) \in U(B)\} \quad \cup \quad
\{F^s_{p_0,q}  (\Rd)\, : \: (s,p_0,q) \in U(F)\}\, .
\]
This space is given either by
$F^{1}_{1,\infty}(\Rd)$ if $ p_0=1$ or by
$B^{1/p_0}_{p_0,1}(\Rd)$ if $ 1 <p_0<\infty$.
If $p_0<1$, then obviously $B^{1/p_0}_{p_0,1}(\Rd)$ is the largest Besov space  and  $F^{1/p_0}_{p_0,\infty}(\Rd)$ is the largest Lizorkin-Triebel space in the above family.
However, these  spaces are incomparable.}
\end{Rem}


\subsection{Decay and boundedness properties of radial functions}
\label{main5a}


We deal with improvements of Strauss'   {\em Radial Lemma}.
Decay can only be expected if we measure smoothness in function spaces built on
$L_p (\Rd)$ with $p<\infty$.\\
It is instructive to have a short look onto the case
of first order Sobolev spaces.
Let $f=g(r(x))\in RC_0^\infty(\R^d)$. Then
\[
\frac{\partial f}{\partial x_i}(x)= g'(r) \, \frac{x_i}{r}, \quad
r=|x|>0 \, , \qquad i=1, \, \dots \, ,d \, .
\]
Hence
\be\label{w-100}
||\, \, \,  |\nabla f(x)|\,\,\,   |L_p(\R^d)|| = c_d \, ||\,  g' \, |L_p(\R,|t|^{d-1})||\, ,
\ee
where $1 \le p <\infty$.
Next we apply  the identity
\[
 g(r)=-\int_r^\infty g'(t)dt
\]
and obtain
\[
|g(r)|\le \int_r^\infty |g'(t)|dt\le r^{-(d-1)}\int_r^\infty t^{d-1}|g'(t)|dt.
\]
This extends to all functions in $RW^1_1 (\Rd)$ by a density argument.
On this elementary way we have proved the inequality
\be\label{w-04}
|x|^{d-1}\,  |f(x)| = r^{d-1} \, |g(r)| \le \frac{1}{c_d} \, \int_{|x|>r } \, |\nabla f(x)| \, dx
 \le \frac{1}{c_d} \, \| \, \nabla f(x)\, \|_1 \, .
\ee
This inequality can be interpreted in several ways:
\begin{itemize}
\item
The possible unboundedness in the origin is limited.
\item
There is some decay, uniformly in   $f$,  if $|x|$ tends to $+\infty$.
\item
We have $\lim_{|x|\to \infty} |x|^{d-1}\,  |f(x)| = 0$ for all $f \in RW^1_1(\Rd)$.
\item
It makes sense to switch to homogeneous function spaces, since in (\ref{w-04})
 only the norm of the  homogeneous Sobolev space occurs (for this, see \cite{CO} and \cite{SSV2}).
\end{itemize}
We shall show that all these phenomena will occur also in the general context of
radial subspaces of Besov and Lizorkin-Triebel spaces.


\subsubsection{The behaviour of radial functions near infinity}
\label{main2}


Suppose $(s,p,q) \in U(A)$. Then
$f \in R{A}^s_{p,q}(\Rd)$ is uniformly continuous near infinity and belongs to $L_p(\Rd)$.
This implies
$\lim_{|x| \to \infty} \,  |f(x)| =0$. However, much more is true.

\begin{T}\label{decay4}
Let $d\ge 2$, $0<  p  < \infty $, and $0 < q \le \infty$.
\\
{\rm (i)} Suppose $(s,p,q) \in U(A)$.
Then there exists a constant $c$ s.t.
\be\label{eq-72c}
|x|^{(d-1)/p}\,  |f(x)| \le c \, \| \, f \, |{A}^s_{p,q} (\Rd)\| \,
\ee
holds for all $|x|\ge 1$ and all $f \in R{A}^s_{p,q} (\Rd)$.
\\
{\rm (ii)} Suppose $(s,p,q) \in U(A)$. Then
\be\label{eq-91b}
 \lim_{|x| \to \infty} \, |x|^{\frac {d-1}{p}}\, |f(x)| =0
\ee
holds for all $ f \in  R{A}^s_{p,q} (\Rd)$.\\
{\rm (iii)} Suppose $(s,p,q) \in U(A)$. Then there exists a constant $c>0$ such that for all  $x$, $|x|>1$, there exists  a smooth radial function
$ f \in  R{A}^s_{p,q} (\Rd)$, $ \|\, f \, |R{A}^s_{p,q} (\Rd)\|=1$, s.t.
\be\label{eq-91}
 |x|^{\frac {d-1}{p}}\, |f(x)| \ge c  \, .
\ee
{\rm (iv)} Suppose $(s,p,q) \not\in U(A)$ and $\frac{1}{p}>\sigma_p(d)$. We assume also that $\frac{1}{p}>\sigma_{q}(d)$ in the $F$-case.
Then, for all sequences $(x^j)_{j=1}^\infty \subset \Rd\setminus \{0\}$ s.t. $\lim_{j \to \infty} |x^j| =\infty$, there exists a  radial function
$ f \in  R{A}^s_{p,q} (\Rd)$, $ \|\, f \, |R{A}^s_{p,q} (\Rd)\|=1$, s.t.
$f$ is unbounded in any neighborhood of $x^j$, $j\in \N$.
\end{T}

\begin{Rem} \label{flambda}
{\rm
(i) Increasing $s$ (for fixed $p$) is not improving the decay rate.
In the case of Banach spaces, i.e.,  $p,q\ge 1$, the additional assumptions in  point (iv) are always fullfiled.  
Hence, the largest spaces, guaranteeing the decay rate $(d-1)/p$, are spaces with $s=1/p$, see Remark \ref{border}.
The dependence of the decay on the parameters $s$ and $p$ is illustrated in Fig.~1 below.
\\
(ii) Observe that in (iii) the function depends on $|x|$.
There is no function in $ R{A}^s_{p,q} (\Rd)$
such that (\ref{eq-91}) holds for all $x$, $|x|\ge 1$, simultaneously.
The naive construction $f(x):= (1-\psi (x)) \, |x|^{\frac{1-d}{p}}$, $x\in \Rd$,
does not belong to $L_p (\Rd)$.
\\
(iii)
If one switches from inhomogeneous spaces to
the larger homogeneous spaces of Besov and Lizorkin-Triebel type, then the decay rate becomes smaller. It will depend
also on $s$, see \cite{CO} and \cite{SSV2} for details.
\\
(iv) Of course, formula (\ref{eq-72c}) generalizes the estimate
(\ref{strauss}). Also Coleman, Glazer and Martin \cite{CGM} have dealt with (\ref{strauss}). P.L.~Lions \cite{Li} proved a $p$-version of the {\em Radial Lemma}.\\
(v) In case $A=B$ the theorem has been proved in \cite{SS}. For $A=F$ and $1<p<\infty$ the result follows 
from the embedding $B^s_{p,1} (\Rd) \hookrightarrow F^s_{p,q} (\Rd) \hookrightarrow B^s_{p,\infty} (\Rd)$,
valid for all $s$ and all $q$. The novelity in Theorem. \ref{decay4}
consists in a discussion of the case $RF^{1/p}_{p,\infty} (\Rd)$, $0 < p \le 1$, 
see also Remark \ref{border}.
\\
(vi)
Originally the {\em Radial Lemma} has been used to prove
compactness of embeddings of radial Sobolev spaces into
$L_p$-spaces, see \cite{CGM}, \cite{Li}. In the framework of
radial subspaces of
Besov and of Lizorkin-Triebel spaces compactness of embeddings has been investigated in \cite{SS}.
There we have given a final answer, i.e., we proved an {\em if, and only if,} assertion.\\
(vii) Compactness of embeddings of radial subspaces of homogeneous
Besov and of Lizorkin-Triebel spaces will be investigated in \cite{SSV2}.
}
\end{Rem}

\beginpicture

\setcoordinatesystem units <0.8cm,0.8cm>
\unitlength0.8cm

\put {\vector (0,1) {8} } [Bl] at 1 0
\put {\vector (1,0) {13} } [Bl] at 0 1

\thicklines
\plot 7 1 13 8 /
\plot 1 1 13 7 /
\plot 7 1 13 4 /
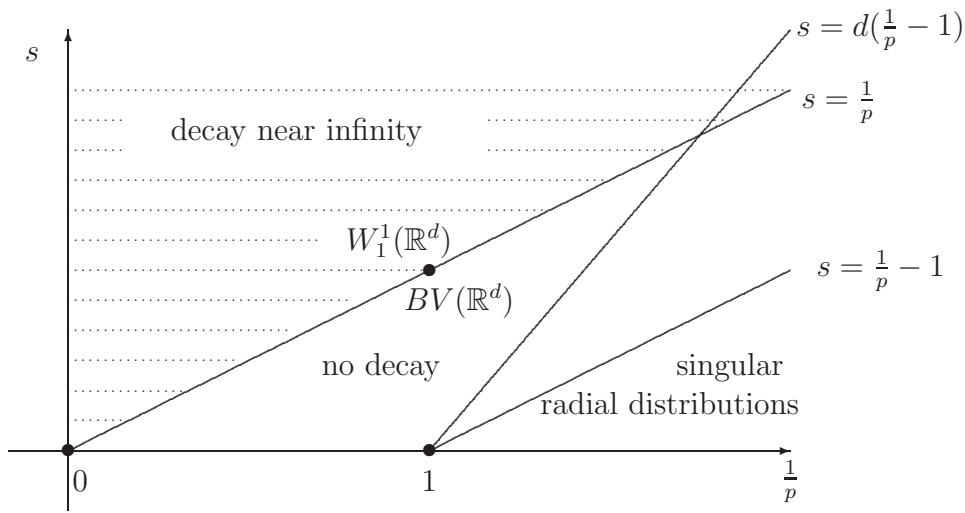
\put {Fig.~2: Decay at infinity} at 6 -1

\put {$0$} at 1.2 0.5
\put {$\frac{1}{p}$} at 13 0.5
\put {$s$} at 0.4 7.6
\put {decay near infinity} at 4.8 6.3
\put {$ W^{1}_{1}(\Rd)$} at 6.5 4.5
\put {$ BV(\Rd)$} at 7.5 3.5
\put {no decay} at 6.2 2.4

\put {$\bullet$} at 1 1
\put {$1$} at 7 0.5
\put {$\bullet$} at 7 1
\put {$\bullet$} at 7 4
\put {singular} at 12 2.4
\put {radial distributions} at 11 1.8

\put {$s=d(\frac{1}{p}-1)$} at 14.5 8
\put {$s = \frac{1}{p}$} at 13.8 6.8
\put {$s = \frac{1}{p}- 1$} at 14.5 4

{\setdots <0.1cm>

\plot 1 1.5 2 1.5 /
\plot 1 2 3 2 /
\plot 1 2.5 3.8 2.5 /
\plot 1 3 4.7 3 /
\plot 1 3.5 5.7 3.5 /
\plot 1 4 7 4 /
\plot 1 4.5 5.2 4.5 /
\plot 1 5 9 5 /
\plot 1 5.5 10 5.5 /
\plot 1 6 2 6 /
\plot 1 6.5 2 6.5 /
\plot 1 7 13 7 /
\plot 8 6.5 12 6.5 /
\plot 8 6 11 6 /
}

\endpicture
\hfill\\


\subsubsection{The behaviour of radial functions near
infinity -- borderline cases}
\label{border1}


As indicated in Remark \ref{border},
within the scales of Besov and Lizorkin-Triebel spaces
the borderline cases  for the decay rate $(d-1)/p$ are
either  $F^{1}_{1,\infty}(\Rd)$ if $p = 1$ or
$B^{1/p}_{p,1}(\Rd)$ if $ 1 <p<\infty$.
Now we turn to spaces  which do not belong to these scales
and where the elements of the  radial subspaces have such a
decay rate.
Hence, we are looking for spaces of radial functions with a simple norm
which satisfy (\ref{eq-72c}).
The Sobolev space $RW^1_1 (\Rd)$ is such a candidate for which
(\ref{eq-72c}) is already known, see \cite{Li}.
But this is not the end of the story. Also for the radial
functions of bounded variation such a decay estimate is true.

\begin{T}\label{decaybv}
Let $d\ge 2$.
Then there exists  constant $c$ s.t.
\be\label{ws-01}
|x|^{d-1}\,  |f(x)| \le c \, \| \, f \, | BV (\Rd)\| \,
\ee
holds for all $|x|>0$ and all $f \in RBV (\Rd)$.
Also
\be\label{ws-02}
 \lim_{|x| \to \infty} \, |x|^{d-1}\, |f(x)| =0
\ee
is true for all $ f \in  RBV(\Rd)$.
\end{T}

\begin{Rem}\label{repre}
{\rm
(i) Both assertions, (\ref{ws-01}) and  (\ref{ws-02}), require an interpretation since, in contrast to the classical definition of $BV (\R)$, the
spaces $BV(\R^d)$, $d\ge 2$, are  spaces of equivalence classes, see
Subsection \ref{spurbv}. Nevertheless, in every equivalence class $[f]\in BV(\R^d)$,
there is a representative $\tilde f\in[f]$, such that
\[
 |\tilde f(x)|\le \limsup_{y\to x} |f(y)|
\]
(simply take  $\tilde f(x):=f(x)$ in every Lebesgue point $x$ of $f$ and
$\tilde f(x):=0$ otherwise).
Hence,  (\ref{ws-01}) and  (\ref{ws-02}) have to be interpreted
as follows: whenever we work with values of the equivalence class
$[f]$ then we mean  the function values of the above representative
$\tilde f$.\\
(ii)
Notice that $F^1_{1,\infty} (\Rd)$ and $BV (\Rd)$ are incomparable.
\\
(iii) Observe, as in case of the {\em Radial Lemma}, that (\ref{ws-01})
holds for $x \neq 0$.
}
\end{Rem}

As a preparation for Theorem \ref{decaybv} we shall characterize the traces of radial elements in $BV(\Rd)$.
This seems to be of independent interest.
For this reason  we are forced to introduce weighted spaces of functions of bounded variation on the positive half axis. We put $\R^+ := (0,\infty)$.
As usual, $|\nu|$ denotes the total variation of the measure $\nu$, see, e.g.,
\cite[Chapt.~6]{Ru}.

\begin{Def}\label{BVDef2}
{\rm (i)} A function $\varphi\in C([0,\infty))$
belongs to $C_c^1([0,\infty))$ if it is continuously differentiable
on $\R^+$, has  compact support, satisfies $\varphi(0)=0$ and $\displaystyle
\lim_{t\to 0^+}\varphi'(t)=\varphi'(0)=\lim_{t\to 0^+}\frac{\varphi(t)}{t}$
exists and is finite.
\\
{\rm (ii)} A function $g\in L_1(\R^+,t^{d-1})$ is said to belong to $BV(\R^+,t^{d-1})$
if there is a signed Radon measure $\nu$ on $\R^+$ such that
\be\label{f-06}
\int_0^\infty g(t)\, [\varphi(s)s^{d-1}]'(t)\, dt
= - \, \int_0^\infty \varphi(t)\, t^{d-1}\, d\nu(t)\, ,\qquad \forall
\varphi\in C_c^1([0,\infty))
\ee
and
\be\label{f-07}
\|\, g\, |BV(\R^+,t^{d-1})\|
:=\|\, g \, |L_1(\R^+,t^{d-1})\|+ \, \int_0^\infty r^{d-1}\, d|\nu|(r)
\ee
is finite.
\end{Def}

By using these new spaces we can prove the following trace theorem.
With a slight abuse of notation $\ext$ denotes here  the radial extension 
of a function defined on $[0,\infty)$ (instead of $\R$).

\begin{T}\label{BVThm1}
Let $g$ be a measurable function on $\R^+$. Then $ \ext g \in BV(\R^d)$ if, and only if $g\in BV(\R^+,t^{d-1})$ and
\[
 \| \, \ext g \, |BV(\R^d)\| \asymp \|\, g \, |BV(\R^+,t^{d-1})\,\| \, .
\]
\end{T}

\subsubsection*{Spaces with $1<p<\infty$}

For $1<p<\infty$ one could use interpolation between $p=1$
and $p=\infty$ to obtain spaces with the decay rate $(d-1)/p$.
The largest spaces with this respect are obtained by the
real method.
Let $M_p (\Rd):= (RC (\Rd), RBV(\Rd))_{\Theta, \infty}$, $\Theta = 1/p $.
Then (\ref{eq-72c}) holds for all elements $f \in M_p (\Rd)$.
The disadvantage of these classes $M_p (\Rd)$ lies in the fact
that elementary descriptions of $M_p (\Rd)$ are not known.
However, at least some embeddings are known.
From
\[
RB^{1/p}_{p,1}(\Rd)= [RB^0_{\infty,1} (\Rd), RB^1_{1,1,}(\Rd)]_{\Theta}
\hookrightarrow (RC (\Rd), RBV(\Rd))_{\Theta, \infty}\, , \qquad
\Theta = 1/p \, ,
\]
(combine Proposition \ref{interpol2} with \cite[Thm.~4.7.1]{BL}),
we get back Theorem \ref{decay4} (i), but only  in case $1<p<\infty$.


\subsubsection{The behaviour of radial functions near the origin}
\label{main3}


At first we mention that the embedding relations with respect to
$L_\infty (\Rd)$ do not change when we switch from $\aspq  (\Rd)$ to its
radial subspace $R\aspq  (\Rd)$.

\begin{Lem}\label{bi}
{\rm (i)}
The embedding $R\bspq (\Rd) \hookrightarrow L_\infty(\Rd)$ holds
if and only if either $s>d/p$ or $s=d/p$ and $q\le 1$.
\\
{\rm (ii)}
The embedding $R\fspq (\Rd) \hookrightarrow L_\infty (\Rd)$ holds
if and only if either $s>d/p$ or $s=d/p$ and $p\le 1$.
\end{Lem}

The explicit counterexamples will be  given in Lemma \ref{sharp} below.
Hence, unboundedness can only happen in case $s \le d/p$.

\begin{T}\label{decay2}
Let $d\ge 2$, $0<  p  < \infty  $ and $0 < q \le \infty$.
\\
{\rm (i)} Suppose $(s,p,q) \in U(A)$
and $s<\frac d p$.
Then there exists a constant $c$ s.t.
\be\label{eq-72}
|x|^{\frac dp -s}\,  |f(x)| \le c \, \| \, f \, |R{A}^s_{p,q} (\Rd)\| \,
\ee
holds for all $0< |x|\le 1$ and all $f \in R{A}^s_{p,q}(\Rd)$. \\
{\rm (ii)}
Let $\sigma_p(d) < s < d/p$. There exists a constant $c>0$ such that
for all  $x$, $0< |x|<1$, there exists  a smooth radial function
$ f \in  R{A}^s_{p,q} (\Rd)$, $ \|\, f \, |R{A}^s_{p,q} (\Rd)\|=1$, s.t.
\be\label{eq-72b}
|x|^{\frac dp -s}\,  |f(x)| \ge c \, 
.
\ee
\end{T}

\noindent
In the next figure (Fig.~3) we summarize our knowledge.\\
{~}\\
\beginpicture

\setcoordinatesystem units <0.8cm,0.8cm>
\unitlength0.8cm

\put {\vector (0,1) {8} } [Bl] at 1 0
\put {\vector (1,0) {13} } [Bl] at 0 1

\thicklines
\plot 1 1 13 8 /
\plot 1 1 13 4 /

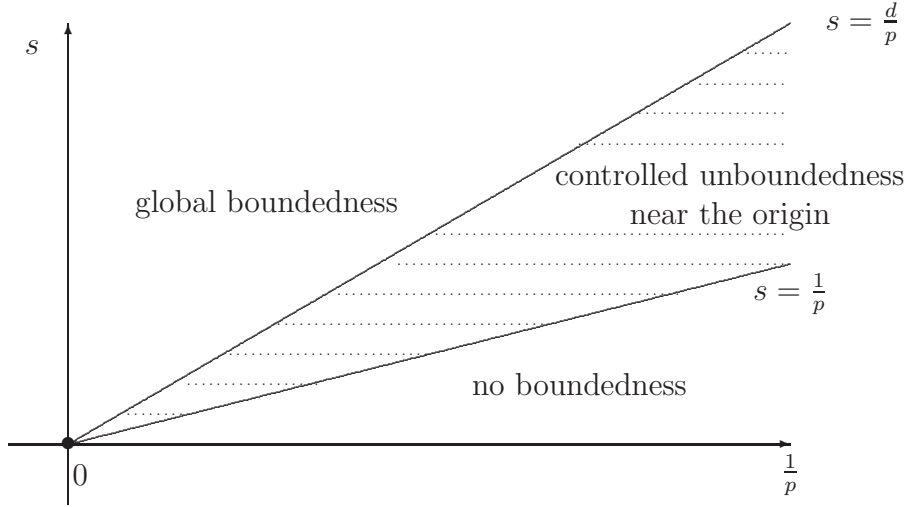
\put {Fig.~3: Unboundedness at the origin} at 6 -0.5

\put {$0$} at 1.2 0.5
\put {$\frac{1}{p}$} at 13 0.5
\put {$s$} at 0.4 7.6
\put {global boundedness} at 4.3 5

\put {controlled unboundedness} at 12 5.5
\put {near the origin} at 12 4.8

\put {no boundedness} at 9.5 2

\put {$\bullet$} at 1 1

\put {$s=\frac{d}{p}$} at 14.2 8
\put {$s = \frac{1}{p}$} at 13 3.5

{\setdots <0.1cm>

\plot 7 4.5 13 4.5 /
\plot 6.5 4 13 4 /
\plot 5.5 3.5 11.2 3.5 /
\plot 4.5 3 9 3 /
\plot 3.7 2.5 7 2.5 /
\plot 3 2 5.2 2 /
\plot 2 1.5 3 1.5 /
\plot 9.5 6 13 6 /
\plot 10.5 6.5 13 6.5 /
\plot 11.5 7 13 7 /
\plot 12.3 7.5 13 7.5 /
}
\endpicture
\hfill

\begin{Rem}{\rm
(i) In case of  $R{B}^s_{p,\infty} (\Rd)$
we have a function which realizes the extremal behaviour for all $|x|<1$ simultaneously.
It is well-known, see e.g. \cite[Lem.~2.3.1/1]{RS}, that the function
\[
f(x):= \psi (x)\, |x|^{\frac{d}{p} -s}\, , \qquad x \in \Rd\, ,
\]
belongs to  $R{B}^s_{p,\infty} (\Rd)$, as long as
$s> \sigma_p(d)$. This function does not belong to
$R{B}^s_{p,q} (\Rd)$, $q<\infty$. Since it is also not contained in
$R{F}^s_{p,q} (\Rd)$, $0 < q\le \infty$
we conclude that in these  cases  there is no function, which realizes this upper
bound for all $x$    simultaneously.
In these cases the function $f$ in (\ref{eq-72b}) has to depend on $x$.\\
(ii) These estimates do not change by switching to the larger homogeneous spaces
$R\dot{A}^s_{p,q} (\Rd)$ of Besov and Lizorkin-Triebel type.
In case of $R\dot{H}^s (\Rd) = R\dot{F}^s_{p,2}(\Rd) $ this has been observed in a recent paper by Cho and Ozawa \cite{CO}, see also Ni \cite{Ni1},  Rother \cite{Ro} and Kuzin, Pohozaev \cite[8.1]{KP}. The general case is treated in \cite{SSV2}.
\\
(iii) Estimates as in (\ref{eq-72}) have been investigated also in \cite{SS}.
However, there we proved a much weaker result only. 
We had overlooked the $s$-dependence of the power of $|x|$.
\\
(iv) In the literature one can find various types of further inequalities
for radial functions.
Many times preference is given to a homogeneous context, see
the inequalities (\ref{strauss}) and (\ref{ws-01}) as examples.
Then one has to deal with the behaviour at infinity and around the origin
simultaneously. That would be not appropriate in the context of inhomogeneous spaces. Inequalities like  (\ref{strauss}) and (\ref{ws-01}) will be investigated systematically in \cite{SSV2}. However, let us refer to \cite{St}, \cite{Li},
\cite{Ni1}, \cite{Ro}, \cite[8.1]{KP} and \cite{CO} for results in this direction.
Sometimes also decay estimates are proved by replacing on the right-hand side
the norm in the space $\aspq (\Rd)$ ($\dot{A}^s_{p,q}(\Rd)$) by products of norms,
e.g., $\|\, f |L_p (\Rd)\|^{1-\Theta}\, \| \, f \, |\dot{A}^s_{p,q}(\Rd)\|^{\Theta}$
for some $\Theta \in (0,1)$, see \cite{Li},
\cite{Ni1}, \cite{Ro} and \cite{CO}. Here we will not deal with those modifications (improvements).
}
\end{Rem}

Finally we have to investigate $s \le 1/p$ and $(s,p,q) \not\in U(A)$.

\begin{Lem}\label{lim2}
Let
$d\ge 2$, $0<  p  < \infty  $ and  $0 < q \le \infty$.
Suppose $(s,p,q) \not\in U(A)$ and $\sigma_{p}(d) < 1/p$. Moreover let  $\sigma_{q}(d) < 1/p$ in the $F$-case.
Then there exists  a  radial function
$ f \in  R{A}^s_{p,q} (\Rd)$, $ \|\, f \, |R{A}^s_{p,q} (\Rd)\|=1$,
and a sequence $(x^j)_j \subset \Rd\setminus \{0\}$ s.t.
 $\lim_{j \to \infty} |x^j| = 0$ and $f$ is unbounded in a neighborhood of
 all $x^j$.
\end{Lem}


\subsubsection{The behaviour of radial functions near
the origin -- borderline cases}
\label{border2}


Now we turn to the remaining limiting situation. We
shall show  that there is also controlled unboundedness near the origin if $s=d/p$ and
$ R{A}^{d/p}_{p,q}(\Rd) \not\subset  L_\infty (\Rd)$.

\begin{T}\label{lim1}
Let
$d\ge 2$, $0<  p  < \infty  $,  $0 < q \le \infty$, and suppose $s=d/p$.
\\
{\rm (i)} Let $1 < q \le \infty$.
Then there exists  constant $c$ s.t.
\be\label{eq-72l}
(- \log |x|)^{-1/q'}\,  |f(x)| \le c \, \| \, f \, |{B}^{d/p}_{p,q} (\Rd)\| \,
\ee
holds for all $0< |x|\le 1/2$ and all $f \in R{B}^{d/p}_{p,q}(\Rd)$. \\
{\rm (ii)}  Let $1 < p  < \infty$.
Then there exists  constant $c$ s.t.
\be\label{eq-72li}
(- \log |x|)^{-1/p'}\,  |f(x)| \le c \, \| \, f \, |{F}^{d/p}_{p,q} (\Rd)\| \,
\ee
holds for all $0< |x|\le 1/2$ and all $f \in R{F}^{d/p}_{p,q}(\Rd)$.
\end{T}

\begin{Rem}
{\rm
Comparing Lemma \ref{sharp} below and Theorem \ref{lim1}
we find the following.
For the case $q=\infty$ in Theorem \ref{lim1}(i)
the function $f_{1,0}$, see (\ref{spezial}),
realizes the extremal behaviour. In all other cases there remains a gap of order $\log \log$ to some power.
}
\end{Rem}


\section{Traces of radial subspaces -- proofs}
\label{trace}


The main aim of this section is to prove Theorem
\ref{main5}.  It expresses the fact that  all information about a
radial function
is contained in its trace onto a straight line through the  origin.
However, a few things more will be done here.
For later use one subsection is devoted to the study of interpolation of
radial subspaces (Subsection \ref{int}) and another one is devoted to the study of certain test functions (Subsection \ref{testfunctions}).


\subsection{Interpolation of radial subspaces}
\label{int}


We mention two different results here, one with respect
to the complex method and one with respect to the real method of interpolation.


\subsubsection{Complex  interpolation of radial subspaces}


In \cite{LSR} one of the authors has proved that in case $p,q\ge 1$
the spaces  $RB^{s}_{p,q} (\Rd) $ ($RF^{s}_{p,q} (\Rd) $) are  complemented subspaces of $\bspq (\Rd)$ ($\fspq (\Rd)$).
By means of the method of retraction and coretraction, see, e.g.,
Theorem 1.2.4 in \cite{TrI}, this allows to transfer the interpolation formulas
for the original spaces $B^{s}_{p,q} (\Rd) $ ($F^{s}_{p,q} (\Rd) $)
to its radial subspaces.
However, we prefer to quote a slightly more general result, proved in
\cite{SSV3}, concerning the complex method.
It is based on the results on complex interpolation for Lizorkin-Triebel spaces from \cite{FJ1} and uses the method of \cite{MM} for an extension to the quasi-Banach space case.

\begin{Prop}\label{interpol2}
Let  $0 <  p_0 ,  p_1 \le \infty$, $0 <  q_0 ,\,  q_1 \le \infty$, $s_0, s_1 \in \R$, and $0 < \Theta < 1$.
Define $s:= (1-\Theta)\, s_0 + \Theta \, s_1$,
\[
\frac 1p := \frac{1-\Theta}{p_0} + \frac{\Theta}{p_1}
\qquad \mbox{and}\qquad \frac 1q := \frac{1-\Theta}{q_0} +
\frac{\Theta}{q_1} \, .
\]
{\rm (i)} Let $\max (p_0,q_0) <\infty$.
Then we have
\[
RB^s_{p,q} (\Rd)  =   \Big[RB^{s_0}_{p_0,q_0} (\Rd),
RB^{s_1}_{p_1,q_1} (\Rd)\Big]_{\Theta} \, .
\]
{\rm (ii)} Let $p_1 <\infty$ and $\min (q_0,q_1) <\infty$.
Then we have
\[
RF^s_{p,q} (\Rd) =   \Big[RF^{s_0}_{p_0,q_0} (\Rd),
RF^{s_1}_{p_1,q_1} (\Rd)\Big]_{\Theta} \, .
\]
\end{Prop}


\subsubsection{Real Interpolation of  radial subspaces}


For later use we also formulate a result with respect to the real method
of interpolation.

\begin{Prop}\label{interpol1}
Let $d\ge 1$, $1 \le q, q_0, q_1 \le \infty$, $s_0, s_1 \in \R$, $s_0 \neq s_1$,
and $0 < \Theta < 1$.
\\
{\rm (i)} Let $1 \le p \le \infty$.
Then, with $s:= (1-\Theta)\, s_0 + \Theta \, s_1$, we have
\[
RB^s_{p,q} (\Rd)  =   \Big(RB^{s_0}_{p,q_0} (\Rd),
RB^{s_1}_{p,q_1} (\Rd)\Big)_{\Theta,q} \, .
\]
{\rm (ii)} Let $1\le p < \infty$.
Then, with $s:= (1-\Theta)\, s_0 + \Theta \, s_1$, we have
\[
RB^s_{p,q} (\Rd) =   \Big(RF^{s_0}_{p,q_0} (\Rd),
RF^{s_1}_{p,q_1} (\Rd)\Big)_{\Theta,q} \, .
\]
\end{Prop}

\noindent
\bpr
As mentioned above, the spaces
$RB^{s}_{p,q} (\Rd) $ ($RF^{s}_{p,q} (\Rd) $) are  complemented subspaces of $\bspq (\Rd)$ ($\fspq (\Rd)$), see \cite{LSR}.
Using the method of retraction and coretraction, see \cite[1.2.4]{Tr1},
the above statements are consequences of the corresponding formulas without $R$, see e.g.
\cite[2.4.2]{Tr1}.
\epr


 \subsection{Proofs of the statements in Subsection \ref{main11}}
\label{sub11}



\subsection*{Proof of Theorem \ref{anfang}}


\noindent
{\em Step 1.} Let $m \in \N_0$.
If $f\in RC^{m}(\R^d)$, then, of course, $f_0 = \tr f \in RC^{m}(\R)$
and  the inequality  
\beq\label{eq:trf}
\|\, \tr f \, |C^m (\R)\| \le \| \, f \, |C^m (\Rd)\|
\eeq
follows immediately.
\\
{\em Step 2.}  Let $g \in RC^{m}(\R)$ for some $m \in \N$.\\
{\em Substep 2.1.} We first deal with the regularity of $f: = \ext g$.
Here we make use of the following elementary observation.
Let $U(x^0)$ denote the neigborhoud of a point $x^0 \in \Rd$. 
If a function $v$ belongs to $C^m (U(x^0)\setminus \{x^0\})$, such that
$v$ is continuous in $x^0$ and the limit
\begin{equation}\label{limit}
\lim_{x \to x^0} \, D^\alpha v (x)
\end{equation}
exists for all $\alpha$, $|\alpha|\le m$, and all $j$, $j=1, \ldots \, , d$,
then $v \in C^m (U(x^0))$.
This is a consequence of the mean-value theorem.
By means of this argument it remains to prove the existence of 
these limits for $D^\alpha f$.
This turns out to be most simple using the extra condition
\be\label{eq:rx4}
g(0)=g'(0)=\dots=g^{(m)}(0)=0 \, .
\ee
This clearly implies for $0 \le \ell \le m$
\be\label{erg}
g^{(\ell)}(t)= o(|t|^{m-\ell}) \qquad \mbox{if}\quad t \to 0\, .
\ee
Furthermore, we shall use, for every $n\in\N_0$ there is a constant $c>0$ such that the function $r= r(x)$ satisfies
\be\label{eq:rx}
|D^\beta r(x)|\le c \, r(x)^{1-|\beta|}
\ee
for all $\beta\in\N_0^d$, $|\beta|\le n$,  and all $x\in\Rd \setminus\{0\}$.
Hence, using chain rule and the  estimates \reff{eq:rx},  (\ref{erg}) we find
 \beq\label{erg1}
|\,(D^\alpha f)(x)\, | & \lsim &  \sum_{\ell=1}^{|\alpha|}
|g^{(\ell)}(r)| \,   \sum_{\beta^1 + \ldots + \beta^\ell=\alpha}
|D^{\beta^1} r(x)| \, \ldots \, |D^{\beta^\ell} r(x)|
\nonumber
\\
& \lsim &  \sum_{\ell=1}^{|\alpha|}
o( r^{m-\ell}) \,   r^{\ell - |\alpha|}
= o( r^{m-|\alpha|}) , \qquad r \downarrow 0
\, ,
 \eeq
where  $|\alpha|\le m$. This proves the existence of the limit in (\ref{limit}) and therefore, $D^\alpha f$, $|\alpha|\le m$, is
continuous everywhere.
\\ 
Next, we wish to remove the restriction (\ref{eq:rx4}).
Suppose that $m$ is even. Hence
\[
g'(0)=g'''(0)=\, \dots \, =g^{(m-1)}(0)=0,
\]
but $g(0), g''(0),\dots, g^{(m)}(0)$ can  be arbitrary.
Let $\psi_0 = \tr \psi$.
We introduce the function
\[
h(t): = g(t)-g_1(t)\, \psi_0(t)\, , \qquad t \in \R\, ,
\]
where
\[
g_1(t):= g(0)+\frac{g''(0)}{2!}t^2+\dots+\frac{g^{(m)}(0)}{(m)!}t^{m},
\qquad t\in\R \, .
\]
The extension of $g_1\, \psi_0$ is a radial function with
compact support and continuous derivatives of arbitrary order.
The function $h$ satisfies   (\ref{eq:rx4}) and this implies 
that $\ext h$ has continuous derivatives up to order $m$ in  $\Rd$.
and therefore, also $f$ has continuous derivatives up to order $m$ in  $\Rd$.
The uniform continuity of the derivatives of $f$ in $|x|\ge 1$ 
follows immediately from the chain rule, the uniform continuity in $|x|\le 1$ 
is obvious.
For odd $m$ the proof is similar.
\\
{\em Substep 2.2}. It remains to estimate the norm of $f$ in $RC^m(\Rd)$.
Again we deal with $m$ being even first.
Using the same notation as in Step 2.1 we find on the one hand  the obvious estimate
\beqq
\|\, \ext  g_1\, \psi_0\, | C^{m}(\R^d)||
& \le & \sum_{j=0}^{m/2}\frac{|g^{(2j)}(0)|}{(2j)!} \, \|\, |x|^{2j} \, \psi(x) \, |C^{m}(\R^d) \|
\\
& \lsim &  \|\, g\, |C^{m}(\R)\| \, ,
\eeqq
and on the other hand we can estimate 
$\| \, \ext h \, | C^{m}(\R^d)\|$ as in (\ref{erg1}). 
This leads to 
\begin{align*}
\| \, \ext h \, | C^{m}(\R^d)\| & \lsim \, \| \, h\, |C^{m}(\R)\|
\\
& \lsim \,  \|\, g\, |C^{m}(\R)\| + \|\, g_1\, \psi_0\, |C^{m}(\R)\|
\\
&\lsim \,  \|\, g\, |C^{m}(\R)\|\, .
\end{align*}
By means of  $\ext g = \ext h + \ext (g_1\, \psi_0)\in RC^{(m)}(\R^d)$
this proves 
\[
\|\, \ext g \, | C^{m}(\R^d)\| \lsim\,  \|\, g\, | C^{m}(\R)\|\, .
\]
As above, for odd $m$ the proof is similar.
\epr

\subsection*{Proof of Theorem \ref{infty}}

For $\tr \in \cl (B^{s}_{\infty,q} (\Rd), B^{s}_{\infty,q} (\R))$
we refer to \cite[2.7.2]{Tr1}.
This immediately gives $\tr \in \cl (RB^{s}_{\infty,q} (\Rd), RB^{s}_{\infty,q} (\R))$.
Concerning $\ext$ we argue by using real interpolation.
Observe, that  $\ext \in \cl  (RC^m(\R), RC^m(\Rd))$ for all $m\in \N_0$,
see Theorem \ref{anfang}. From the interpolation property of the real interpolation method we derive
\[
\ext \in  \, \cl \Big( (RC^{m} (\R), RC(\R))_{\Theta,q},
(RC^m (\Rd), RC (\Rd))_{\Theta,q}\Big)
\]
for all $\Theta \in (0,1)$ and all $q \in (0,\infty]$.
Using Proposition  \ref{interpol1}, simple monotonicity properties of the real method and 
\[
B^n_{\infty,1}(\Rd) \hookrightarrow C^n (\Rd) \hookrightarrow B^n_{\infty,\infty}(\Rd)\, , \qquad n \in \N_0\, ,
\]
the claim follows.
\epr


 \subsection{Proofs of the assertions in Subsection  \ref{main12}}
\label{main11b}



 \subsubsection{Proof of Lemma  \ref{simple}}
\label{simpleb}


\noindent
Recall, for $f \in RL_p (\R)$ we  have
\[
\int_{\Rd} |f(x)|^p \, dx = 2\, \frac{\pi^{d/2}}{\Gamma(d/2)} \, \int_0^\infty |f_0(r)|^p \, r^{d-1}\, dr\, .
\]
Using
\[
 \int_0^\infty |f_0(r)|^p \, r^{d-1}\, dr = \lim_{\varepsilon \downarrow 0}
 \int_\varepsilon^\infty |f_0(r)|^p \, r^{d-1}\, dr\, ,
\]
which implies the density of the test functions in $ L_p ([0, \infty), r^{d-1})$,
we can read this formula also from the other side, it means
\[
\int_{\Rd} | \ext g (x)|^p \, dx = 2\, \frac{\pi^{d/2}}{\Gamma(d/2)} \, \int_0^\infty |g(r)|^p \, r^{d-1}\, dr
\]
for all $g \in  L_p ([0, \infty), r^{d-1})$. This proves (i).
Part (ii) is obvious.
\epr


 \subsubsection{Characterizations of radial subspaces by atoms}
\label{rsa}


As mentioned above our proof of the trace theorem
relies on atomic decompositions of radial distributions on $\Rd$.
We recall our characterizations of $R\aspq (\Rd)$ from \cite{SS}, see also
 \cite{KLSS}.\\
In this paper we shall consider two different versions of atoms.
They are not related to each other.
We hope that it will be always clear from the context
with which type of atoms we are working.
For the following definition of an atom we refer to \cite{FJ1} or \cite[3.2.2]{Tr2}.
For an open set $Q$ and $r>0$ we put
$r\, Q= \{x \in \Rd \, : \: \dist (x,Q)<r\}$. Observe that $Q$ is always
a subset of $r \, Q$ whatever $r$ is.

\begin{Def}
\label{atomunw}
Let $s \in \R$ and let $0<  p \le \infty $.
Let $L$ and $M$ be integers such that
$L \ge 0$ and $M \ge -1$. Let $Q \subset \Rd$ be an open connected set with $\diam Q =r$.
\\
{\rm (a)}
A smooth function $a(x)$ is called an $1_L$-atom centered in $Q$ if
\beqq
\supp a & \subset & \frac r 2 \, Q\, , \\
\sup_{y \in \Rd}\, |D^\alpha a(y)| & \le & 1 \, ,
\qquad |\alpha | \le L\, .
\eeqq
{\rm (b)}
A smooth function $a(x)$ is called an $(s,p)_{L,M}$-atom centered in $Q$ if
\beqq
\supp a & \subset & \frac r 2 \, Q\, , \\
\sup_{y \in \Rd}\, |D^\alpha a(y)| & \le &
r^{s-|\alpha|-\frac dp} \, , \qquad |\alpha | \le L\, , \\
\int_{\Rd} a(y)\, y^\alpha \, dy  & = & 0\, , 
\qquad
|\alpha | \le M \, .
\eeqq
\end{Def}

\begin{Rem}{\rm
If $M=-1$, then the interpretation is that no  moment condition is required.}
\end{Rem}

\noindent
In \cite{SS} and \cite{KLSS}
we constructed a regular sequence of  coverings
with certain special properties which we now recall.
Consider the annuli (balls if $k=0$)
\[
P_{j,k} := \bigg\{ x \in \Rd \, : \quad k \, 2^{-j} \le |x| <
(k+1)\, 2^{-j}\bigg\}\, , \qquad j =0,1,\ldots \, , \quad k=0,1,
\ldots \, .
\]
Then there  is a sequence
$(\Omega_j)_{j=0}^\infty= ((\Omega_{j,k,\ell})_{k,\ell}
)_{j=0}^\infty$
of coverings of $\Rd$ such that
\begin{itemize}
\item[(a)]
all $\Omega_{j,k,\ell}$ are balls with center in $x_{j,k,\ell}$ s.t.
$x_{j,0,1}=0$ and $|x_{j,k,\ell}|= 2^{-j} (k+1/2)$ if $k \ge 1$;
\item[(b)]
$\diam \Omega_{j,k,\ell} = 12 \, \cdot \,  2^{-j}$ for all $k$ and all $\ell$;
\item[(c)]
$P_{j,k} \subset  \bigcup\limits_{\ell=1}^{C(d,k)}
\Omega_{j,k,\ell}\, , \qquad j=0,1,\ldots \, , \quad k=0,1,
\ldots \, ,$ where the numbers $C(d,k)$ satisfy the relations
$C(d,k)\le  \, (2k+1)^{d-1}$, $\, \, C(d,0)=1$.
\item[(d)]
the sums $\sum\limits_{k=0}^\infty \sum\limits_{\ell=1}^{C(d,k)}
{\cal X}_{j,k,\ell}(x) $  are uniformly bounded in $x\in \Rd$ and
$j =0,1,\ldots $ (here ${\cal X}_{j,k,\ell}$ denotes the
characteristic function of $\Omega_{j,k,\ell}$);
\item[(e)]
$\Omega_{j,k,\ell}= \{x\in \Rd\, : \: 2^j\, x \in
\Omega_{0,k,\ell}\}$ for all  $j, k $ and $\ell$;
\item[(f)] There exists a natural  number $K$ (independent of $j$ and $k$)
such that
\be\label{K}
\{(x_1,0,\ldots , 0)\, : \quad x_1 \in \R\} \: \cap  \:
\frac{\diam(\Omega_{j,k,\ell})}{2}\,  \, \Omega_{j,k,\ell} = \emptyset\qquad \mbox{if}
\quad \ell >K
\ee
(with an appropriate enumeration).
\end{itemize}

\noindent
We collect some properties of  related atomic decompositions.
To do this it is convenient to introduce some  sequence spaces.

\begin{Def}
Let $0 <q \le \infty$.
\\
{\rm (i)}
If $0 < p \le \infty$, then we define
\[
b_{p,q,d} := \bigg\{
s= (s_{j,k})_{j,k}\, : \quad
\| \,s \, | b_{p,q,d}\| =
\left(\sum_{j=0}^\infty \left(
\sum_{k=0}^\infty (1+k)^{d-1}\,
|s_{j,k}|^p\right)^{q/p}
\right)^{1/q}<\infty \bigg\}\, .
\]
{\rm (ii)} By
$\widetilde{\chi}_{j,k}$  we denote the  characteristic function of
the set $P_{j,k}$.
If $0 <p< \infty$ we define
\beqq
f_{p,q,d} & := & \bigg\{ s= (s_{j,k})_{j,k}\, : \\
\nonumber
\\
\| \,s \, | f_{p,q,d}\| & = &
\left\| \bigg(
\sum_{j=0}^\infty \,
\sum_{k=0}^\infty  \, |s_{j,k}|^q \, 2^{\frac{jdq}{p}}\, \widetilde{\chi}_{j,k}(\cdot)
\bigg)^{1/q}|L_p(\Rd)\right\|\, <\infty \bigg\}\, .
\eeqq
\end{Def}

\begin{Rem}{\rm
Observe $b_{p,p,d} = f_{p,p,d}$ in the sense of equivalent quasi-norms.}
\end{Rem}

\noindent
Atoms have to satisfy moment and regularity conditions. With this respect  we
suppose
\be\label{LM1}
L \ge \max(0, [s]+1)\; ,\qquad  M\geq \max([\sigma_p(d) - s], -1)
\ee
in case of Besov spaces
and
\be\label{LM2}
L \ge \max(0, [s]+1)\; ,\qquad
M\geq \max([\sigma_{p,q}(d)- s], -1)
\ee
in case of Lizorkin-Triebel spaces. Under these restrictions the following assertions
are known to be true:
\\
\begin{itemize}
\item[(i)]
Each $f \in R\bspq (\Rd)$ ( $f \in R\fspq (\Rd)$) can be decomposed into
\be\label{atom1}
f= \sum_{j=0}^\infty \sum_{k=0}^\infty \sum_{\ell=1}^{C(d,k)}
s_{j,k}\, a_{j,k,\ell}
\qquad
(\, \mbox{convergence in}\quad \cs'(\Rd)\, ),
\ee
where the functions $a_{j,k,\ell}$ are $(s,p)_{L,M}$-atoms with respect to
$\Omega_{j,k,\ell}$ $(j\ge 1)$,
and the functions $a_{0,k,\ell}$ are $1_L$-atoms with respect to
$\Omega_{0,k,\ell}$.
\item[(ii)]
Any formal series  $\sum_{j=0}^\infty
\sum_{k=0}^\infty \sum_{\ell=1}^{C(d,k)}
s_{j,k}\, a_{j,k,\ell}$
converges in $\cs' (\Rd)$ with limit in $\bspq (\Rd)$ if the sequence
$s=(s_{j,k})_{j,k}$ belongs to $b_{p,q,d}$ and if
the $a_{j,k,\ell}$ are $(s,p)_{L,M}$-atoms with respect to
$\Omega_{j,k,\ell}$ $(j\ge 1)$,
and the $a_{0,k,\ell}$ are $1_L$-atoms with respect to
$\Omega_{0,k,\ell}$.
There exists a universal constant such that
\be\label{atom2}
 \| \,  \sum_{j=0}^\infty \sum_{k=0}^\infty \sum_{\ell=1}^{C(d,k)}
s_{j,k}\, a_{j,k,\ell} \, | \bspq (\Rd)\|
\le c \, \| \, s \, |b_{p,q,d}\|
\ee
holds for all sequences $s=(s_{j,k})_{j,k}$. \\
\item[(iii)]
There exists a constant $c$ such that for any $f \in R\bspq (\Rd)$
there exists an atomic decomposition as in (\ref{atom1}) satisfying
\be\label{atom3}
\|\, (s_{j,k})_{j,k} \, |b_{p,q,d}\| \le c \, \| \, f \, |
\bspq (\Rd)\| \, .
\ee
\item[(iv)]
The infimum on the left-hand side in (\ref{atom2})
with respect to all admissible representations
(\ref{atom1})
yields  an equivalent norm on $R\bspq (\Rd)$.
\item[(v)]
Any formal series  $\sum_{j=0}^\infty
\sum_{k=0}^\infty \sum_{\ell=1}^{C(d,k)}
s_{j,k}\, a_{j,k,\ell}$
converges in $\cs' (\Rd)$ with limit in  $\fspq (\Rd)$ if the sequence
$s=(s_{j,k})_{j,k}$ belongs to $f_{p,q,d}$
and if
the functions $a_{j,k,\ell}$ are $(s,p)_{L,M}$-atoms with respect to
$\Omega_{j,k,\ell}$ $(j\ge 1)$,
and the functions $a_{0,k,\ell}$ are $1_L$-atoms with respect to
$\Omega_{0,k,\ell}$.
There exists a universal constant such that
\be\label{atom2b}
 \| \,  \sum_{j=0}^\infty \sum_{k=0}^\infty \sum_{\ell=1}^{C(d,k)}
s_{j,k}\, a_{j,k,\ell} \, | \fspq (\Rd)\|
\le c \, \| \, s \, |f_{p,q,d}\|
\ee
holds for all sequences $s=(s_{j,k})_{j,k}$. \\
\item[(vi)]
There exists a constant $c$ such that for any $f \in R\fspq (\Rd)$
there exists an atomic decomposition as in (\ref{atom1}) satisfying
\be\label{atom3b}
\|\, (s_{j,k})_{j,k} \, |f_{p,q,d}\|  \le c \, \| \, f \, |
\fspq (\Rd)\| \, .
\ee
\item[(vii)]
The infimum on the left-hand side in (\ref{atom2b})
with respect to all admissible representations
(\ref{atom1})
yields  an equivalent norm on $R\fspq (\Rd)$.
Such  decompositions as in (\ref{atom3}) and (\ref{atom3b}) we  shall call
optimal.
\end{itemize}

\begin{Rem}{\rm
For proofs of all these facts (even with respect to more general decompositions of $\Rd$) we refer to \cite{SS} and \cite{LSR}.
A different approach to atomic decompositions of radial subspaces has been
given by Epperson and Frazier \cite{EF1}.
}
\end{Rem}


 \subsubsection{Proof of Theorem \ref{main5}}


\noindent
{\em Step 1.}
Let $f \in R\bspq (\Rd)$.  Then there exists an optimal atomic decomposition, i.e.
\beq\label{ws-3}
f & = &
\sum_{j=0}^\infty \sum_{k=0}^\infty \sum_{\ell=1}^{C(d,k)}
s_{j,k}\, a_{j,k,\ell} \\
\| \, f \, |
\bspq (\Rd)\|  & \asymp &
\|\, (s_{j,k})_{j,k} \, |b_{p,q,d}\|\, ,
\nonumber
\eeq
see (\ref{atom1}) - (\ref{atom3}).
Since $f$ is even we obtain
\be\label{ws-2}
f (x) =
\sum_{j=0}^\infty \sum_{k=0}^\infty \sum_{\ell=1}^{C(d,k)}
s_{j,k}\, \frac{a_{j,k,\ell} (x) + a_{j,k,\ell} (-x)}{2} \, .
\ee
We define
\[
g_{j,k,\ell} (t):= 2^{j(s-d/p)}\,
\Big(\tr  \frac{a_{j,k,\ell} (\, \cdot\,) + a_{j,k,\ell} (-\, \cdot\,)}{2}\Big)(t)\,, \qquad t \in \R\, ,
\]
and $d_{j,k}:= 2^{-j(s-d/p)}\, s_{j,k}$.
Of course, $a_{j,k,\ell} (\, \cdot\,) + a_{j,k,\ell} (-\, \cdot\,)$ is not a radial function.
But it is an even and continuous. So, $\tr$ means simply the restriction to the $x_1$-axis.
Clearly,
\[
f_N (x) := \sum_{j=0}^N \sum_{k=0}^N \sum_{\ell=1}^{C(d,k)}
s_{j,k}\, \frac{a_{j,k,\ell} (x) + a_{j,k,\ell} (-x)}{2}\, , \qquad x \in \Rd, \quad N \in \N\, ,
\]
is an even (not necessarily radial) function in $C^L (\Rd)$.
By means of property (f) of the particular coverings of $\Rd$, stated
in the previous subsection, we obtain
\[
\tr f_N \,  = \sum_{j=0}^N \sum_{k=0}^N \sum_{\ell=1}^{\min(C(d,k),K)}
d_{j,k}\, g_{j,k,\ell}
\]
(here $K$ is the natural number in (\ref{K})).
Furthermore
\[
\max_{0 \le n \le L}\, \sup_{t \in \R} \, |(g_{j,k,\ell})^{(n)}(t)|\le
12^{s-d/p}\,  2^{jn}\, .
\]
Obviously
\beqq
\| \, (s_{j,k})_{j,k} \, | b_{p,q,d}\| & = &
\left(\sum_{j=0}^\infty \left(
\sum_{k=0}^\infty (1+k)^{d-1}\,
|s_{j,k}|^p\right)^{q/p}  \right)^{1/q}
\\
& = &
\left(\sum_{j=0}^\infty  2^{j(s - \frac dp)q}
\left(
\sum_{k=0}^\infty (1+k)^{d-1}\,
|d_{j,k}|^p\right)^{q/p}  \right)^{1/q}\,.
\eeqq
This implies
\[
\| \, \tr f_N \, | TB^s_{p,q}(\R,L,d)\| \le K \, c \,  \| \,  f\, | B^s_{p,q}(\Rd)\|
\]
where $c$ and $K$ are independent of $f$ and $N$. \\
Next we comment on the convergence
of the sequences $(f_N)_N$ and $(\tr f_N)_N$.
Of course, $f_N$ converges in $\cs' (\Rd)$ to $f$.
For the investigation of the convergence of $(\tr f_N)_N$
we choose $s'$ such that  $s> s' > \sigma_p(d)$ and conclude with $N>M$
\beqq
\| \, \tr f_N  &&\hspace{-0.9cm} -   \tr f_M \, |TB^{s'}_{p,p} (\R,L,d)\|
 \lsim
\Big\| \, \sum_{j=M+1}^N \sum_{k=0}^N \sum_{\ell=1}^{\min (C(d,k),K)}
d_{j,k}\,  g_{j,k,\ell} \, \Big| TB^{s'}_{p,p}(\R)\Big\|
\\
&& \qquad + \qquad \Big\| \, \sum_{j=0}^M \sum_{k=M+1}^N \sum_{\ell=1}^{\min (C(d,k),K)}
d_{j,k}\,  g_{j,k,\ell} \, \Big| TB^{s'}_{p,p}(\R)\Big\|
\\
& \lsim &
 \Big(\sum_{j=M+1}^\infty   \sum_{k=0}^\infty
(1+ k)^{d-1} |2^{j(s'-s)} \, s_{j,k}|^p  \Big)^{1/p}
\\
&  & \qquad + \qquad \Big(\sum_{j=0}^M   \sum_{k=M+1}^\infty
(1+ k)^{d-1} |2^{j(s'-s)} \, s_{j,k}|^p  \Big)^{1/p}\, ,
\eeqq
by taking into account the different normalization of the atoms in
$RB^{s'}_{p,p}(\Rd)$ and in $RB^{s}_{p,q}(\Rd)$, respectively.
The right-hand side in the previous inequality tends to zero if $M$ tends to infinity since
$\| \, (s_{j,k})_{j,k}\, |b_{p,q,d}\|<\infty $.
Lemma \ref{simple} in combination with
$B^s_{p,q} (\Rd) \subset {{L_{\max(1,p)}(\Rd)}}$ implies the continuity of
$\tr : \,  RB^{s}_{p,1}(\Rd) \to  L_{\max(1,p)}(\R,t^{d-1}) $ as well as the existence of
$\tr f \in L_{\max(1,p)} (\R,|t|^{d-1})$. Consequently
\[ \lim_{N\to \infty} \tr f_N = \tr (\lim_{N\to \infty} f_N)= \tr f
\]
with convergence in $ L_{\max(1,p)} (\R,|t|^{d-1})$.
This proves that
$\tr $ maps $R\bspq (\Rd)$ into $T\bspq (\R,L,d)$ if
$L$ satisfies (\ref{LM1}). Observe, that $M$ can be chosen {{ $-1$ in (\ref{LM1})}}.
\\
{\em Step 2.}
The same type of arguments proves
 that
$\tr $ maps $R\fspq (\Rd)$ into $T\fspq (\R)$, in particular the convergence analysis
is the same. Furthermore, observe
\beqq
\bigg\| \bigg(
\sum_{j=0}^\infty \, & 2^{jsq} &\,
\sum_{k=0}^\infty  \, |s_{j,k}|^q \, {\chi}^\#_{j,k}(\cdot)
\bigg)^{1/q}|L_p(\R, t^{d-1})\bigg\|
\\
& = & c_d \,
\bigg\| \bigg(
\sum_{j=0}^\infty \, 2^{jsq}\,
\sum_{k=0}^\infty  \, |s_{j,k}|^q \, \widetilde{\chi}_{j,k}(\cdot)
\bigg)^{1/q}|L_p(\Rd)\bigg\|\, .
\eeqq
This proves that
$\tr $ maps $R\fspq (\Rd)$ into $T\fspq (\R,L,d)$ if $L$ satisfies
(\ref{LM2}) (again we use  $M=-1$). \\
{\em Step 3.} Properties of $\ext$.
Let $g$ be an even function with a decomposition as in (\ref{new})
and
\[
\| \, g \, |TB^s_{p,q} (\R,L,d)\| \asymp
\| \, (s_{j,k})\, |b^s_{p,q,d}\|\, .
\]
We define
\[
a_{j,k} (x) :=  g_{j,k}(|x|) \, , \qquad x\in \Rd \, .
\]
The functions $a_{j,k}$ are compactly supported, continuous, and  radial.
Obviously
\[
\supp a_{j,k} \subset \{ x: \quad 2^{-j}k-2^{-j-1}\le |x|\le
 2^{-j}(k+1) + 2^{-j-1}\} \, , \qquad k \in \N \, ,
\]
and
\[
\supp a_{j,0} \subset \{ x: \quad  |x|\le 3\, \cdot \,  2^{-j-1}\} \, .
\]
From Theorem \ref{anfang} we derive
\be
|D^\alpha a_{j,k}(x)|  \le  \| \, a_{j,k}\, |C^{|\alpha|}(\Rd)\|
\lsim  \, \| \, g_{j,k}\, |C^{|\alpha|}(\R)\| \lsim \, 2^{j|\alpha|}\, ,
\label{res4}
\ee
if $|\alpha|\le L$.
Here the constants behind $\lsim$ do not depend on $j,k$ and $g_{j,k}$.
We continue with an investigation of the sequence
\beq\label{extension}
h_N (x) := \sum_{j=0}^N  \sum_{k=0}^\infty
s_{j,k}\, a_{j,k}(x)\, , \qquad x \in \Rd\, , \quad N \in \N\, .
\eeq
Related to our decomposition $(\Omega_{j,k,\ell})_{j,k,\ell}$
of $\Rd$, see Subsection \ref{rsa}, there is a sequence of  decompositions of unity
$(\psi_{j,k,\ell})_{j,k,\ell}$, i.e.
\beq\label{res}
\sum_{k=0}^\infty \sum_{\ell=1}^{C(d,k)}
\psi_{j,k,\ell} (x)  & = & 1 \qquad
\mbox{for all} \quad x \in \Rd\, , \quad j=0,1,\ldots \,  ,
\\
\label{res2}
\supp \psi_{j,k,\ell} & \subset &  \Omega_{j,k,\ell} \, ,
\\
|D^\alpha \psi_{j,k,\ell}| & \le &  C_L\, 2^{j|\alpha|} \qquad
 |\alpha|\le L \, ,
\label{res3}
\eeq
see \cite{SS}. Hence
\beqq
h_{N} (x) & = & \sum_{j=0}^N
\sum_{k = 0}^\infty s_{j,k} \, a_{j,k}(x)
\, \left(
\sum_{m=0}^\infty \sum_{\ell=1}^{C(d,m)}
\psi_{j,m,\ell} (x) \right)
\\
& = & \!\!\!   \sum_{j=0}^N
\sum_{m = 0}^\infty  \,  \left(
\sum_{\ell=1}^{C(d,m)} \sum_{k=-7}^7 s_{j,m+k} \, a_{j,m+k}(x)
\psi_{j,m,\ell} (x) \right)\\
& = & \!\!\!   \sum_{j=0}^N
\sum_{m = 0}^\infty  \lambda_{j,m} \, \sum_{\ell=1}^{C(d,k+m)}
\, {e}_{j,m,\ell}(x)  \, .
\eeqq
where
\beqq
\lambda_{j,m} & := & 2^{j(s- \frac{d}{p})}\,
\max_{-1 \le k \le 1}\, |s_{j,m+k}|
\\
{e}_{j,m,\ell}(x) & := & 2^{-j(s-\frac{d}{p})} \,
t_{j,m}\, \sum_{k=-7}^7 s_{j,m+k} \, a_{j,m+k}(x) \,
\psi_{j,m,\ell} (x)
\\
t_{j,m}& : = &
\begin{cases}
1& \text{if} \quad \max\limits_{k=-7, \, \ldots \, , 7}\, |s_{j,m+k}|= 0, \\
& \\
(\max\limits_{k=-7, \, \ldots \,, 7}|s_{j,m+k}|)^{-1}  & \text{ otherwise .}
\end{cases}
\eeqq
We claim that the functions ${e}_{j,m,\ell}$ are $(s,p)_{L,-1}$-atoms
($1_{L}$-atoms if $j=0$) on $\Rd$ related to the covering
$(\Omega_{j,k,\ell})_{j,k,\ell}$ (up to a universal constant).
But this follows immediately from (\ref{res2}), (\ref{res3}), and
(\ref{res4}).
Finally we show that the sequence $\lambda = (\lambda_{j,m})_{j,m}$
belongs to $b_{p,q,d}$. The estimate
\beqq
\| \, \lambda \, |b_{p,q,d}\| & = &
\Big(\sum_{j=0}^N \Big(\sum_{m=0}^\infty (1+m)^{d-1} \,
|\lambda_{j,m}|^p\Big)^{q/p} \Big)^{1/q}
\\
 & \lsim  &
\Big(\sum_{j=0}^N 2^{j(s-\frac dp)q} \Big(\sum_{m=0}^\infty (1+m)^{d-1} \,
|s_{j,m}|^p\Big)^{q/p} \Big)^{1/q}
\eeqq
is obvious. Hence $\ext $ maps ${T}\bspq (\R,L,d)$ into $R\bspq (\Rd)$.
Here we need that the pair $(L,-1) $ satisfies (\ref{LM1}).
\\
{\em Step 4.} The proof of the F-case is similar.
Here we need that the pair $(L,-1) $ satisfies (\ref{LM2}).
The proof is complete.
\epr


\subsubsection{Proof of Theorem \ref{spur1}}


Since $RF^s_{p,q}(\Rd) \hookrightarrow RB^s_{p,\infty}(\Rd) $
it will be enough to deal with radial Besov spaces.\\
{\em Step 1.} Let $1 \le p <\infty$. Then $\sigma_p (d) = \sigma_p (1)=0$.
From $s>0$ we derive $B^s_{p,q}(\Rd) \subset L_p (\Rd)$.
Hence $f$ is a regular distribution.\\
{\em Step 2.} Let $0 <p<1$.
Since $0 \not\in \supp f$ there exists some $\varepsilon >0$ s.t.
the ball with radius $\varepsilon$ and centre in the origin
has an empty intersection with $\supp f$. Let $\lambda >0$.
Since $f$ is a regular distribution if, and only if
$f(\lambda \, \cdot \, )$ is a regular distribution
we may assume $\varepsilon =2$.
Let $\varphi \in RC^\infty (\Rd)$ be a function s.t. $\varphi (x)=1$
if $|x|\ge 2 $ and $\varphi (x)=0$
if $|x|\le 1$.
Again we shall work with an optimal atomic decomposition of $f \in R\bspq (\Rd)$,
see (\ref{ws-3}).
Obviously
\[
f=  f \varphi = \sum_{j=0}^\infty \sum_{k=0}^\infty \sum_{\ell=1}^{C(d,k)}
s_{j,k}\, (\varphi \, a_{j,k,\ell})\, .
\]
By checking the various support conditions we obtain
\[
f=  \sum_{j=0}^{3} s_{j,0}\, (\varphi \, a_{j,0,1})
+  \sum_{j=0}^\infty \sum_{k= \max (1, 2^j-9)}^\infty \sum_{\ell=1}^{C(d,k)}
s_{j,k}\, (\varphi \, a_{j,k,\ell})\, .
\]
We consider the following splitting
\beqq
f_1 (x) & := &   \sum_{j=0}^{3} s_{j,0}\, \varphi (x) \, \frac{a_{j,0,1}(x)+ a_{j,0,1}(-x)}{2}
\\
 f_2 (x) & := & \sum_{j=0}^\infty
\sum_{k=\max (1,2^j-9)}^\infty \sum_{\ell=1}^{C(d,k)}
s_{j,k}\, \varphi (x)\, \frac{a_{j,k,\ell} (x) + a_{j,k,\ell} (-x)}{2}\, .
\eeqq
Concerning the first part $f_1$  we observe that $\tr f_1$ is a compactly supported
even $C^L$ function.
Now we  concentrate on $f_2$.
Let
\[
f_{2,N} (x)  :=  \sum_{j=0}^N \sum_{k= \max (1,2^j-9)}^{2^N} \sum_{\ell=1}^{C(d,k)}
s_{j,k}\, \varphi (x)\, \frac{a_{j,k,\ell} (x) + a_{j,k,\ell} (-x)}{2}\, , \qquad N \in \N\, .
\]
We put
\[
g_{j,k,\ell} (t):= \tr  \Big(\varphi (\, \cdot\, )\,
\frac{
a_{j,k,\ell}(\, \cdot\, ) + a_{j,k,\ell}(-\,\cdot\, )}{2}\Big) (t)\,, \qquad t \in \R\, ,
\]
by using the same convention concerning $\tr$ as in Step 1 of the proof of Thm. \ref{main5}.
Since
\[
\sup_{t\in \R} \, |g_{j,k,\ell}^{(n)} (t)| \le c_\varphi \, (12 \,\cdot\,  2^{-j} )^{s-n-d/p}\, , \qquad 0 \le n\le L\, ,
\]
and $|\supp g_{j,k,\ell} |\lsim 2^{-j}$
we obtain for a natural number $m$
\beqq
 \int_m^{m+1} && \hspace{-1cm} |\tr f_{2,N} (t)|\,  dt
 \lsim  \sum_{j=0}^N \sum_{k=\max (1,2^j-9)}^{2^N} \sum_{\ell=1}^{\min (C(d,k),K)}
|s_{j,k}|\, \int_m^{m+1} |g_{j,k,\ell} (t)| \, dt
\\
& \lsim & \sum_{j=0}^N  2^{-j} \, 2^{-j (s-d/p)}
\sum_{k=2^jm-9}^{2^j(m+1)+6}
|s_{j,k}|\\
& \lsim &  m^{-(d-1)/p}
\sum_{j=0}^N  2^{-j} \, 2^{-j (s-d/p)} \, 2^{-j(d-1)/p}\,
\Big(\sum_{k=2^jm-9}^{2^j(m+1)+6} (1+k)^{d-1}
|s_{j,k}|^p\Big)^{1/p}\, .
\eeqq
Hence
\beqq
\|\, \tr f_N \, |L_1 (\R)\| & = & 2 \,
\int_1^{\infty} |\tr f_N (t)|\, dt
\lsim \| \, (s_{j,k})_{j,k}\, |b_{p,\infty,d}\|\,
\sum_{m=1}^\infty m^{-(d-1)/p}
\\
& \lsim & \| \, f \, |\bspq (\Rd)\|
\eeqq
since $s>0$ and $0< p<1$. Let $M \le N$. Then
the same type of argument yields
\beqq
\int_{m}^{m+1} |\tr f_{2,N} (t)  &- &\tr f_{2,M} (t)|\, dt
\\
& \lsim &
 \sum_{j=M+1}^N \sum_{k=\max (1,2^j-9)}^{2^N} \sum_{\ell=1}^{\min (C(d,k),K)}
|s_{j,k}|\, \int_m^{m+1} |g_{j,k,\ell} (t)| \, dt
\\
&& + \quad
 \sum_{j=0}^N \sum_{k=\max (2^M, 2^jm-9)}^{2^N} \sum_{\ell=1}^{\min (C(d,k),K)}
|s_{j,k}|\, \int_m^{m+1} |g_{j,k,\ell} (t)| \, dt
\\
& \lsim &
2^{-Ms} \sup_{j=M+1, \ldots \,  }  {{2^{j\frac{d-1}{p}}}}
\sum_{k=2^jm-9}^{2^j(m+1)+6}
|s_{j,k}|
\\
 &&  \quad + \quad  \sum_{j=0}^N  2^{-j} \, 2^{-j (s-d/p)}
\sum_{k= \max (2^M, 2^jm-9)}^{2^j(m+1)+6}
|s_{j,k}|
\eeqq
\beqq
\qquad
& \lsim &  m^{-(d-1)/p}\, \Big(2^{-Ms} \, \| \, (s_{j,k})_{j,k}\, |b_{p,\infty,d}\|
\\
&& \qquad + \quad
\sup_{j=0,1, \ldots }  \Big(\sum_{k=\max (2^M, 2^jm-9)}^{2^j(m+1)+6} (1+k)^{d-1}
|s_{j,k}|^p\Big)^{1/p} \Big)\, .
\eeqq
Since
\[
\lim_{M\to \infty}\,
\sup_{j=0,1, \ldots }  \Big(\sum_{k=\max (2^M, 2^jm-9)}^{2^j(m+1)+6} (1+k)^{d-1}
|s_{j,k}|^p\Big)^{1/p} = 0
\]
for all $m \in \N$ we conclude
\[
\|\, \tr f_{2,N} -\tr f_{2,M}\, |L_1 (\R)\| \longrightarrow  0 \qquad \mbox{if} \quad M \to \infty\, .
\]
Hence
\[
\sum_{j=0}^\infty \sum_{k= \max (1,2^j-9)}^{\infty} \sum_{\ell=1}^{\min(C(d,k),K)}
s_{j,k}\, g_{j,k,\ell}   \in L_1 (\R)\, .
\]
Let $\theta \in \Rd$, $|\theta|=1$. We denote by
$\Tr_\theta$ the restriction of a continuous function to the line $\Theta:= \{ t\, \theta: \quad t \in \R\}$.
Now we repeat, what we have done with respect to the $x_1$-axis, for such a line.
As the outcome we obtain
\[
\Tr_\theta \Big( \sum_{j=0}^{3} s_{j,0}\, (\varphi \, a_{j,0,1})
+  \sum_{j=0}^N \sum_{k= \max (1, 2^j-9)}^{2^N} \sum_{\ell=1}^{C(d,k)}
s_{j,k}\, (\varphi \, a_{j,k,\ell})\Big)\, , \qquad N \in \N\, ,
\]
is a Cauchy sequence in $L_1 (\Theta)$ and the limit satisfies
\beqq
\Big\| \Tr_\theta \Big( \sum_{j=0}^{3} s_{j,0}\, (\varphi \, a_{j,0,1})
 & + &  \sum_{j=0}^\infty \sum_{k= \max (1, 2^j-9)}^{\infty} \sum_{\ell=1}^{C(d,k)}
s_{j,k}\, (\varphi \, a_{j,k,\ell})\Big)\, \Big|L_1 (\Theta)\Big\|
\\
& \lsim &
\| \, f \, |B^s_{p,q}(\Rd)\|
\eeqq
with a constant independent of $\theta$
(of course, here, by a slight abuse of notation, $\Tr_\theta$ denotes the continuous extension
of the previously defined mapping). Using spherical coordinates this yields
\begin{eqnarray*}
\int_{\Rd} |f(x)|\, dx & = &  \int_{|\theta|=1} \int_0^\infty \, |f(t\, \theta)|\, dt \,  d\theta\\
&\lsim & \| \, f \, |B^s_{p,q}(\Rd)\| \, .
\end{eqnarray*}
But this means $f$ is a regular distribution.
\epr

\begin{Rem}\rm
We have proved a bit more than stated.
Under the given restrictions
the pointwise trace $\tr f$ of a distribution $f \in R\bspq(\Rd)$,
$0 \not\in \supp f$, makes sense and
belongs to $L_1(\R)$.
\end{Rem}


\subsubsection{Proof of Remark \ref{smean}}
\label{Daubechies}


We shall argue by using the wavelet characterization of
$B^{\frac 1p - 1}_{p,\infty}(\Rd)$, see, e.g.,  \cite[Thm.~1.20]{Tr08}.
Let $\phi$ denote an appropriate univariate scaling function
and $\Psi$ an associated Daubechies wavelet of sufficiently high order.
The tensor product ansatz yields $(d-1)$ generators $\Psi_1, \ldots \, , \Psi_{2^d-1}$ for the
wavelet basis in $L_2 (\Rd)$. Let
$\Phi$ denote the $d$-fold tensor product of the univariate scaling function.
We shall use the abbreviations
\[
\Phi_{k} (x) := \Phi (x-k)\, , \qquad k \in \Zd\, ,
\]
and
\[
\Psi_{i,j,k}(x):= 2^{jd/2}\, \Psi_{i} (2^j x -k)\, , \qquad k
\in \Zd\, , \quad j \in \N_0\,, \quad i=1, \ldots \, , 2^d-1\,  .
\]
An equivalent norm in $B^{\frac 1p - 1}_{p,\infty}(\Rd)$
is given by
\[
\| \, f \, |B^{\frac 1p - 1}_{p,\infty}(\Rd)\| =
\Big(\sum_{k \in \Zd} \,|\langle  f, \Phi_{k} \rangle |^p
\Big)^{1/p} +
\sup_{j=0,1, \ldots }\,  2^{j(\frac 1p-1 + d(\frac 12 - \frac 1p))}\, \Big(\sum_{i=1}^{2^d-1}
\sum_{k \in \Zd} \,|\langle  f, \Psi_{i,j,k} \rangle |^p
\Big)^{1/p} \, .
\]
Daubechies wavelets have compact support. This implies
\[
\supp \Psi_{i,j,k} \subset C \, \{x\in \Rd: \: 2^{-j}(k_\ell-1) \le x_\ell \le 2^{-j}(k_\ell+1)\, , \: \ell=1, \ldots \, d\}
\]
and
\[
\supp \Phi_{k} \subset C \, \{x\in \Rd: \: (k_\ell-1) \le x_\ell \le (k_\ell+1)\, , \: \ell=1, \ldots \, d\}
\]
for an appropriate $C>1$. By employing these relations we conclude that
for fixed $j$ the cardinality of the set of those functions $\Psi_{i,j,k}$,
which do not vanish identically on $|x|=1$ is $\lsim 2^{j(d-1)}$.
There is the general estimate
\[
|\langle  f, \Psi_{i,j,k} \rangle |=\Big| \int_{|x|=1}
2^{jd/2}\, \Psi_i (2^jx-k)\, dx\Big|
\lsim 2^{jd/2} \, 2^{-j(d-1)}\, ,
\]
by using the information on the size of the support.
Inserting this we find
\beqq
 2^{j(\frac 1p-1 + d(\frac 12 - \frac 1p))}\, \Big(\sum_{i=1}^{2^d-1}
\sum_{k \in \Zd} \,|\langle  f, \Psi_{i,j,k} \rangle |^p
\Big)^{1/p}
&\lsim &
 2^{j(\frac 1p-1 + d(\frac 12 - \frac 1p))}\, 2^{j(d-1)/p}\,
2^{jd/2} \, 2^{-j(d-1)}
\\
& \lsim & 1 \, .
\eeqq
This proves the claim.
\epr


\subsubsection{Proof of Theorem \ref{spur0}}


From Thm. \ref{spur1} we already know that for $f \in R\aspq (\Rd)$,
$0 \not\in \supp f$,
the trace $\tr f$ makes sense and that
$\tr f \in L_1 (\R)$.
\\
{\em Step 1.} Let $f \in R\bspq (\Rd)$.
Since $0 \not\in \supp f$ there exists some $\varepsilon >0$ s.t.
the ball with radius $\varepsilon$ and centre in the origin
has an empty intersection with $\supp f$.
Without loss of generality we assume $\varepsilon <1$.
Let $\varphi \in RC^\infty (\Rd)$ be a function s.t. $\varphi (x)=1$
if $|x|\ge \varepsilon $ and $\varphi (x)=0$
if $|x|\le \varepsilon/2$.
Again we shall work with an optimal atomic decomposition of $f$,
see (\ref{ws-3}).
It follows
\[
f=  \sum_{j=0}^{m} s_{j,0}\, (\varphi \, a_{j,0,1})
+  \sum_{j=0}^\infty \sum_{k=k_j}^\infty \sum_{\ell=1}^{C(d,k)}
s_{j,k}\, (\varphi \, a_{j,k,\ell})
\]
where
\[
m := 1 + [\log_2 (18 \, \varepsilon^{-1})]\qquad \mbox{and}\qquad
k_j := \max(1, [2^{j-1}\, \varepsilon] - 10)\, .
\]
As in the previous proof we introduce the splitting $f= f_1 + f_2$,
where
\[
f_1 (x):= \sum_{j=0}^{m} s_{j,0}\, \varphi (x) \, \frac{a_{j,0,1}(x)+ a_{j,0,1}(-x)}{2}\, .
\]
Obviously, $\tr f_1$ is a compactly supported even $C^L$ function.
Let
\[
f_{2,N} (x)  :=  \sum_{j=0}^N \sum_{k= k_j}^{2^N} \sum_{\ell=1}^{C(d,k)}
s_{j,k}\, \varphi (x)\, \frac{a_{j,k,\ell} (x) + a_{j,k,\ell} (-x)}{2}\, , \qquad N \in \N\, .
\]
As above we use the notation
\[
g_{j,k,\ell} (t):= \tr  \Big(\varphi (\, \cdot\, )\,
\frac{
a_{j,k,\ell}(\, \cdot\, ) + a_{j,k,\ell}(-\,\cdot\, )}{2}\Big) (t)\,, \qquad t \in \R\, .
\]
Hence
\[
\tr f_{2,N} (t) = \sum_{j=0}^N \sum_{k= k_j}^{2^N} \sum_{\ell=1}^{\min (C(d,k),K)}
s_{j,k}\, g_{j,k,\ell} (t)\,
\]
Since
\[
\sup_{t \in \R}\, |g_{j,k,\ell}^{(n)} (t)| \le c_\varphi \, (12 \, \cdot \,  2^{-j} )^{s-n-d/p}
= c_\varphi \,  (12 \, \cdot \,  2^{-j})^{-(d-1)/p} \, (12 \, \cdot \,  2^{-j} )^{s-n-1/p}
\]
the functions $2^{-j(d-1)/p} 12^{(d-1)/p} \, g_{j,k,\ell}/c_\varphi$ are $(s,p)_{L,-1}$-atoms
in the sense of Subsection \ref{rsa} (in the one-dimensional context).
Applying property (ii) from this subsection we find
\beqq
\| \, \tr f_{2,N} \, | B^s_{p,q}(\R)\| & \lsim &
\Big(\sum_{j=0}^N  \, 2^{j(d-1)q/p}\, \Big(\sum_{k=k_j}^{2^N} |s_{j,k}|^p\Big)^{q/p}\Big)^{1/q}
\\
& \lsim & \Big(\sum_{j=0}^N  \, \Big(\sum_{k=k_j}^{2^N} (1+k)^{d-1} |s_{j,k}|^p\Big)^{q/p}\Big)^{1/q}
\\
& \lsim & \|\, f \, |R\bspq (\Rd)\| \, .
\eeqq
Now we consider convergence of the sequence $\tr f_N$. Let $\sigma_p (1) < s' < s$.
Arguing as before (but taking into account the different normalization of the atoms with respect to $B_{p,p}^{s'} (\R)$) we find
\beqq
\| \, \tr f_{2,N} - \tr f_{2,M}\, |L_1 (\R)\|
& \le  &  \| \, \tr f_{2,N} - \tr f_{2,M}\, |B_{p,p}^{s'} (\R)\|
\\
& \lsim & \Big(\sum_{j=M+1}^N \sum_{k=k_j}^{2^N}
(1+ k)^{d-1} |2^{j(s'-s)}  s_{j,k}|^p\Big)^{1/p}
\\
&& + \Big(\sum_{j=0}^M \sum_{k=\max (2^M, k_j)}^{2^N}
(1+ k)^{d-1} |2^{j(s'-s)}  s_{j,k}|^p\Big)^{1/p}
\eeqq
Since $\| (s_{j,k})_{j,k}\, | b_{p,q,d}\| < \infty$ it follows that the right-hand side tends to $0$ if $M \to \infty$.
The uniform boundedness of
$(\tr f_{2,N})_N$ in $ B^s_{p,q}(\R)$ in combination with the weak convergence of this sequence
yields $\lim_{N \to \infty} \tr f_{2,N} \in \bspq (\R)$
by means of the so-called Fatou property, see \cite{BM,Fr}.
Hence, $\tr f_2 \in \bspq (\R)$. In combination with our knowledge about $f_1$
the claim in case of Besov spaces follows.
\\
{\em Step 2.} Let $f \in R\fspq (\Rd)$.
One can argue as in Step 1.
For the Fatou property of the spaces $\fspq (\R)$ we refer to \cite{Fr}.
\epr


 \subsection{Proofs of the statements  in Subsection  \ref{main16}}


 \subsection*{Proof of Theorem  \ref{main10}}

{\em Step 1.}
The proof of Theorem \ref{main10}(i) follows from formula (\ref{w-100})
and the density of $RC_0^\infty (\Rd)$ in $RW_p^1 (\Rd)$.
\\
{\em Step 2.} Let $f \in RC_0^\infty (\Rd)$.
This is equivalent to $\tr f = f_0 \in RC_0^\infty (\R)$, see Thm. \ref{anfang}.
Observe, that
\[
\frac{\partial^2 f}{\partial x_i\partial x_j}(x)=
\begin{cases}
f_0''(r)\cdot\frac{x_i\cdot x_j}{r^2} - f_0'(r)\, \frac{x_i\cdot x_j}{r^3},
\qquad&\text{if}\quad i\not=j \, ,
\\
& \\
f_0''(r)\cdot\frac{x^2_i}{r^2}-f_0'(r)\cdot \frac{r^2-x^2_i}{r^3},
\qquad&\text{if}\quad i=j.\\
\end{cases}
\]
We fix  $j\in\{1,2,\dots,d\}$ and sum up
\[
\sum_{i=1}^d \left(\frac{\partial^2 f}{\partial x_i\partial x_j}(x)\right)^2=
\frac{f_0''(r)^2}{r^2} \, x_j^2+\frac{f_0'(r)^2}{r^4}\cdot(r^2-x_j^2)\, .
\]
Now we sum up with respect to $j$ and find
\[
\sum_{i,j=1}^d \left(\frac{\partial^2 f}{\partial x_i\partial x_j}(x)\right)^2=
f_0''(r)^2+\frac{d-1}{r^2}\cdot f_0'(r)^2.
\]
Since the terms on the right-hand side are nonnegative this proves the claim
for smooth $f$.
As above the  density argument completes the proof.
\epr

 \subsection*{Proof of Theorem  \ref{main1}}

The formulas (\ref{eq-31})-(\ref{laplace})
have to be combined with the  density of $RC_0^\infty (\Rd)$ in $RW_p^{2m} (\Rd)$.
\epr


 \subsection{Proof of the statements in Subsection
             \ref{main30}}



 \subsubsection{Proof of Lemma \ref{dp}}


{\em Step 1.} Necessity of $p >d$. Let $\varphi \in C_0^\infty (\R)$ be
an even function s.t. $\varphi (0) \neq 0$ and $\supp \varphi \subset [-1/2, 1/2]$.
Since $d \ge 2$ the function
$g_1(t):= \varphi (t)\, |t|^{-1}$, $t\in \R$,  belongs to $RL_p (\R, |t|^{d-1})$ if $p<d$.
Hence $RL_p (\R, |t|^{d-1}) \not \subset S'(\R)$ if $p < d$.
Let $p=d$ and take
$g_2(t):= \varphi (t)\, |t|^{-1} \, (-\log |t|)^{-\alpha}$, $t \in \R\setminus\{0\}$, for $\alpha >0$.
In case $\alpha \, d >1$ we have $g_2 \in RL_d (\R, |t|^{d-1})$.
However, if $\alpha <1$ then $g_2 \not \in S'(\R)$. With $1/d < \alpha < 1$
the claim follows.\\
{\em Step 2.} Sufficiency of $p>d$.
Using H\"older's inequality we find
\[
\int_{-1}^1 |g(t)|\, dt \le \Big(\int_{-1}^1 |g(t)|^p\, |t|^{d-1}\, dt\Big)^{1/p}\,
\Big(\int_{-1}^1 |t|^{-\frac{(d-1)p'}{p}}\, dt\Big)^{1/p'}\, .
\]
The second factor on the right-hand side is finite if, and only if,
\[
(d-1)(p'-1) < 1 \qquad \Longleftrightarrow \qquad d <p \, .
\]
Complemented by the obvious inequality
\[
\int_{|t|>1} |g(t)|^p\, dt \le \int_{|t|>1} |g(t)|^p\, |t|^{d-1} dt
\]
we conclude $L_p (\R, |t|^{d-1}) \hookrightarrow L_1 (\R) + L_p (\R) \subset S'(\R)$.
\epr


 \subsubsection{Proof of Theorem \ref{main8}}


Thm. \ref{main5} implies the equivalence of (ii) and (iii).
Also (i) and (ii) are obviously equivalent.\\ 
{\em Step 1.} We shall prove that (iv) implies (i) and (ii).\\
{\em Substep 1.1}. The $B$-case.
It will be enough to deal with the limiting case.
Let $s = d(\frac 1p - \frac 1d)>0$ ($s>\sigma_p(d)$) and $q=1$.
In addition we assume $1 \le p < d$, where the upper bound results from the previous restriction on $s$, see Fig. 1 in Subsection \ref{main31}.
For $f \in R\bspq (\Rd)$ we select an optimal atomic decomposition of the trace in the sense of Theorem \ref{main5}.
Let $\varphi \in S(\R)$. Then
\beqq
\Big|\int_{-\infty}^\infty && \hspace{-0.7cm} \sum_{j=0}^\infty \sum_{k=0}^\infty
s_{j,k}\, b_{j,k}(t)\, \varphi (t)\, dt\Big|
\\
& \le & 4 \sum_{j=0}^\infty \sum_{k=0}^\infty |s_{j,k}|\, 2^{-j}
\, \| \, b_{j,k}\,|L_\infty (\R)\| \,  \|\, \varphi \, |L_\infty (\R)\|
\\
& \le & 4 \|\, \varphi \, |L_\infty (\R)\| \,
\sum_{j=0}^\infty 2^{j(s-d/p)} \, \sum_{k=0}^\infty |s_{j,k}|\,
\\
& \le & 4 \|\, \varphi \, |L_\infty (\R)\| \,
\Big(\sum_{k=0}^\infty (1+k)^{-(d-1)\frac{p'}{p}} \Big)^{1/p'}
\\
&& \hspace{2cm} \times \quad \sum_{j=0}^\infty 2^{j(s-d/p)} \, \Big(\sum_{k= 0}^\infty (1+k)^{d-1} |s_{j,k}|^p\Big)^{1/p}\, .
\eeqq
Since
\[
\Big(\sum_{k= 0}^\infty (1+k)^{-(d-1)\frac{p'}{p}} \Big)^{1/p'}\,
<\infty
\]
if $1\le p< d$, we obtain
\beqq
\Big|\int_{-\infty}^\infty  \sum_{j=0}^\infty \sum_{k= 0}^\infty
s_{j,k}\, b_{j,k}(t)\, \varphi (t)\, dt\Big|
& \le &  c_1\,  \|\, \varphi \, |L_\infty (\R)\| \, \| \, (s_{j,k})_{j,k}\, |b^s_{p,1,d}\|
\\
& \le &  c_2\,  \|\, \varphi \, |L_\infty (\R)\| \, \| \, f\, |B^s_{p,1} (\Rd)\| \, ,
\eeqq
see Theorem \ref{main5}.
This proves sufficiency for $1 \le p < d$ and
 $s = d(\frac 1p - \frac 1d)$. Now, let $0< p<1$.
Then it is enough to apply the continuous embedding
\[
B^{\frac dp - 1}_{p,1} (\Rd) \hookrightarrow
B^{d - 1}_{1,1} (\Rd) \, ,
\]
see e.g. \cite[2.7.1]{Tr1} or \cite{SiTr}.
\\
{\em Substep 1.2}.
Now we turn to the same implication in case of the $F$-spaces.
Also here an  embedding argument turns out to be sufficient.
For $0 <p\le 1$ and $p < p_1 < \infty$
we have
\[
F^{\frac dp - 1}_{p,\infty} (\Rd) \hookrightarrow B^{\frac{d}{p_1} - 1}_{p_1,1} (\Rd) \, ,
\]
see \cite{Ja0} or \cite{SiTr}. Now the claim follows from Substep 1.1.\\
{\em Step 2.}
Since $\tr$ is an isomorphism of $R\aspq (\Rd)$ onto
$T\aspq (\R, L,d)$ we deduce from Step 1 the implication (iv)
$\Longrightarrow$ (iii).\\
{\em Step 3.} It remains to prove the implication
(i) $\Longrightarrow$ (vi).
We argue by contradiction.
\\
{\em Substep 3.1.} The $B$-case.
Let $s = \frac dp - 1$ and suppose $q>1$.
Oriented at our investigations in Lemma  \ref{test2}
we will use as test functions
\be\label{eq-99}
f_\alpha (x):= \varphi (|x|)\, |x|^{-1} (-\log |x|)^{-\alpha}\, ,
\qquad x\in \Rd\, .
\ee
It is known, see e.g. \cite[Lem.~2.3.1]{RS},
that
\[
f_\alpha \in B^{\frac dp - 1}_{p,q} (\Rd) \qquad \mbox{if, and only if,}
\qquad q\, \alpha >1 \, .
\]
Since $\tr f_\alpha \not\in S'(\R)$ if $\alpha < 1$,
we obtain that $\tr$ does not map into $S'(\R)$ as long as
$1/q < \alpha < 1$.
\\
{\em Substep 3.2.} The $F$-case.
This time it holds
\[
f_\alpha \in F^{\frac dp - 1}_{p,\infty} (\Rd) \qquad \mbox{if, and only if,}
\qquad p\, \alpha >1\, ,
\]
 see  \cite[Lem.~2.3.1]{RS}. Choosing $1/p < \alpha < 1$ we obtain
that $\tr$ does not map $ F^{\frac dp - 1}_{p,\infty} (\Rd) $ into $S'(\R)$.
\epr


 \subsection{Proof of the assertions in Subsection \ref{main31}}


\subsection*{Proof of Theorem \ref{main32}}

Comparing our atomic decomposition with that one for weighted spaces obtained in
\cite{HP} it is essentially a question of renormalization of the atoms.
This is enough to prove
${T}\aspq (\R,L,d) \hookrightarrow R\aspq (\R, w_{d-1})$.
To see the converse one has to start with the fact that
$f \in R\aspq (\R, w_{d-1})$ is even.
This allows to decompose  $f$ into sum of atoms that are even as
well, see (\ref{ws-3}) and (\ref{ws-2}) for this argument.
\epr

\subsection*{Proof of Remark \ref{delta}}

The regularity of the $\delta$ distribution is calculated at several places,
see e.g. \cite[Remark~2.2.4/3]{RS}.
The argument, used in this reference, comes from Fourier analysis
and transfers to the weighted case.
For the Fourier analytic characterization of $\aspq (\R,w_{d-1})$
we refer to \cite{bui-1,bui-2} and \cite{HP}.
\epr


 \subsection{Proof of the assertions in Subsection \ref{main14}}



 \subsubsection*{Proof of  Corollary \ref{nice}}


We shall only prove part (i). The proof for the Lizorkin-Triebel spaces is similar.
By our trace theorem we have
\[
\| \, f_0 \, | {T}\bspq (\R,L,d)\|  \lsim   \| \, f \, | \bspq (\Rd) \|
\]
if  $L>[s]+1$, cf. Theorem \ref{main5}. Thus, it is sufficient to prove that
\begin{equation}\label{LSnew}
\| \, f_0 \, | \bspq(\R)\|  \lsim \,  \tau^{-(d-1)/p}\,\| \, f_0 \, | {T}\bspq (\R,L,d)\| \,.
\end{equation}
The trace $f_0\in {T}\bspq (\R,L,d)$  can be represented in the form
\begin{equation}
\label{newaa}
f_0(t)= \sum_{j=0}^\infty \sum_{k=0}^\infty s_{j,k}\, g_{j,k} (t)
\end{equation}
(convergence in $L_{\max(1,p)} (\R, |t|^{d-1})$), where  the sequence $(s_{j,k})_{j,k}$ belongs to $b^s_{p,q,d}$, cf. \reff{new}.
Let $\varphi\in C^\infty (\R)$ be an even function such that  $\varphi(t)=0$ if $|t|\le \frac{1}{2}$ and $\varphi(t)=1$ if $|t|\ge 1$. For any $\tau > 0$
we define  $\varphi_\tau (t) = \varphi (\tau^{-1}t)$.
We will consider two cases: $\tau \ge 2$ and $0<\tau < 2$.
\\
{\em Case 1}. Let $\tau \ge 2$. Under this assumption
any function $\varphi_\tau \, g_{j,k}$ is an even $L$-atom
centered at the same interval as $g_{j,k}$ itself (up to a general constant depending on $\varphi$), see Definition
\ref{latom}. For any $j \in \N_0$ we define a nonnegative integer $k_j$ by
\[
k_j:= \max \{k \in \N_0: \quad 2^{-j}\, (k+1) \le  \tau /2\}
\, .
\]
Hence,  $\varphi_\tau \,  g_{j,k} =0$ if $k < k_j$. Furthermore,
the functions $2^{-j(s-1/p)} \, \varphi_\tau \,  g_{j,k} $, $k \ge 1$,
restricted either to the positive or negative half axis,
are $(s,p)_{L,-1}$-atoms in the sense of Definition \ref{atomunw}
up to a universal constant $c$. The  functions
$2^{-j(s-1/p)} \, \varphi_\tau \,  g_{j,0} $
are $(s,p)_{L,-1}$-atoms as well (again up to a universal constant). We obtain
\[
f_0(t) =\varphi_\tau (t) \, f_0(t) = \sum_{j=0}^{\infty} \sum_{k=k_j}^{\infty} s_{j,k}\, \varphi_\tau (t) \, g_{j,k} (t)
\]
and  applying (\ref{atom2})
(which is also valid for $d=1$) we arrive at the estimate
\begin{eqnarray}\label{LSnew-1}
\|\, f_0 \, |\bspq(\R)\| & \lsim &
\bigg(\sum_{j=0}^\infty 2^{j(s - \frac 1p)q}\, \bigg(
\sum_{k=k_j}^\infty \,|s_{j,k}|^p\bigg)^{q/p} \bigg)^{1/q}
\nonumber
\\
& \lsim &   \tau^{\frac{1-d}{p}}\,
\bigg(\sum_{j=0}^\infty 2^{j(s - \frac dp)q}\, \bigg(
\sum_{k=0}^\infty (1+k)^{d-1}\,|s_{j,k}|^p\bigg)^{q/p} \bigg)^{1/q}
\nonumber
\\
& = & \, \tau^{\frac{1-d}{p}} \, \| \,s \, | b^s_{p,q,d}\|
\end{eqnarray}
since $k_j\sim 2^j\, \tau $. Taking the infimum with respect to all atomic representations of $f_0$ we  have proved \reff{LSnew}.
\\
{\em Case 2}. Let $0 <\tau < 2$. \\
{\em Step 1.} We assume $s < d/p$.
Then we define $j_0\in \N_0$ via the relation
$2^{-j_0} \le \tau <  2^{-j_0+1}$.
Further, we put $K_j:= \max (1, 2^{j-j_0-1}-1)$.
Now we decompose $f_0 $ into four sums
\begin{eqnarray*}
f_0 (t)  =  \varphi_\tau (t) \, f_0(t) & = &
\sum_{j=0}^{j_0+1}  s_{j,0}\, \varphi_\tau (t) \, g_{j,0} (t)
\quad  + \quad  \sum_{j=j_0}^{\infty}
\sum_{k=K_j}^{2^{j-j_0+1}} s_{j,k}\, \varphi_\tau(t)
\, g_{j,k} (t) \\
 & + &
 \sum_{j=j_0}^{\infty} \sum_{k=2^{j-j_0+1}+1}^{\infty} s_{j,k}\,  g_{j,k} (t)
\quad + \quad
\sum_{j=0}^{j_0-1} \sum_{k=1}^{\infty} s_{j,k}\,  g_{j,k} (t)
 \nonumber
\\
&& \\
& = & f_1(t) + \ldots \, + f_4(t) \, ,
\nonumber
\end{eqnarray*}
with $f_4 =0$ if $j_0 =0$.
Observe
\[
\supp   f_i   \subset  \{t: \: |t| \ge \,  \tau \}\, , \qquad
i=3,4 \, ,
\]
whereas the supports of the functions $\varphi_\tau  \, g_{j,k}$,
occuring in the defintions of $f_1$ and $f_2$, may have  nontrivial intersections with the interval $ (\tau /2, \tau)$.
The function $f_1$ belongs to $C^L$ and has compact support.
The functions $f_2, \, f_3$, and $f_4$ are supported on
$\{ t:\, |t|\ge \tau / 2 \}$. Thus, the known convergence in $L_{\max(1,p)} (\R, |t|^{d-1})$ implies the convergence in $\cs'(\R)$.
As in {\em Case 1} the functions
$2^{-j(s-1/p)} \,   g_{j,k} $, $k \ge 1$,
restricted either to the positive or negative half axis,
are $(s,p)_{L,-1}$-atoms in the sense of Definition \ref{atomunw}.
An easy calculation shows
that also  the functions $2^{-j(s-1/p)} \, \varphi_\tau\,   g_{j,k} $,
$j \ge j_0$, are $(s,p)_{L,-1}$-atoms (up to a universal constant). Hence we may employ (\ref{atom2})
and obtain
\beqq
\|\, f_2 + f_3 \,  |\bspq (\R)\|
&\lsim &
\bigg(\sum_{j=j_0}^{\infty} 2^{j(s - \frac 1p)q}\, \bigg(
\sum_{k=K_j}^\infty \,|s_{j,k}|^p\bigg)^{q/p}\bigg)^{1/q}
\\
& \lsim & \tau^{\frac{1-d}{p}}
\bigg(\sum_{j=j_0}^{\infty} 2^{j(s - \frac dp)q}\, \bigg(
\sum_{k=K_j}^\infty (1+k)^{d-1}\,
|s_{j,k}|^p\bigg)^{q/p} \bigg)^{1/q}
\eeqq
as well as
\beqq
\|\, f_4  \,  |\bspq (\R)\|
&\lsim &
\bigg(\sum_{j=0}^{j_0-1} 2^{j(s - \frac 1p)q}\, \bigg(
\sum_{k=1}^\infty \,|s_{j,k}|^p\bigg)^{q/p}\bigg)^{1/q}
\\
& \lsim &  2^{j_0 \frac{d-1}{p}}\,
\bigg(\sum_{j=0}^{\infty} 2^{j(s - \frac dp)q}\, \bigg(
\sum_{k=0}^\infty (1+|k|)^{d-1}\,
|s_{j,k}|^p\bigg)^{q/p} \bigg)^{1/q}
\\
& \lsim & \tau^{\frac{1-d}{p}}\,
\bigg(\sum_{j=0}^{\infty} 2^{j(s - \frac dp)q}\, \bigg(
\sum_{k=0}^\infty (1+|k|)^{d-1}\,
|s_{j,k}|^p\bigg)^{q/p} \bigg)^{1/q}\, .
\eeqq
Now we turn to the estimate of $f_1$.
First we deal with the estimate of the quasi-norm of the functions $\varphi_\tau  \, g_{j,0}$.
Let in addition $s \ge 1/p$.
Employing  the Moser-type estimate of Lemma 5.3.7/1  in \cite{RS}
(applied with $r=\infty$)
we obtain
\beq\label{ws-1111}
\|\,  \varphi_\tau  \, g_{j,0}\, | \bspq(\R)\|\,
& \lsim & \|\,  \varphi_\tau  \,  | \bspq(\R)\|\,\|\,   g_{j,0}\, | L_\infty(\R)\|
+ \|\,  \varphi_\tau  \, | L_\infty (\R)\|\,\|\, g_{j,0}\, | \bspq(\R)\|
\nonumber
\\
& \lsim &  \tau^{-(s-1/p)}\, \|\,  \varphi \, | B^s_{p,q} (\R)\|
+ \|\,  \varphi \, | L_\infty (\R)\|\,\|\,   g_{j,0}\, | B^s_{p,q} (\R)\|
\nonumber
\\
& \lsim &  \tau^{-(s-1/p)} + 2^{j(s-1/p)}
\nonumber
\\
& \lsim & 2^{j_0(s-1/p)} \, ,
\eeq
since the functions $2^{-j(s-1/p)} \,   g_{j,0} $ are atoms and
$j \le j_0+1$.
If $s < 1/p$, we argue by using real interpolation.
Because of
\[
B^s_{p,q} (\R)  =   \Big(B^{s_0}_{p,q_0} (\R),
L_{p} (\R)\Big)_{\Theta,q} \, , \qquad  s= (1-\Theta)s_0 > \sigma_p (1)\, ,
\]
see \cite{DVX}, an application of the interpolation inequality
\[
\|\,  \varphi_\tau  \, g_{j,0}\, | \bspq(\R)\|  \lsim  \|\,  \varphi_\tau  \, g_{j,0}\, | B^{s_0}_{p,q}(\R)\|^{1-\Theta}\,
\|\,  \varphi_\tau  \, g_{j,0}\, | L_p(\R)\|^{\Theta}\,
\]
yields (\ref{ws-1111}) for all $s> \sigma_p (1)$.
With  $r:= \min (1,p,q)$  and $\sigma_p (d) < s < d/p$ we conclude
\beqq
\|\, f_1\, | \bspq(\R)\|^r  & \le &  \sum_{j=0}^{j_0+1}  |s_{j,0}|^r \, \| \, \varphi_\tau  \, b_{j,0} \, |\bspq (\R)\|^r
\\
& \lsim &  \, 2^{j_0(s-\frac{1}{p})r} \,
\sum_{j=0}^{j_0+1}  |s_{j,0}|^r \, .
\\
& \lsim &  \tau^{\frac{r(1-d)}{p}}
\sum_{j=0}^{j_0+1}  2^{r (j_0-j)(s-\frac d p)} 2^{j(s-\frac d p)r} \,  |s_{j,0}|^r
\\
& \lsim &  \tau^{\frac{r(1-d)}{p}} \,
\, \Big(\sup_{j=0, \ldots \, ,j_0+1}  2^{j(s-\frac d p)} \,  |s_{j,0}|\Big)^r
\\
& \lsim &  \tau^{\frac{r(1-d)}{p}} \,
\, \| \, f_1 \, |TB^s_{p,\infty} (\R)\|^r \, .
\eeqq
This proves the claim for $s<d/p$.
\\
{\em Step 2.} Let $s \ge d/p$.
As in Step 1 we define $j_0\in \N$ via the relation
$2^{-j_0} \le \tau <  2^{-j_0+1}$.
We shall use  ${T}\aspq (\R,L,d)= R\aspq (\R,w_{d-1})$,  cf. Theorem \ref{main32}.
Alternatively one could use interpolation, see Propositions
\ref{interpol2}, \ref{interpol1}.
The spaces $R\aspq (\R,w_{d-1})$ allow a characterization
by Daubechies wavelets, see \cite{HSk} for Besov spaces
and \cite{IS} for Lizorkin-Triebel spaces.
The same is true with respect to the ordinary spaces
$\aspq (\R)$, see e.g. \cite[Thm.~1.20]{Tr08}.
Let $\phi$ denote an appropriate scaling function
and $\Psi$ an associated Daubechies wavelet of sufficiently high order.
Let
\[
\phi_{0,\ell} (t) := \phi (t-\ell) \qquad \mbox{and}
\qquad \Psi_{j,\ell}(t):= 2^{j/2}\, \Psi_{j,\ell} (2^j t -\ell)\, , \qquad \ell
\in \Z\, , \quad j \in \N_0\, .
\]
Since  $\Psi$ has compact support, say $\supp \Psi \subset [-2^N,2^N]$
for some $N \in \N$, and $\supp f_0 \subset \{t\in \R: |t|\ge \tau\}$ we find
that
\[
 \langle  f_0,\Psi_{j,\ell} \rangle =0 \qquad \mbox{if}
\qquad  j-j_0 \ge N
\qquad \mbox{and}
\qquad
|\ell|\le \,  2^{j-j_0}- 2^N \, .
\]
Hence, $f_0$ has a wavelet expansion given by
\beqq
f_0  & = &  \sum_{\ell \in \Z} \langle f_0,\phi_{0,\ell} \rangle \, \phi_{0,\ell}
+    \sum_{j=0}^ {j_0+N-1}    \sum_{\ell \in \Z}
\langle f_0, \Psi_{j,\ell} \rangle \,   \Psi_{j,\ell}
  + \sum_{j=j_0 + N}^ {\infty}    \sum_{|\ell| > 2^ {j-j_0}-2^N}
\langle f_0, \Psi_{j,\ell} \rangle \,   \Psi_{j,\ell}\,.
\\
& = & f_1 + f_2 + f_3 \, .
\eeqq
By the references given above it follows
\beqq
\| \, f_1 \, | \bspq(\R,w_{d-1})\| & \asymp &
\bigg(
\sum_{\ell \in \Z} \,|
\langle f_0,\phi_{0,\ell}\rangle |^p (1+|\ell|)^{d-1} \bigg)^{1/p}
\\
\| \, f_2 \, | \bspq(\R,w_{d-1})\|
& \asymp &
 \bigg(\sum_{j=0}^{j_0 + N-1} 2^{j(s + \frac 12 - \frac dp)q}\, \bigg(
\sum_{\ell \in \Z} \,|\langle  f_0, \Psi_{j,\ell} \rangle |^p
(1+|\ell|)^{d-1} \bigg)^{q/p} \bigg)^{1/q}\,
\\
\| \, f_3 \, | \bspq(\R,w_{d-1})\|
& \asymp &
 \bigg(\sum_{j=j_0+N}^\infty 2^{j(s + \frac 12 - \frac dp)q}\, \bigg(
\sum_{|\ell| \ge 2^{j-j_0} -2^N} \,|\langle  f_0, \Psi_{j,\ell} \rangle |^p
(1+|\ell|)^{d-1} \bigg)^{q/p} \bigg)^{1/q}\, .
\eeqq
The quasi-norm in the unweighted spaces is obtained by
deleting the factor $ 2^{-j(d-1)/p}\, (1+|\ell|)^{d-1}$, see \cite[Thm.~1.20]{Tr08}.
This  immediately implies
\beqq
\| \, f_1 \, | \bspq(\R)\| & \lsim &
\| \, f_1 \, | \bspq(\R,w_{d-1})\|\, ,
\\
\| \, f_2 \, | \bspq(\R)\| & \lsim &
2^{(j_0 + N)(d-1)/p}\,  \| \, f_2 \, | \bspq(\R,w_{d-1})\|\, .
\eeqq
Moreover, we also obtain
\[
\| \, f_3 \, | \bspq(\R)\| \lsim
2^{(j_0 + N)(d-1)/p}\,  \| \, f_3 \, | \bspq(\R,w_{d-1})\|\, .
\]
This proves \reff{LSnew} in case  $s \ge d/p$ and $0 < \tau < 2$.
\epr


 \subsubsection*{Proof of  Corollary \ref{ok}}


We concentrate on the proof in case of Besov spaces.
The proof for Lizorkin-Triebel spaces is similar.
\\
{\em Step 1.}
We claim that $g \in T\bspq (\R,L,d)$.
We argue as in {\em Case 2, Step 2} of the proof of Corollary \ref{nice}.
\\
Under the given restrictions
$g \in R\bspq(\R)$ has a wavelet expansion of the form
\[
g   =   \sum_{|\ell| \le c_1} \langle g,\phi_{0,\ell} \rangle \, \phi_{0,\ell}
  + \sum_{j=0}^ {\infty}    \sum_{|\ell| \le c_1  2^ {j}}
\langle g_0, \Psi_{j,\ell} \rangle \,   \Psi_{j,\ell}
\]
with an appropriate constant $c_1$.
Since $g$ is even we obtain
\[
g   =   \sum_{|\ell| \le c_1} \langle g,\phi_{0,\ell} \rangle \,
\frac{\phi_{0,\ell} (t) + \phi_{0,\ell} (-t)}{2}
  + \sum_{j=0}^ {\infty}    \sum_{|\ell| \le c_1 2^ {j}}
\langle g_0, \Psi_{j,\ell} \rangle \,
\frac{\Psi_{j,\ell}(t) + \Psi_{j,\ell} (-t)}{2}\, .
\]
The functions $2^{-j/2}(\Psi_{j,\ell}(t) + \Psi_{j,\ell} (-t))$ are even
$L$-atoms (up to a universal constant) centered at $c_2\, I_{j,\ell}$, where
\[
I_{j,k}:= [-2^{-j}(\ell+1),-2^{-j}\ell] \cup [2^{-j}\ell,2^{-j}(\ell +1)]
\]
(modification if $\ell =0$, see Definition  \ref{tracedef}).
The constant $c_2>1$ depends on the size of the supports of the generators
$\phi$  and $\Psi$. Without proof we mention that Theorem \ref{main5}
remains true also for those more general decompositions.
This implies
\beqq
\| \, g \, | && \hspace*{-1.5cm}  T\bspq (\R,L,d)\|  \asymp
\Big(\sum_{|\ell| \le c_1 b} (1+|\ell|)^{d-1}
|\langle g,\phi_{0,\ell} \rangle|^p \Big)^{1/p}
\\
& & \qquad + \quad \Big(\sum_{j=0}^ {\infty}   2^{j(s-d/p)q}\Big( \sum_{|\ell| \le c_1 2^ {j}}
(1+|\ell|)^{d-1} |2^{j/2} \, \langle g, \Psi_{j,\ell} \rangle|^p \Big)^{q/p}\Big)^{1/q}
\\
& \lsim &
\Big(\sum_{|\ell| \le c_1 b} |\langle g,\phi_{0,\ell} \rangle|^p \Big)^{1/p}
   +  \Big(\sum_{j=0}^ {\infty}   2^{j(s+ \frac 12 - \frac1p)q}
\Big( \sum_{|\ell| \le c_1 2^ {j}}
|\langle g, \Psi_{j,\ell} \rangle|^p \Big)^{q/p}\Big)^{1/q}
\\
& \lsim & \| \, g \, |\bspq (\R)\|\, ,
\eeqq
see e.g. \cite[Thm.~1.20]{Tr08} for the last step. This proves the claim.
\\
{\em Step 2.} Since $g$ belongs to $T\bspq (\R,L,d)$
we derive by means of Theorem \ref{main5} that $f : = \ext g$ is an element
of $R\bspq (\Rd)$ and
\[
\| \, f\, |R\bspq(\Rd)\| \lsim \| \, g \, | T\bspq (\R,L,d)\|
\lsim \| \, g \, | \bspq (\R)\|\, .
\]
Since $\supp f \subset \{x: \: |x|\ge a\}$ Corollary \ref{nice}
yields
\[
\| \, g \, | \bspq(\R)\|  \lsim \,  a^{-(d-1)/p}\,
\| \, f \, | \bspq (\Rd) \|\, ,
\]
because of $f_0 =g$.
This completes the proof.
\epr

\begin{Rem}\label{equi}
\rm
A closer look onto the proof shows that
\[
a^{(d-1)/p}\,
\| \, g \, | \aspq(\R)\|  \lsim \| \, f\, |R\aspq(\Rd)\|\lsim
b^{(d-1)/p}\,
\| \, g \, | \aspq(\R)\|
\]
with constants independent of $g$, $a>0$ and $b\ge 1$.
\end{Rem}


 \subsubsection*{Proof of  Corollaries \ref{smooth}, \ref{uc}}


{\em Step 1.} Proof of Cor. \ref{smooth}.
The function $\varphi$ is a pointwise multiplier for  the spaces $\aspq (\Rd)$, see e.g. \cite[4.8]{RS}.
Hence, with $f$ also the product $\varphi \, f $ belongs to
$R\aspq (\Rd)$ and we can apply Thm.  \ref{spur0} with respect to this product.
Concerning the sharp  embedding relations for the spaces
 $\aspq (\R)$ into H\"older-Zygmund spaces
we refer to \cite{SiTr}  and the references given there.
This proves the assertion for $\varphi_0 \, f_0$. A further application of Theorem
\ref{infty} finishes the proof.
\\
Observe, that we do not need the assumption $s>\sigma_{p,q}(d)$ in case of Lizorkin-Triebel spaces. We may argue with $RF^s_{p,\infty} (\Rd)$ first and use the elementary
embedding  $RF^s_{p,q} (\Rd) \hookrightarrow RF^s_{p,\infty} (\Rd)$ afterwards.
\\
{\em Step 2.} Proof of Cor. \ref{uc}. The arguments are as above.
Concerning te embedding relations of the spaces $\aspq(\R)$ into the space
of uniformly continuous and bounded functions we also refer to
\cite{SiTr} and the references given there.
\epr

\begin{Rem}
\rm A different proof of Cor. \ref{uc}, restricted to Besov spaces,  has been given in \cite{SS}.
\end{Rem}


\subsection{Test functions}
\label{testfunctions}


Using our previous results, in particular Corollary \ref{ok},
we shall investigate the regularity of certain families of radial
test functions.

\begin{Lem}\label{test2}
Let $0 <\alpha < \min (1, 1/p)$.
Let $\varphi \in C_0^\infty (\R)$ be an
even function such that
$\supp \varphi \subset  [-2,-1/2]\cup[1/2,2]$ and
$\varphi (1)\neq 0$.\\
{\rm (i)} The function
\be\label{eq-97}
f_\alpha (x):= \varphi (|x|)\,\,  |\, |x|-1\, |^{-\alpha}\, , \qquad x \in \Rd\, ,
\ee
belongs to $B^{\frac 1p - \alpha}_{p,\infty} (\Rd)$ if
\be\label{eq-97q}
\alpha < \frac 1p - \sigma_p (d).
\ee
{\rm (ii)} Suppose $\frac 1p -\alpha > \sigma_p (1)$. Then $f_\alpha$ does not belong to
$B^{\frac 1p - \alpha}_{p,q} (\Rd)$ for any $q<\infty$.
\\
{\rm (iii)} Under the same restriction as in (ii) we have that
$f_\alpha$ does not belong to
$F^{\frac 1p - \alpha}_{p,\infty} (\Rd)$.
\end{Lem}

\noindent
\bpr
{\em Step 1.} Proof of (i).
Let $\wt{\varphi} \in C_0^\infty (\R)$ be a function such that $\supp
\wt{\varphi} \subset  [1/2,2]$.
Then the regularity of
\[
g_\alpha (t):= \wt{\varphi} (t)\, |t-1|^{-\alpha}\, , \qquad t \in \R\, ,
\]
is well understood, cf. e.g. \cite[Lem.~2.3.1/1]{RS}.
One has
$g_\alpha \in B^{\frac 1p - \alpha}_{p,\infty} (\R)$
as long as $0 <\alpha <  \min (1, 1/p)$.
An application of Corollary \ref{ok} yields the claim.\\
{\em Step 2}. Proof of (ii) and (iii).
It is also known, see again  \cite[Lem.~2.3.1/1]{RS}, that
\[
g_\alpha \not\in (B^{\frac 1p - \alpha}_{p,q} (\R)
\cup F^{\frac 1p - \alpha}_{p,\infty} (\R))\, , \qquad
0 < q<\infty\, , \quad 0 <\alpha < \min (1, 1/p)\, .
\]
These properties do not change when we ''add'' the reflection of $g_\alpha$
to the left half of the real axis. With other words, if we replace
$\wt{\varphi}$ by $\varphi$ itself we do not change the regularity properties.
Now we use Thm. \ref{spur0}.
\epr

\begin{Rem}{\rm
Let $\delta>0$. Then also the regularity of functions like
\be\label{eq-97b}
f_{\alpha, \delta} (x):= \varphi (|x|)\,\,
|\, |x|-1|^{-\alpha}\, (-\log |\, |x|-1| )^{-\delta}
, \qquad x \in \Rd\, ,
\ee
can be checked in this way. With the help of the parameter $\delta$ one can see the microscopic index $q$. We refer to \cite[5.6.9]{Stein1} or
\cite[Lem.~2.3.1/1]{RS} for details.}
\end{Rem}

\begin{Lem}
Let $\alpha >0$.\\
{\rm (i)}
Then the function
\be\label{eq-98}
\Phi_\alpha (x):= \max (0,(1-|x|^2))^{\alpha}\, , \qquad x \in \Rd\, ,
\ee
belongs to $B^{\frac 1p + \alpha}_{p,\infty} (\Rd)$ if
\be\label{eq-98q}
\frac 1p + \alpha > \sigma_p (d)\, .
\ee
{\rm (ii)} Suppose $\frac 1p + \alpha > \sigma_p (1)$.
Then $\Phi_\alpha$ does not belong to
$B^{\frac 1p + \alpha}_{p,q} (\Rd)$ for any $q<\infty$.
\\
{\rm (iii)} Under the same restrictions as in (ii) we have that
$\Phi_\alpha$ does not belong to
$F^{\frac 1p + \alpha}_{p,\infty} (\Rd)$.
\end{Lem}

\noindent
\bpr
{\em Step 1.} Proof of (i).
First we investigate the one-dimensional case.
Let $\psi_1, \psi_2 \in C_0^\infty (\R)$ be such
that
\[
\supp \psi_1 \subset [-1/2,\infty), \quad
\supp \psi_2 \subset (-\infty, 1/2] \quad \mbox{and}\quad \psi_1 (t) + \psi_2(t) =1
\]
for all $t\in \R$. We put
$\Phi_{i,\alpha}:= \psi_i \, \Phi_\alpha $, $i=1,2$.
Then $\Phi_{1,\alpha} $ behaves near $1$ like
\[
\phi_\alpha (t):= \left\{
\begin{array}{lll}
t^\alpha & & \mbox{if} \quad t>0\, ,
\\
0 && \mbox{if} \quad t<0\, ,
\end{array}
\right.
\]
near the origin. The regularity of $\phi_\alpha$ is well understood, we refer
to \cite[Lem.~2.3.1]{RS}. As above the transfer to general dimensions
$d>1$ is done by Corollary \ref{ok}.
\\
{\em Step 2.}
To prove the statements in (ii) and (iii) we argue by contradiction.
If $\Phi_\alpha$ belongs to $R\aspq (\Rd)$, then also
$ \varphi \, \Phi_\alpha$ belongs to $R\aspq (\Rd)$
for any smooth radial $\varphi$.
Choosing $\varphi$ s.t. $0 \not\in \supp \varphi$, we my apply Thm. \ref{spur0}
to conclude that $\tr (\varphi \, \Phi_\alpha) \in \aspq(\R)$.
But this implies $\Phi_{2,\alpha} \in \aspq (\R)$.
In the one-dimensional case necessary and sufficient conditions are known,
we refer again to \cite[Lem.~2.3.1]{RS}.
\epr

Next we shall consider smooth functions supported in thin annuli.

\begin{Lem}\label{test}
Let $d\ge 2$, $0< p,q \le \infty$ and $s> \sigma_p(d)$.
Let $\varphi \in C_0^\infty (\R)$ be  an even  function such that
$\varphi(1)=1$ and $\supp \varphi  \subset [-2,-1/2]\cup [1/2, 2]$.
Then the  functions
\[
f_{j,\lambda}(y):= \varphi (2^j|y|-\lambda)\, , \qquad y \in \Rd\, ,
\quad j \in \N\, , \quad \lambda >0 \, .
\]
have the following properties:
\beq\label{eq-100}
\supp f_{j,\lambda} & \subset & \{y: \quad (\lambda -2)\, 2^{-j} \le |y|
\le (2+\lambda)\, 2^{-j}\}\, , \\
\label{eq-93}
\| \, f_{j,\lambda} \, |R\bspq(\Rd)\| \,  & \asymp  &
 \, 2^{j(s-\frac dp)} \, \lambda^{(d-1)/p}\,
\eeq
with constants in $\asymp$ independent of $\lambda >2$ and $j\in \N$.
\end{Lem}

\noindent
\bpr
{\em Step 1.} Estimate from above in (\ref{eq-93}).
It will be convenient to use
the atomic characterizations described in Subsection \ref{rsa}.
Therefore we shall use the decompositions of unity from
(\ref{res})-(\ref{res3}).
Thanks to the support restrictions for the functions $\psi_{j,k,\ell}$ we obtain
\[
f_{j,\lambda}(y)= \sum_{\max (0,\lambda - 2-n_0) \le k \le \lambda + 2+n_0} \sum_{\ell=1}^{C(d,k)}
( \varphi (2^j|y|-\lambda) \, \psi_{j,k,\ell} (y) )
\]
where $n_0$ is a fixed number ($n_0 \ge 18$ would be sufficient).
The functions
\[
a_{j,k,\ell} (y):= 2^{-j(s-\frac dp)}\, \varphi (2^j|y|-\lambda) \, \psi_{j,k,\ell} (y)
\]
are $(s,p)_{M,-1}$-atoms for any $M$ (up to a universal constant).
Hence
\beqq
\| \, f_{j,\lambda} \, |R\bspq(\Rd)\| & \lsim &
 \Big(\sum_{\max (0,\lambda - 2-n_0) \le k \le \lambda + 2+n_0} k^{d-1}
2^{j(s-\frac dp)p}\Big)^{1/p}
\\
& \lsim &
2^{j(s-\frac dp)} \lambda^{(d-1)/p}\, .
\eeqq
{\em Step 2.} Estimate from below.
\\
{\em Substep 2.1}.
First we deal with $p=\infty$.
By construction
$f_{j,\lambda}(y) =1 $ if $|y|=(1+\lambda)\, 2^{-j}$.
Furthermore, calculating the derivatives of
$f_{j,\lambda}(y_1, 0, \ldots , 0)$, $y_1\in \R$, it is immediate that
\be\label{eq-101}
\| \, f_{j,\lambda} \, |C^m (\Rd)\| \asymp 2^{jm}
\ee
for all $m \in \N_0$.
Now we argue by contradiction. We fix $s>1$, $q_1 \in (0,\infty]$
and assume that
\[
\| \, f_{j,\lambda} \, |B^s_{\infty,q_1} (\Rd)\| \le
 \, \phi (j,\lambda)\, 2^{js}\, ,
\]
where $\phi :\, \N \times [1,\infty) \to (0,1) $
and $\lim_{\ell \to \infty} \, \phi (j_\ell, \lambda_{\ell}) =0$
for some sequence $(j_\ell, \lambda_{\ell})_{\ell} \subset
\N \times [1,\infty)$.
We choose $\Theta \in (0,1)$ s.t. $m = \Theta \, s$ and
$q=1$.
Real interpolation between $C(\Rd)$ and $B^s_{\infty,q_1} (\Rd)$ yields
\[
\| \, f_{j,\lambda} \, |B^m_{\infty,1} (\Rd)\| \le c\, 2^{jm} \,
(\phi (j,\lambda))^\Theta \, ,
\]
where $c$ is independent of $j$ and $\lambda$, see the proof of Thm. \ref{infty}.
The continuous embedding $B^m_{\infty,1} (\Rd) \hookrightarrow C^m (\Rd)$
leads to a contradiction with (\ref{eq-101}).
\\
Now let $0< s < 1$.
We interpolate between $B^s_{\infty,q_1} (\Rd)$ and $B^2_{\infty,q} (\Rd)$.
By arguing as above we could improve the estimate from above with respect to the space $B^1_{\infty,1} (\Rd)$. Since
$B^1_{\infty,1} (\Rd) \hookrightarrow C^1 (\Rd)$
this contradicts again (\ref{eq-101}).
Hence the claim is proved with $p=\infty$, $0< q \le \infty$, and $s>0$.
\\
{\em Substep 2.2.}
Also obvious is the behaviour in $L_p (\Rd)$.
For $0 < p \le \infty$ we have
\be\label{eq-103}
\| \, f_{j,\lambda} \, |L_p (\Rd)\| \asymp 2^{-jd/p}\,  \lambda^{(d-1)/p}\, .
\ee
A few more calculations yield
\be\label{eq-104}
\| \, f_{j,\lambda} \, |W^1_p (\Rd)\| \asymp 2^{j(1-d/p)}\,  \lambda^{(d-1)/p}\, ,
\ee
as long as $1\le p \le \infty$.\\
{\em Substep 2.3}.
Let $p_1<\infty$.  We assume that for some fixed $s_1>\sigma_{p_1} (d)$ and
$q_1 \in (0,\infty]$
\[
\| \, f_{j,\lambda} \, |B^{s_1}_{p_1,q_1} (\Rd)\| \le
 \, \phi (j,\lambda)\, 2^{j(s_1-d/p_1)}\, \lambda^{(d-1)/p_1} 
\]
holds, where $\phi $ is as above.
Complex interpolation between $B^{s_1}_{p_1,q_1} (\Rd)$ and
$B^{s_2}_{\infty,q_0} (\Rd)$, $s_2 >0$,
yields an improvement of our estimate with respect to
$B^s_{p,q} (\Rd)$, where  $p>p_1$ is at our disposal.
For $s_2$ large  we can choose $p>1$
s.t. $s= (1-\Theta)s_2 + \Theta \, s_1 >1$.
Now we need a further interpolation,
this time real, between
$B^s_{p,q} (\Rd)$ and $L_p (\Rd)$, improving the estimate for
$B^1_{p,1} (\Rd)$ in this way.
But $B^1_{p,1} (\Rd) \hookrightarrow W^1_p (\Rd)$ and so we found a contradiction to (\ref{eq-104}).
\epr

\begin{Rem}\label{testf}
\rm
Obviously there is no $q$-dependence in Lemma \ref{test}.
As an immediate consequence of the elementary embeddings
\[
B^s_{p, \min (p,q)} (\Rd) \hookrightarrow \fspq (\Rd) \hookrightarrow
B^s_{p, \max (p,q)} (\Rd)\, ,
\]
see \cite[]{Tr1}, and (\ref{eq-93}) we obtain
\[
\| \, f_{j,\lambda} \, |\fspq(\Rd)\| \,   \asymp
 \, 2^{j(s-\frac dp)} \, \lambda^{(d-1)/p}\, , \qquad \lambda >2, \quad j \in \N\, .
\]
\end{Rem}
{~}\\

Some extremal functions in ${A}^{d/p}_{p,q} (\Rd)$
have been investigated by Bourdaud \cite{Bou2}, for
${B}^s_{p,p} (\Rd)$ see also Triebel \cite{Trex}. We recall the result obtained in
 \cite{Bou2}.
For $(\alpha,\sigma)\in \R^2$ we define
\be\label{spezial}
f_{\alpha,\sigma} (x) :=  \psi (x) \, \Big| \, \log |x|\, \Big|^{\alpha}\,
\Big| \, \log |\, \log |x|\, |  \Big|^{-\sigma}\,
, \qquad x \in \Rd\, .
\ee
Furthermore we define a set $U_t \subset \R^2$
as follows:
\[
U_t :=
 \left\{
\begin{array}{lll}
(\alpha = 0  \quad \mbox{and}  \quad \sigma >0)  \quad \mbox{or}  \quad
\alpha <0  & \mbox{if} & \quad t=1\, ,
\\
&& \\
(\alpha = 1-1/t  \quad \mbox{and}  \quad \sigma >1/t)  \quad \mbox{or}  \quad
\alpha < 1-1/t   \quad
& \mbox{if} & \quad 1 < t  < \infty\, ,
\\
&& \\
(\alpha = 1  \quad \mbox{and}  \quad \sigma \ge 0)  \quad \mbox{or} \quad
\alpha < 1  \quad & \mbox{if} & \quad t  =\infty\,,
\end{array}
\right.
\]

\begin{Lem}\label{sharp}
{\rm (i)} Let $0 <p \le \infty$ and $1 < q \le \infty$.
Then $f_{\alpha, \sigma}$
belongs to  $R{B}^{d/p}_{p,q} (\Rd)$ if, and only if
$(\alpha,\sigma) \in U_q$.
\\
{\rm (ii) }
Let $1 < p < \infty$.
Then $f_{\alpha, \sigma}$
belongs to  $R{F}^{d/p}_{p,q} (\Rd)$ if, and only if
$(\alpha,\sigma) \in U_p$.
\end{Lem}

\begin{Rem}{\rm
Let us mention that in    \cite{Bou2} the result is stated for
$p\ge 1$ only. However, the proof extends to $p<1$ nearly without changes
(in his argument which follows formula (9) in \cite{Bou2}
one has to choose $k> d/(2p)$).
}
\end{Rem}


 \section{Decay properties of radial functions -- proofs}
 \label{inhom}



 \subsection{Proof of Theorem  \ref{decay4}}


{\em Step 1.} Proof of (i).
Following Remark \ref{border} it will be enough to prove the decay
estimate (\ref{eq-72c}) for $RB^{1/p}_{p,1} (\Rd)$, $0< p < \infty$,
and for $RF^{1/p}_{p,\infty} (\Rd)$, $0 < p \le 1$.
A proof in  case $RB^{1/p}_{p,1} (\Rd)$ has been given in
\cite{SS}.
So we are left with the proof for the Lizorkin-Triebel spaces.
We will follow the ideas of the proof of Cor. \ref{uc}. Let $f \in RF^{1/p}_{p,\infty}(\R^d)$.
Let
\[
f = \sum_{j=0}^\infty \,  s_{j,0}\, a_{j,0} +
\sum_{j=0}^\infty \, \sum_{k=1}^\infty \,
\sum_{\ell = 1}^{C_{d,k}} s_{j,k}\, a_{j,k,\ell},
\]
be an atomic decomposition such that
$\|s_{j,k}|f_{p,\infty ,d}\|\asymp \| \, f \, | F^{1/p}_{p,\infty}(\R^d)\| $.
We fix $x$, $|x|>1$. Observe, that for all $j\ge 0$ there exists $k_j\ge 1$ such that
\be\label{e1}
k_j \, 2^{-j} \le |x| < (k_j + 1)\, 2^{-j} \, .
\ee
Then the main part of $f$ near $(|x|, 0, \ldots \, ,0)$ is given by the function
\be\label{mainf}
f^M (y) = \sum_{j=0}^\infty  \, s_{j,k_j}\,
a_{j,k_j,0}(y) \, , \qquad y \in \Rd\, ,
\ee
(in fact, $f$ is a finite sum of functions of type
\[
\sum_{j=0}^\infty  \, s_{j,k_j+r_j}\,
a_{j,k_j+r_j, t_j}(y) \, ,
\]
and $|r_j|$ and $|t_j|$ are uniformly bounded). For convenience we
shall derive an estimate of the main part $f^M$ only.
Because of (\ref{e1}) and the normalization of the atoms
we obtain
\be\label{e2}
|f^M (y)| \, \lsim \,  \sum_{j=0}^\infty |s_{j,k_j}|\, 2^{j\frac{d-1}{p}}\,
\lsim |x|^{\frac{1-d}{p}} \, \Big(\sum_{j=0}^\infty |s_{j,k_j}|^p\, k_j^{d-1} \Big)^{1/p}\, ,
\ee
since $p\le 1$. On the other hand
\begin{eqnarray}
\nonumber \| \,s \, | f_{p,\infty,d}\| & = &
\| \, \sup_{j=0,1,\ldots } \,
\sup_{k \in \N_0}  \, |s_{j,k}| \, 2^{\frac{jd}{p}}\, \widetilde{\chi}_{j,k}(\cdot)\, |L_p(\Rd)\|
\\
&\ge&  \| \, \sup_{j=0,1,\ldots }\, |s_{j,k_j}| \, 2^{\frac{jd}{p}}\, \widetilde{\chi}_{j,k_j}(\cdot)\, |L_p(\Rd)\|\, .
\label{e3ls}
\end{eqnarray}
Using $P_{j+1,k_{j+1}} \subset P_{j,k_{j}}$ we obtain the identity
\[
\sup_{j}\, |s_{j,k_j}| \, 2^{\frac{jd}{p}}\, \widetilde{\chi}_{j,k_j}(\cdot)\, = \, \sum_{j=0}^\infty \max_{i=0,\ldots ,j}|s_{i,k_i}| 
\big( \widetilde{\chi}_{j,k_j}(\cdot) - \widetilde{\chi}_{j+1,k_{j+1}}(\cdot)\big) \, .
\]
By the pairwise disjointness of the sets $P_{j,k_{j}} \setminus P_{j+1,k_{j+1}}$
this implies
\begin{eqnarray}\label{e3als}
\|\, \sup_{j=0,1,\ldots }\, |s_{j,k_j}| \, 2^{\frac{jd}{p}}\, \widetilde{\chi}_{j,k_j}(\cdot)\, |L_p(\Rd)\|\,\asymp \,
\bigg(\sum_{j=0}^\infty \max_{i=0,\ldots ,j}\big(|s_{i,k_i}|^p 2^{id}\big)\, 2^{-jd}k_j^{d-1}\bigg)^{1/p}\, .
\end{eqnarray}
Obviously
\begin{equation}
\Big(\sum_{j=0}^\infty |s_{j,k_j}|^p\, k_j^{d-1} \Big)^{1/p}\, \le \, \bigg(\sum_{j=0}^\infty \max_{i=0,\ldots ,j}\big(|s_{i,k_i}|^p 2^{id}\big)\, 2^{-jd}k_j^{d-1}\bigg)^{1/p}\, .
\label{e4}
\end{equation}
Combining  \reff{e2} - \reff{e4} we have proved (\ref{eq-72c})
in case of  Lizorkin-Triebel spaces.
\\
{\em Step 2.} Proof of (ii).
As in Step 1 it will be sufficient to deal with the limiting cases.
\\
{\em Substep 2.1.} Let $f \in RF^{1/p}_{p,\infty}(\R^d)$.
Let $B_r(0)$ be the ball in $\Rd$ with center in the origin and radius $r$.
Then \reff{e2}-\reff{e4} yield
\beqq
|f^M (x)| \, & \lsim & \, |x|^{\frac{1-d}{p}}\,
\| \, \sup_{j=0,1,\ldots }\, |s_{j,k_j}| \, 2^{\frac{jd}{p}}\, \widetilde{\chi}_{j,k_j}(\cdot)\, |L_p(\Rd \setminus B_r(0))\|\,
\\
& \lsim & \, |x|^{\frac{1-d}{p}}\,
\| \, \sup_{j=0,1,\ldots }\, \sup_{k \in \N_0} |s_{j,k}| \, 2^{\frac{jd}{p}}\, \widetilde{\chi}_{j,k}(\cdot)\, |L_p(\Rd \setminus B_r(0))\|\,
\eeqq
where $r=|x|-18>0$, see property (b) of the covering $(\Omega_{j,k,\ell})$
in Subsection \ref{rsa}.
In view of this inequality an application of Lebesgue's theorem on dominated convergences proves
(\ref{eq-91b}).
\\
{\em Substep 2.2.} Let $f \in RB^{1/p}_{p,1}(\R^d)$.
We argue as in Substep 2.1 by using the notation from Step 1.
Since
\[
\lim_{r \to \infty}\quad
\sum_{j=0}^\infty \Big(\sum_{k\ge r} |s_{j,k}|^p \, (1+k)^{d-1}\Big)^{1/p}
=0
\]
we conclude from (\ref{e2}) that (\ref{eq-91b}) holds in this case as well.
\\
{\em Step 3.} Proof of (iii). We shall use the test functions constructed in Lemma
\ref{test}. We choose $s_0>\max\{\sigma_p(d), s\}$. For simplicity we consider  $|x|=2^r$ with   $r \in \N$. We choose $\lambda$ s.t. $|x|=(1+\lambda)/2$. Hence $f_{1,\lambda} (x) =1$.
This implies
\[
|x|^{\frac{d-1}{p}}  \, | 2^{-r(d-1)/p}  \, f_{1,\lambda} (x)| = \, 1 \, ,
\]
and
\[
\|\, 2^{-r(d-1)/p}  \, f_{1,\lambda} \, |\aspq (\Rd)\|\, \lsim \, \|\, 2^{-r(d-1)/p}  \, f_{1,\lambda} \, |B^{s_o}_{p,q} (\Rd)\|
\asymp \,  1\, ,
\]
see Rem.~\ref{testf}, which proves the claim.
\\
{\em Step 4.} Proof of (iv).
It will be enough to study the case $s=1/p$.
\\
{\em Substep 4.1} Let  $q>1$.
According to Lemma \ref{sharp}(i) there exists a compactly supported function $g_0$ which belongs to $RB^{1/p}_{p,q} (\R)$ and  is unbounded near the origin.
By multiplying with a smooth cut-off function if necessary
we can make the support of this functions as small as we want.
For the given sequence $(x^j)_j$ we define
\[
g(t):= \sum_{j=1}^\infty \frac{1}{\max(|x^j|,j)^\alpha} \, g_0(t-|x^j|)\, , \qquad t \in \R \, ,
\]
where we will choose $\alpha >0$ in dependence on $p$.
The function $g$ is unbounded near $|x^j|$ and by means of the translation invariance of the Besov spaces $B^{1/p}_{p,q} (\R)$ we obtain
\[
\| \, g \, |B^{1/p}_{p,q} (\R)\|^{\min(1,p)}
\le \| \, g_0\, | B^{1/p}_{p,q} (\R)\|^{\min(1,p)}
\sum_{j=1}^\infty j^{-\alpha \, \min (1,p)} \lsim \, \| \, g_0\, | B^{1/p}_{p,q} (\R)\|^{\min(1,p)}\, ,
\]
if $\alpha\, \cdot\, \min (1,p) >1$.
We employ Rem.~\ref{equi} (Cor.~\ref{ok})
with respect to each summand. This yields
\beqq
& & \hspace*{-0.8cm}
\| \, \ext g \, |B^{1/p}_{p,q} (\Rd)\|^{\min(1,p)}
\\
& \le &
\sum_{j=1}^\infty \max(|x^j|, j)^{-\alpha \, \min (1,p)} \, \| \, \ext (g_0 (\, \cdot \, -|x^j|))\, | B^{1/p}_{p,q} (\Rd)\|^{\min(1,p)}
\\
& \lsim & \, \| \, g_0\, | B^{1/p}_{p,q} (\R)\|^{\min(1,p)}\,
\sum_{j=1}^\infty \Big(\max (|x^j|, j)^{-\alpha} \, |x^j|^{(d-1)/p}\Big)^{\min (1,p)}
\\
& \lsim & \, \| \, g_0\, | B^{1/p}_{p,q} (\R)\|^{\min(1,p)}\,
\, ,
\eeqq
if $(\alpha - (d-1)/p)\, \min (1,p) >1$.
\\
{\em Substep 4.2.} We turn to the $F$-case. It will be enough to study the situation
$s=1/p$ and $1<p<\infty$.
We argue as above using this time Lemma \ref{sharp}(ii).
An application of Rem. \ref{equi} (Cor. \ref{ok}) yields the result
but with the extra condition $1/p>\sigma_{q}(d)$.
\epr


\subsection{Traces of $BV$-functions and consequences
for the decay}
\label{spurbv}


We recall a definition of the space $BV(\Rd)$, $d\ge 2$, which will be convenient for us,
see \cite[5.1]{EG} or \cite[5.1]{Zi}.

\begin{Def}\label{BVDef1}
Let $g\in L_1(\R^d)$. We say, that $g\in BV(\R^d)$ if for every $i=1,\dots, d$ there
is a  signed Radon measure $\mu_i$ of finite total variation such that
$$
\int_{\R^d} g(x)\frac{\partial}{\partial x_i}\phi(x) dx
=-\int_{\R^d}\phi(x)d\mu_i(x),\qquad \phi\in C^1_c(\R^d),
$$
where $C^1_c(\R^d)$ denotes the set of all continuously differentiable functions
on $\R^d$ with compact support.
The space $BV(\R^d)$ is equipped with the norm
\[
\|\, g\, |BV(\R^d)\|=\|\, g \, |L_1(\R^d)\| +
\sum_{i=1}^d \|\, \mu_i\, |{\mathcal M}\|,
\]
where $\|\, \mu_i\, |{\mathcal M}\|$ is the total variation of $\mu_i$.
\end{Def}


\subsubsection{Proof of Theorem \ref{BVThm1}}
\label{spurbv1}


We need some preparations. Recall, the space $C_c^1([0,\infty))$
has been defined in Definition \ref{BVDef2}.
By $\omega_{d-1}$ we denote the surface area of the unit sphere in $\R^d$
and by $\sigma$ the $(d-1)$-dimensional Hausdorff
measure in $\R^d$, i.e. $\omega_{d-1} :=\sigma(\{x\in\R^d: |x|=1\}).$
As above $r= r(x):= |x|$.

\begin{Lem}\label{BV:Lem1}
(i) If $\varphi\in C_c^1([0,\infty))$, then all the functions
\begin{equation}\label{BV:eq:1}
\phi_i(x):= \left\{\begin{array}{lll}
\varphi(r(x))\cdot\frac{x_i}{r(x)} &\quad  & \mbox{if} \quad x=(x_1,\dots,x_d)\neq 0\, ,
\\
&& \\
0 && \mbox{if} \quad x=0\, ,
\end{array}
\right.
\end{equation}
$i=1,\dots,d$, belong to $C_c^1(\R^d)$.
\\
(ii) If $\phi\in C_c^1(\R^d)$, then all functions
\begin{equation}\label{BV:eq:2}
\varphi_i(t) :=  \left\{\begin{array}{lll}
\frac{1}{\omega_{d-1} \, t^{d-1}}\, \int_{|x|=t} \, \phi(x)\cdot\frac{x_i}{r(x)}\, d\sigma(x)
&\qquad & \mbox{if}\quad t>0\, ,
\\
&& \\
0 && \mbox{if} \quad t=0\, ,
 \end{array}\right.
\end{equation}
$i=1,\dots,d$, belong to $C_c^1([0,\infty))$.
\end{Lem}

\bpr
{\em Step 1.} Proof of (i). Under the given assumption we immediately get
$\phi_i\in C^1(\R^d\setminus\{0\})$ and $\supp \phi_i$ is compact. Hence, we have to study the regularity properties in the origin.
Obviously, $\phi_i(0)=0$ and $\displaystyle \lim_{x\to 0}\phi_i(x)=0.$
We claim that
$$
\frac{\partial\phi_i}{\partial x_j}(0)=\begin{cases}
\varphi'(0)\quad &\text{if}\ i=j,\\
0\quad &\text{otherwise}.
\end{cases}
$$
Let $e_1, \ldots\, e_d$ denote  the elements of the canonical basis of $\Rd$.
If $i=j$, then
\[
\lim_{t\to 0} \, \frac{\phi_i(te_i)}{t}=
\lim_{t\to 0} \, \frac{\varphi(|t|)}{|t|} = \varphi'(0) \, .
\]
Hence $\displaystyle \frac{\partial\phi_i}{\partial x_i}(0)=\varphi'(0).$
The cases $i\not=j$ are obvious.
\\
Next, we show, that the functions
$\displaystyle \frac{\partial \phi_i}{\partial x_j}$
are continuous in the origin. To begin with we investigate  the case  $i=j$.
Then
\begin{align*}
\left|\frac{\partial \phi_i}{\partial x_i}(x)-\varphi'(0)\right|
&=\left|\varphi'(r(x))\cdot \frac{x_i^2}{r^2(x)}+ (\varphi(r(x))-\varphi(0))
\cdot \frac{r^2(x)-x_i^2}{r^3(x)}-\varphi'(0)\right|\\
&=\left|\varphi'(r(x))\cdot \frac{x_i^2}{r^2(x)}+\varphi'(\theta r(x))
\cdot \frac{r^2(x)-x_i^2}{r^2(x)}-\varphi'(0)\right|,
\end{align*}
where we have used the Mean Value Theorem with a suitable $0< \theta = \theta(x)<1$. The continuity
of $\varphi'$ implies, that the  expression on the right-hand side
tends to zero if $x\to 0$. If $i\not=j$, we write
\begin{align}
\frac{\partial \phi_i}{\partial x_j}(x) &= \varphi'(r(x))\frac{x_ix_j}{r^2(x)}
-\varphi(r(x))\frac{x_ix_j}{r^3(x)}
\nonumber
\\
&=\varphi'(r(x))\frac{x_ix_j}{r^2(x)}+(\varphi(0)-\varphi(r(x)))\frac{x_ix_j}{r^3(x)}
\nonumber
\\
&=\frac{x_ix_j}{r^2(x)}\left[\varphi'(r(x))-\varphi'(\theta r(x))\right]\, ,
\label{f-05}
\end{align}
with some  $0<\theta<1$. Again the continuity
of $\varphi'$ implies, that the  expression in (\ref{f-05})
tends to zero as $x\to 0$. Hence, $\phi_i \in C_c^1 (\Rd)$.
\\
{\em Step 2.} Proof of (ii).
The regularity and support properties of $\varphi_i$
on $(0,\infty)$ are obvious.
Hence, we are left with the study of the behaviour near $0$.
We shall use the identity
$$
\int_{|x|=t}\, \frac{x_i}{r(x)} \, d\sigma(x) = 0\, , \qquad t>0 \, .
$$
In case $t>0$ this leads to the estimate
\beqq
|\varphi_i(t)| & \le &  \frac{1}{\omega_{d-1}t^{d-1}} \Big(
\left|\int_{|x|=t}\, \phi(0)\, \frac{x_i}{r(x)}d\sigma(x)
\right|
+ \left|\int_{|x|=t}\, (\phi(x)-\phi(0))\, \frac{x_i}{r(x)}\,
d\sigma(x)\right|\Big)
\\
&\le & 0+\sup_{|x|=t}|\phi(x)-\phi(0)|.
\eeqq
Hence, $\varphi_i(t)$ tends to $0$ if $t\to 0^+$.
Furthermore,
\begin{align*}
 \frac{\varphi_i(t)}{t}
& = \frac{1}{\omega_{d-1}\,t^d}\, \int_{|x|=t} \, (\phi(x)-\phi(0))\, \frac{x_i}{r(x)}\, d\sigma(x)
\\
&= \frac{1}{\omega_{d-1}\, t^{d+1}}\, \int_{|x|=t}\, (\nabla\phi(\eta_x x)\cdot x) \, x_i \, d\sigma(x)
\\
& = \frac{1}{\omega_{d-1}\, t^{d+1}}\, \sum_{j=1}^d\, \int_{|x|=t}
\left(\frac{\partial\phi}{\partial x_j}(\eta_x x)-\frac{\partial\phi}{\partial x_j}(0)\right)\, x_j\, x_i\, d\sigma(x)
\\
& \quad \hspace{3cm} + \quad \frac{1}{\omega_{d-1}\, t^{d+1}}\, \sum_{j=1}^d\, \frac{\partial\phi}{\partial x_j}(0)\, \int_{|x|=t}
\, x_j\, x_i \, d\sigma(x)
\end{align*}
with some  $0<\eta_x<1$. The first term on the right-hand side  tends always to zero as $t\to 0^+$ (since $\displaystyle \frac{\partial\phi}{\partial x_j}$ is continuous). From the  second term those summands with  $i\not=j$ are vanishing for all $t$. If $i=j$, then the integrand is homogeneous.
We obtain, by taking the limit with respect to $t$,
\begin{equation}\label{BV:eq:3}
\varphi_i'(0)= \lim_{t \to 0^+} \,  \frac{\varphi_i(t)}{t} = \frac{1}{\omega_{d-1}}\, \frac{\partial\phi}{\partial x_i}(0)\, \int_{|y|=1}\, y_i^2 \, d\sigma(y) \, .
\end{equation}
It remains to check the limit of $\varphi_i'(t)$ if $t$ tends to $0^+$.
Observe
\begin{align*}
\varphi_i'(t)= \frac{1}{\omega_{d-1}}\, \frac{d}{ dt}\int_{|y|=1}\, \phi(ty) \, y_i \, d\sigma(y)
= \frac{1}{\omega_{d-1}}\, \sum_{j=1}^d\, \int_{|y|=1}
\, \frac{\partial \phi}{\partial y_j}(ty)\, y_j\, y_i\,  d\sigma(y)\, .
\end{align*}
Since $\displaystyle \int_{|y|=1}\, y_j\, y_i\,  d\sigma(y) =0$ if $i \neq j$
those summands (with $i \neq j$) tend to $0$ if  $t$ tends to $0^+$.
Hence
\[
\lim_{t \to 0^+} \, \varphi_i' (t) = \lim_{t \to 0^+}
\frac{1}{\omega_{d-1}}\, \int_{|y|=1}
\, \frac{\partial \phi}{\partial y_i}(ty)\,  y_i^2\,  d\sigma(y) =
\varphi'(0)\, ,
\]
see \eqref{BV:eq:3}. The proof is complete.
\epr

\noindent
{\bf Proof of Theorem \ref{BVThm1}.}\\
{\em Step 1.} Let $f (x)= g(|x|)\in BV(\R^d)$. We claim that $g\in BV(\R^+,t^{d-1})$.
Let $\mu_1,\dots,\mu_d$ denote the corresponding signed Radon measures according to
Definition \ref{BVDef1}. By means of
\[
d\nu := \sum_{i=1}^d \, \frac{x_i}{r(x)}\, d\mu_i
\]
we define the measure $\nu$ on $\R^d$.
Since  $g(r(x))$ is radial we  conclude  $\mu_i(\{0\})=0$, $i=1,\ldots \, ,d$.
Hence the measure $d\nu$ is well defined. In addition we  introduce a
measure $\nu^+$ on $\R^+$ by
$$
\omega_{d-1}\, \int_{A} \, t^{d-1} \, d\nu^+(t):= \int_{\{|x|\in A\}}\, d\nu(x),
$$
for any Lebesgue measurable subset $A\subset \R^+$.
We fix $\varphi\in C_c^1([0,\infty))$. Since
$\varphi(r(x))\frac{x_i}{r(x)}\in C_c^1(\R^d)$, cf. Lemma \ref{BV:Lem1}(i), we  calculate
\begin{eqnarray*}
\omega_{d-1}\, \int_{0}^\infty \, g(t)\,  && \hspace{-0.8cm}
[\varphi(s)\, s^{d-1}]'(t)\, dt
=\omega_{d-1}\, \int_{0}^\infty \, t^{d-1}\, g(t) \left[\varphi'(t)+\frac{d-1}{t}\varphi(t)\right]\, dt
\\
& = & \int_{\R^d}\, g(r(x))\, \left[\varphi'(r(x))+\varphi(r(x))\cdot\frac{d-1}{r(x)}\right]\, dx
\\
&= & \sum_{i=1}^d\, \int_{\R^d}g(r(x))\, \left[\varphi'(r(x))\frac{x_i^2}{r^2(x)}+\varphi(r(x))\cdot\frac{r^2(x)-x_i^2}{r^3(x)}\right]\, dx
\end{eqnarray*}
\begin{eqnarray*}
\qquad & = &\sum_{i=1}^d\, \int_{\R^d}\, g(r(x))\, \frac{\partial}{\partial x_i}\, \left[\varphi(r(x))\cdot\frac{x_i}{r(x)}\right]\, dx
\\
&= & -\sum_{i=1}^d\, \int_{\R^d}\, \varphi(r(x)) \, \frac{x_i}{r(x)}\, d\mu_i=
-\int_{\R^d}\, \varphi(r(x))\, d\nu(x)\\
&= & -\omega_{d-1}\, \int_{0}^\infty t^{d-1}\, \varphi(t)\, d\nu^+(t)\, .
\end{eqnarray*}
This proves (\ref{f-06}).
Moreover, (\ref{f-07}) follows from
\[
\|\,g(r(x))\, |L_1(\R^d)\| = \omega_{d-1}\, \|\, g\, |L_1(\R,|t|^{d-1})\|
\]
and
\begin{align*}
\omega_{d-1}\, \int_{0}^\infty\,  t^{d-1} \, d|\nu^+|(t) &= \, \int_{\R^d} \, d|\nu|(x)
\le
\sum_{i=1}^d \, \int_{\R^d}\, \frac{|x_i|}{r(x)}\, d|\mu_i|
\\
& \le \sum_{i=1}^d \, \int_{\R^d}\, d|\mu_i|
\le \| \, g(r(x))\, |BV(\R^d)\|\, .
\end{align*}
{\em Step 2.} Let $g$ be a function in $BV(\R^+,t^{d-1})$. We claim, that $g(r(x))\in BV(\R^d)$.
Let $\nu^+$ be the signed Radon measure associated to $g$ according to (\ref{f-06}).
We define
$$
 \nu(A):=\int_0^\infty \, \sigma(\{x:|x|=t\}\cap A)\, d\nu^+(t)
$$
for any Lebesgue measurable set $A\subset \R^d$.
Further we put $\mu_i :=\frac{x_i}{r(x)}\nu$, $i=1, \ldots \, d$.
Let $\chi_A$ denote the characteristic function of $A$.
Then
\[
\nu (A) = \int_{\R^d} \, \chi_A(x)\,  d\nu(x)= \int_0^\infty \Big[\int_{|x|=t}\, \chi_A(x)\, d\sigma(x)\Big]\, d\nu^+(t)
\]
and this identity can be extended to
\[
\int_{\R^d} \, \phi(x) \, d\nu(x)= \int_0^\infty \Big[\int_{|x|=t}\, \phi(x)\, d\sigma(x)\Big]\, d\nu^+(t)\, ,\qquad \phi\in L_1(\R^d)\, ,
\]
by using some standard arguments.
Next we want to show, that $\mu_i, i=1,\dots, d$, are the weak derivatives of $g(r(x))$. Let $\phi\in C_c^1({\R^d})$ and let $\varphi_i$ be the associated
functions, see (\ref{BV:eq:2}).
According to Lemma \ref{BV:Lem1} (ii) we know that $\varphi_i \in C_c^1([0,\infty))$.
Using the normalized outer  normal with respect to the surface $\{x:|x|=T\}$, which is obviously given by
\[
n(x) = (n_1(x),\dots,n_d(x)) = \frac{1}{r(x)} \, \left(x_1,\dots,x_d\right)\, ,
\]
and the Gauss Theorem, we obtain
\begin{align*}
-\varphi_i(T)\, \omega_{d-1}\, T^{d-1} &= -\int_{|x|=T}\, \phi(x)\, \frac{x_i}{r(x)}\, d\sigma(x)
=-\int_{|x|=T}\, \phi(x)\, n_i(x)\, d\sigma(x)\\
&=-\int_{|x|\le T}\, \frac{\partial \phi}{\partial x_i}(x)\, dx
=\int_{|x|\ge T}\, \frac{\partial \phi}{\partial x_i}(x)\, dx
\\
&=\int_T^\infty \Big[
\int_{|x|=t}\, \frac{\partial \phi}{\partial x_i}(x) \, d\sigma(x)\Big] \, dt.
\end{align*}
Hence
$$
\omega_{d-1}\, \left[\varphi_i(t)\, t^{d-1}\right]'(T)=\int_{|x|=T}\, \frac{\partial \phi}{\partial x_i}(x)\, d\sigma(x)\, ,\qquad T>0.
$$
This formula justifies the identity
\beqq
\int_{\R^d}g(r(x))\, \frac{\partial \phi(x)}{\partial x_i}\, dx & = &
\int_{0}^\infty g(t)\, \int_{|x|=t} \frac{\partial \phi(x)}{\partial x_i}\, d\sigma (x) \, dt
\\
& = &
\omega_{d-1}\, \int_0^\infty g(t)[\varphi_i(s)\, s^{d-1}]'(t)dt \, .
\eeqq
Next we use $g \in BV(\R^+,t^{d-1})$. This implies
\begin{align*}
\int_{\R^d}g(r(x))\frac{\partial \phi(x)}{\partial x_i}dx
& = - \omega_{d-1}\, \int_0^\infty \varphi_i(t)\, t^{d-1} \, d\nu^+(t)
\\
&=-\int_0^\infty \int_{|x|=t}\, \phi(x)\cdot\frac{x_i}{r(x)}\,
d\sigma(x)\,  d\nu^+(t)\\
&=-\int_{\R^d}\, \phi(x) \, \frac{x_i}{r(x)}\, d\nu(x)=-\int_{\R^d}\, \phi(x)\, d\mu_i(x)\, ,
\end{align*}
which proves that the  $\mu_i$ are the weak derivatives of $g(r(x))$.
\\
It remains to prove the estimates for the related norms.
The required  estimate follows easily by
\beqq
\int_{\R^d}\,d|\mu_i| & = & \int_{\R^d}\, \frac{|x_i|}{r(x)}\, d|\nu|(x)
= \int_0^\infty \Big[\int_{|x|=t}\, \frac{|x_i|}{r(x)}\, d\sigma(x)\Big]\, d|\nu^+|(t)
\\
&\le& \omega_{d-1}\, \int_0^\infty\,  t^{d-1}\, d|\nu^+|(t)\, .
\eeqq
The proof is complete.
\epr


\subsubsection{Proof of Theorem \ref{decaybv}}
\label{spurbv2}


Recall, that we will work with the particular representative  $\tilde f$ of the equivalence class $[f]$,
see Remark \ref{repre}.
For convenience we will drop the tilde.
We shall apply standard mollifiers.
Let $\varphi \in C_0^\infty (\R)$ be a function such that $\varphi \ge 0$,
$\supp \varphi \subset [0,1]$,
and $\int \varphi (t)\, dt =1 $.
For $R>0$ and $\varepsilon >0$ we define
\be\label{f-08}
\varphi_\varepsilon (t):= \varepsilon^{-1}\, \int_{R}^\infty \, \varphi
\Big(\frac{t-y}{\varepsilon}\Big)\, dy =
\int_{-\infty}^{\frac{t-R}{\varepsilon}} \, \varphi (z)\, dz
\ee
(which is nothing but the mollification of the characteristic function of the interval $(R,\infty)$).
In addition we need a cut-off function. Let $\eta \in C_0^\infty (\R)$
s.t. $\eta (t)=1$ if $|t|\le 1$ and $\eta (t)=0$ if $|t| \ge 2$.
For $M \ge 1$ we define $\eta_M (t):= \eta (t/M)$, $t \in \R$.
It is easily checked that the functions
\[
\phi_{M,\varepsilon} (t):= t^{1-d} \, \varphi_\varepsilon (t)\,
\eta_M (t)\, , \qquad t \in \R\, ,
\]
belong to $ C^1_c([0,\infty))$.
For $g \in BV(\R^+,t^{d-1})$ this implies
\be\label{f-10}
\int_0^\infty g(t) \, [\phi_{M,\varepsilon} (s) \, s^{d-1}]'(t)\, dt
= - \int_0^\infty  \, \varphi_{\varepsilon} (t)\, \eta_M (t) \, d\nu^+(t)\, ,
\ee
see (\ref{f-06}).
Since for $M > R + \varepsilon$ we have
\beqq
\int_0^\infty g(t) \, [\phi_{M,\varepsilon} (s) \, s^{d-1}]'(t)\, dt & = &
\int_{R}^{R + \varepsilon} g(t) \,
\varepsilon^{-1}\, \varphi \Big(\frac{t-R}{\varepsilon}\Big) \,  dt
\\
&& \hspace{1cm}
+ M^{-1}\,  \int_M^\infty g(t)\, \varphi_\varepsilon (t) \, \eta ' (t/M)
\, dt
\eeqq
and
\[
\lim_{M \to \infty} M^{-1}\,  \int_M^\infty g(t)\, \varphi_\varepsilon (t)
\, \eta ' (t/M) \, dt =0
\]
($g \in L_1 (\R^+,t^{d-1})$),
we get
\be\label{f-09}
\lim_{\varepsilon \downarrow 0, M \to \infty} \, \int_0^\infty g(t) \, [\phi_{M,\varepsilon} (s) \, s^{d-1}]'(t)\, dt  = g(R)\, ,
\ee
if $R$ is a Lebesgue point of $g$. But
\[
 \Big|\, \int_0^\infty  \, \varphi_{\varepsilon} (t)\, \eta_M (t) \, d\nu^+(t)\Big|
\le \int_R^{\infty}  \, \varphi_{\varepsilon} (t)\, d|\nu^+|(t)
\le R^{1-d}\, \int_R^{\infty}  \, t^{d-1}\, d|\nu^+|(t) \, .
\]
Combining (\ref{f-09}), (\ref{f-10}) with this estimate  we have proved
(\ref{ws-01}) and (\ref{ws-02}) simultaneously.
\epr


\subsection{Proof of the assertions in
Subsection \ref{main3}}



 \subsubsection{Proof of Lemma  \ref{bi}}


Sufficiency of the conditions is obvious, see e.g. \cite{SiTr}.
Necessity follows from the examples investigated in Lemma \ref{sharp}.
\epr


 \subsubsection{Proof of Theorem  \ref{decay2}}


We argue by using the atomic characterizations in Subsection
\ref{rsa}.
\\
{\em Step 1.} Proof of (i).
For simplicity let $|x|=2^{-r}$, $r\in \N$.
If $y$ satisfies $2^{-r-1} \le |y| \le 2^{-r+1}$, then, using  the support condition for atoms,
we know that
$f$ allows an optimal atomic decomposition such that
\be\label{f-01}
f(y) = \sum_{j=0}^\infty \sum_{k=\max(0,[2^{j-r}]-n_0)}^{[2^{j-r}]+n_0} \sum_{\ell =1}^{C(d,k)}
s_{j,k}\, a_{j,k,\ell} (y) \, .
\ee
Here $n_0$ is a general natural number depending on the decomposition
$\Omega$, but not on $r$.
\\
{\em Substep 1.1.} We assume first that $\frac 1 p < s < \frac d p$.
From the $L_\infty$-estimate of the atoms, property (f) of the coverings $(\Omega_{j,k,\ell})_{j,k,\ell}$, and
the inequality (\ref{atom3}) we derive
\beqq
&&
\hspace*{-1cm}
|f(y_1,0, \ldots , 0)|  \le   \, \sum_{j=0}^\infty \sum_{k=\max(0,[2^{j-r}]-n_0)}^{[2^{j-r}]+n_0} \sum_{\ell =1}^{K}
|s_{j,k}|\, |a_{j,k,\ell} (y_1,0, \ldots , 0)|
\\
& \lsim &  \Big( \sum_{j=0}^{r+n_1} \sum_{k=0}^{n_2} |s_{j,k}|\,
2^{-j (s-d/p)} +   \sum_{j=r+n_1+1}^{\infty} \sum_{k=2^{j-r}-n_0}^{2^{j-r}+n_0} |s_{j,k}|\, 2^{-j (s-d/p)}\Big)
\\
& \lsim & \, \| \, f\, |RB^s_{p,\infty} (\Rd)\|\,
\Big( \sum_{j=0}^{r+n_1} 2^{-j (s-d/p)} +
\sum_{j=r+n_1 + 1}^{\infty} 2^{-(j-r)(d-1)/p}\,
2^{-j (s-d/p)}\Big)
\\
& \lsim & 2^{r(\frac dp-s)}\, \| \, f\, |RB^s_{p,\infty} (\Rd)\|\,
\, ,
\eeqq
for appropriate natural numbers $n_1, n_2$ (independent of $r$).
For the last two steps of the estimate we used $1/p < s < d/p$.
Taking into account the elementary embedding
$\aspq (\Rd) \hookrightarrow B^s_{p,\infty} (\Rd)$
we obtain the inequality (\ref{eq-72}).
\\
{\em Substep 1.2.} Now, let $s=\frac 1 p$. For the Besov spaces $RB^{1/p}_{p,1} (\Rd)$ the inequality (\ref{eq-72}) was proved in \cite{SS}.  So it remains to consider the Lizorkin -Triebel spaces  $RF^{1/p}_{p,\infty} (\Rd)$ with $0< p \le 1$. For simplicity we regard  the main part of $f\in RF^{1/p}_{p,\infty} (\Rd)$, cf. \reff{f-01} and compare with (\ref{mainf}).
Now $k_j=1$ if $0\le j\le r+ n_1$ and $k_j\sim |x|2^j$ if $j > r+ n_1$.
Hence, we obtain
\beqq
&&
\hspace*{-1cm}
|f^M(|x|,0, \ldots , 0)|  \le   \, \sum_{j=0}^\infty 
|s_{j,k_j}|\, |a_{j,k_j,0} (y)| 
\,\lsim\, \Big( \sum_{j=0}^{\infty} 
|s_{j,k_j}|^p\,
2^{j (d-1)} \Big)^{1/p}
\\
&\lsim &  |x|^{\frac{1-d}{p}}\, \Big( \sum_{j=0}^{r+n_1} 
|s_{j,k_j}|^p\, 
+  \, \sum_{j=r+n_1+1}^{\infty} |s_{j,k}|^p\, k_j^{d-1}\Big)^{1/p}
\, \lsim \, |x|^{\frac{1-d}{p}} \, \| \, f\, |RF^{1/p}_{p,\infty} (\Rd)\|\,
\ ,
\eeqq
where we used  $p\le 1$  for the second step and \reff{e3ls}-\reff{e3als} for the last one.
\\
{\em Step 2.} Proof of (ii).
By using elementary embeddings it will be enough to prove
(\ref{eq-72b}) with $R\aspq (\Rd) = R\bspq(\Rd)$ and $q$ small.
Again we concentrate on $|x|=2^{-r}$, $r\in \N$.
Our test function is taken to be $2^{-r(s-d/p)} \, f_{j,\lambda}$, see Lemma \ref{test},
where we choose $j:= 2+ r$ and $\lambda := 3$.
Then it follows
\[
\| \, 2^{-r(s-d/p)}\, f_{j,\lambda}\, |B^s_{p,q}(\Rd)\| \asymp 1\,
\qquad \mbox{and}\qquad 2^{-r(s-d/p)} \,  f_{j,\lambda} (x) = 2^{-r(s-d/p)} = |x|^{s-d/p}\, .
\]
The proof is complete.
\epr


 \subsubsection{Proof of Lemma \ref{lim2}}


The arguments are the same as in proof of Theorem \ref{decay4}(iv).
\epr


 \subsection{Proof of Theorem  \ref{lim1}}


We shall use the notation from the proof of Theorem \ref{decay2}, Step 1.
Again we employ the formula (\ref{f-01}) and obtain
\[
|f(y_1,0, \ldots , 0)|   \lsim   \Big( \sum_{j=0}^{r+n_1} \sum_{k=0}^{n_2} |s_{j,k}|  +   \sum_{j=r+n_1+1}^{\infty} \sum_{k=2^{j-r}-n_0}^{2^{j-r}+n_0} |s_{j,k}|\Big)
\]
since $s=d/p$. \\
{\em Step 1.} Proof of (i). We shall use the standard abbreviation
$q' := q/(q-1)$.
Since $q>1$ we can use H\"older's inequality and conclude
\be\label{f-02}
\sum_{j=0}^{r+n_1} \sum_{k=0}^{n_2} |s_{j,k}| \lsim r^{1/q'} \,
\Big(\sum_{j=0}^{r+n_1} \sum_{k=0}^{n_2} |s_{j,k}|^q \Big)^{1/q}
\lsim (-\log |x|)^{1/q'} \,
\| \, f \, |B^{d/p}_{p,q}(\Rd)\|
\ee
as well as
\beq
&& \hspace*{-0.8cm}
\sum_{j=r+n_1+1}^{\infty} \sum_{k=2^{j-r}-n_0}^{2^{j-r}+n_0} |s_{j,k}|
 \lsim  \,
\sum_{j=r+n_1 + 1}^{\infty} 2^{-(j-r)(d-1)/p}\,
\Big(\sum_{k \ge 2^{j-r}-n_0} (1+|k|)^{d-1}\, |s_{j,k}|^p \Big)^{1/p}
\nonumber
\\
& \lsim &
\Big( \sum_{j=r+n_1 + 1}^{\infty}
\Big(\sum_{k \ge 2^{j-r}-n_0} (1+|k|)^{d-1}\, |s_{j,k}|^p \Big)^{q/p} \Big)^{1/q}
\, \Big( \sum_{j=r+n_1 + 1}^{\infty} 2^{-(j-r)q'(d-1)/p}\, \Big)^{1/q'}
\nonumber
\\
\label{f-03}
& \lsim & \| \, f \, |B^{d/p}_{p,q}(\Rd)\|
\, .
\eeq
The inequalities (\ref{f-02}) and (\ref{f-03}) yield (\ref{eq-72l}).
\\
{\em Step 2.} Proof of (ii). Let $1 <p< p_0 < \infty$.
The Jawerth-Franke embedding $F^{d/p}_{p,1} (\Rd) \hookrightarrow
B^{d/p_0}_{p_0,p} (\Rd)$, see \cite{Ja0} or \cite{SiTr}, combined with (\ref{eq-72l}) proves (\ref{eq-72li}).
\epr



\bigskip
\vbox{\noindent Leszek Skrzypczak\\
Adam Mikiewicz University Poznan\\
Faculty of Mathematics and Computer Science\\
Ul. Umultowska 87\\
61-614 Pozna{\'n}, Poland\\
e--mail: {\tt lskrzyp@amu.edu.pl}}

\bigskip
\smallskip
\vbox{\noindent  Winfried Sickel\\
Friedrich-Schiller-Universit\"at Jena\\
Mathematisches Institut\\
Ernst-Abbe-Platz 2\\
07743 Jena, Germany\\
e-mail:
{\tt winfried.sickel@uni-jena.de}}

\bigskip
\smallskip
\vbox{\noindent  Jan Vybiral\\
Austrian Academy of Sciences\\
Johann Radon Institute \\
for Computational and  Applied Mathematics \\
Altenberger Str. 69\\
A-4040 Linz, Austria\\
e-mail: {\tt jan.vybiral@oeaw.ac.at}}

\end{document}